\documentclass[11pt,reqno]{amsart}
\usepackage{amsthm,amsfonts,amssymb,euscript}
\usepackage{empheq}
\usepackage{esint}

\usepackage{mathrsfs}
\usepackage{color}

\newcommand{\bea}{\begin{eqnarray}}
\newcommand{\eea}{\end{eqnarray}}
\def\beaa{\begin{eqnarray*}}
\def\eeaa{\end{eqnarray*}}
\def\ba{\begin{array}}
\def\ea{\end{array}}
\def\be#1{\begin{equation} \label{#1}}
\def \eeq{\end{equation}}

\def\beq{\begin{equation}}

\def\be{{\beta}}

\def\MM{{\mathcal M}}

\def\FF{{\mathcal F}}

\def\RR{{\mathcal R}}

\def\R{{\mathbb{R}}}

\def\N{{\bf N}}

\def\Z{{\mathbb{Z}}}

\DeclareMathOperator*{\slim}{s-lim}

\newtheorem{theorem}{Theorem}[section]
\newtheorem{lemma}[theorem]{Lemma}
\newtheorem{proposition}[theorem]{Proposition}
\newtheorem{corollary}[theorem]{Corollary}

\newtheorem{remark}[theorem]{Remark}

\numberwithin{equation}{section}

\title[Resonances in the presence of a potential]{Nonlinear resonances with a potential: Multilinear estimates and an application to NLS}

\author[P. Germain]{Pierre Germain}

\address
   {Courant Institute of Mathematical Sciences, 251 Mercer Street, New York 10012
NY, USA}
\email{pgermain@cims.nyu.edu}
\thanks{P. Germain is partially supported by NSF grant DMS-1101269, a start-up grant from the Courant Institute, and a Sloan fellowship}

\author[Z. Hani]{Zaher Hani}
\email{hani@cims.nyu.edu}
\thanks{Z. Hani is partially supported by a Simons Postdoctoral Fellowship.}
\author[S. Walsh]{Samuel Walsh}
\email{walsh@cims.nyu.edu}


\usepackage{hyperref} 

\begin{document}

\maketitle

\begin{abstract}
This paper considers the question of global in time existence and asymptotic behavior of small-data solutions of nonlinear dispersive equations with a real potential $V$. The main concern is treating nonlinearities whose degree is low enough as to preclude the simple use of classical energy methods and decay estimates. In their place, we present a systematic approach that adapts the space-time resonance method to the non-Euclidean setting using the spectral theory of the Schr\"odinger operator $-\Delta+V$.  We start by developing tools of independent interest, namely multilinear analysis (Coifman--Meyer type theorems) in the framework of the corresponding distorted Fourier transform. As a first application, this is then used to prove global existence and scattering for a quadratic Schr\"odinger equation.
\end{abstract}

\tableofcontents

\section{Introduction} \label{intro:sec}

\subsection{The homogeneous setting}

Consider a prototypical homogeneous nonlinear dispersive PDE of the form
\begin{equation}
\label{eqdebaseE}
\left\{ \begin{array}{l}
i \partial_t u + L(D) u = \mathcal N(u), \qquad D:= -i\nabla \\
u(t=0) = u_0 \end{array} \right.
\end{equation}
where the equation is set in $\mathbb{R}^d$, $u$ is a complex-valued scalar or vector, $L(D)$ is a real Fourier multiplier, and $\mathcal N$ is superlinear of order $p$ in $u$ (e.g., $|u|^p$).  
Taking $L = \Delta$, for example, gives the nonlinear Schr\"odinger equation (NLS); $L = |D|$ corresponds to the wave equation; and $L = \langle D \rangle$ is the Klein--Gordon equation.  We call \eqref{eqdebaseE} \emph{homogeneous} due to the homogeneity of the linear part $i \partial_t + L(D)$.  We might also refer to it, variously, as the \emph{flat}, \emph{unperturbed}, or \emph{Euclidean} problem.  

As far as the global theory of this equation goes, the most basic question to ask is: {With $u_0$ sufficiently small, smooth and localized, does there exist a global solution? 
If so, does it scatter?} For high degree nonlinearities, the answer is yes, almost 
regardless of the structures of $L$ or $\mathcal N$: dispersive or Strichartz estimates for the linear part $L$ are enough to close the argument.  In contrast, classical energy and decay estimates fail for low degree nonlinearities ($p$ less than the Strauss exponent), and the structures of $L$ and $\mathcal N$, or more precisely resonances, start playing a decisive role. In fact, global existence is not assured, even for small smooth data (cf., e.g., \cite{john1979blow}). 

As it turns out, many physically interesting problems (in general relativity, plasma physics, water waves, etc.) fall into the second category. Consequently, this regime has been studied intensively over the past forty years. Most notably, two main methods were devised in the 1980s to deal with resonances: the normal forms approach of Shatah \cite{shatah1985normal}, and the vector fields approach of Klainerman \cite{Klainerman}.  Recently, the first author, Masmoudi and Shatah introduced a new method based on the concept of \emph{space-time resonances}  \cite{germain2011resonance,germain2012global,gustafson2009scattering}
which brings together (and goes beyond in some cases) the normal forms and vector fields methods. 

\subsection{State of the art in the inhomogeneous setting}

In this paper, we seek to give a systematic treatment to perturbations of \eqref{eqdebaseE} by a time-independent potential, namely equations of the form\begin{equation}
\label{eqdebase}
\left\{ \begin{array}{l}
i \partial_t u + L(\sqrt{-\Delta +V}) u = \mathcal N(u) \\
u(t=0) = u_0. \end{array} \right.
\end{equation}
Because of the appearance of the operator $\sqrt{-\Delta + V}$, we refer to this type of problem as \emph{inhomogeneous}, \emph{distorted}, \emph{perturbed}, or \emph{non-Euclidean}.  

There are many reasons, both mathematical and physical, to study equations of this form.   First, in a general sense, it is important to understand how properties of solutions to the homogeneous problem react to various perturbations, such as external forcing which is commonly represented by introducing a potential.  Second, even working with a homogeneous model, establishing the stability of a bound state or traveling wave often requires studying an inhomogeneous problem.  This is a simple consequence of the fact that the linearized operator around a ground state automatically involves a potential.  Third, many important physical systems come with a potential inhomogeneity intrinsic to the model. For example, one encounters these issues in studying the stability of the equilibrium state for the water waves equation, since the resulting linearized problem will be homogeneous only if the bottom is perfectly flat \cite{germain2012global}.  In general relativity, as well, one often encounters nonlinear wave equations with a potential (e.g., when studying Schwarzschild metric in Regge--Wheeler coordinates \cite{blue2008decay}.)

Regarding \eqref{eqdebase}, again the fundamental questions to ask are: {For $u_0$ sufficiently smooth, localized, and small, does there exist a global solution? If so, what is its asymptotic behavior?}  As in the homogeneous setting, if the nonlinearity is of high enough degree, global existence and scattering follow easily.  The only additional proviso is that the operator $i \partial_t + L(\sqrt{-\Delta+V})$ must enjoy similar dispersive estimates as the unperturbed operator $i \partial_t + L(D)$.  Sufficient conditions for this have been studied by many authors (cf., e.g., ~\cite{rodnianski2004decay,schlag2007dispersive}). 

On the other hand, comparatively little is known for low degree nonlinearities.  A first line of research in this direction was pursued by Keel, Smith, and Sogge~\cite{keel2002global}, Metcalfe, Nakamura, and Sogge~\cite{metcalfe2005global},  Metcalfe, and Sogge~\cite{metcalfe2007global}. 
These authors considered the nonlinear wave equation in dimension 3 on an exterior domain, and assuming that the nonlinearity satisfies the null condition.  Though this is not explicitly in the framework of \eqref{eqdebase}, note that one can think of exterior domains as a limit of potential perturbations. In these works, small-data global existence was proved by following either of two distinct approaches: the conformal method of Christodoulou, or the vector fields method of Klainerman. To adapt the latter approach, the authors had to (smoothly) cutoff the standard full-space vector fields on a compact set containing the forbidden domain and prove weighted space-time energy estimates.  Very recently, S. Yang~\cite{yang2012global} proved global existence for a nonlinear wave equation involving linear perturbation terms like $L^\mu \partial_\mu u$.  This was done under either smallness or mild (time-)decay assumptions on $L^\mu$ by using the vector field method and a new scheme developed by Dafermos and Rodnianski \cite{dafermos2010new}. We remark that none of these authors address the question of asymptotic behavior of solutions, which is something that we are able to do here.  

Turning to the Schr\"odinger equation, Cuccagna, Georgiev, and Visciglia \cite{cuccagna2012decay} considered the nonlinear Schr\"odinger equation in dimension 1, with an added potential $V$. They were able to prove global existence and describe the asymptotic behavior by adapting the (one-dimensional) approach in \cite{mckean1991nonlinear}.
All the works which we have discussed so far examine (as we will be doing for the rest of this manuscript) the situation where the linear operator spectrum is continuous, but one cannot get global existence simply by means of the 
dispersive estimates, since the decay rate given by the linear group is too weak and the power of the nonlinearity is not high enough to compensate for it. An interesting related problem is the situation where the strength of 
the dispersion is sufficient, but the potential exhibits linear, or nonlinear, bounds states. We refer for instance to~\cite{SW,tsai2001asymptotic,gustafson2004asymptotic,BC} 
for more on this direction of research.

\subsection{Space-time resonances in the inhomogeneous setting} \label{intro:spacetimeres}
The central objective of this work is to develop a space-time resonance method in the inhomogeneous regime.  
At a conceptual level, this would yield a more general understanding of resonant structures in the distorted setting. In a more practical sense, this will allow us to address both the global existence \emph{and} the asymptotic behavior for a number of problems of the form  \eqref{eqdebase}. 

To make matters more concrete, let us consider the equation:
\begin{equation}\label{QNLS}
\left\{ \begin{aligned}
& i \partial_t u - \Delta u + V u = \bar u^2\\
& u(t=0)=\; u_0.
\end{aligned} \right.
\end{equation}
set on $\mathbb{R}^3$. 
 Note that a quadratic nonlinearity qualifies as low degree in $\R^3$ as the Strauss exponent there is exactly 2. In particular, this indicates that resonances need to be taken into account. We shall argue here in an informal fashion. We start by introducing the generalized eigenvectors $e(x,\xi)$ which diagonalize the self-adjoint operator $H = - \Delta + V$. They are given by
$$
H e(x,\xi) = |\xi|^2 e(x,\xi) \qquad \mbox{and} \qquad e(x,\xi) \sim e^{ix\cdot \xi} \quad \mbox{if $|x| \rightarrow \infty$}.
$$
The properties and necessary conditions for the existence of these generalized eigenfunctions are reviewed in Section \ref{prelim:schrodingersec}. They naturally give rise to a distorted Fourier transform, which we denote by $\mathcal{F}^\sharp$ or $\cdot^\sharp$:
$$
(\mathcal{F}^\sharp g)(\xi) = g^\sharp(\xi) := \frac{1}{(2\pi)^{3/2}}\int g(x) \overline{e(x,\xi)}\,dx.
$$
(cf. Theorem \ref{prelim:Fsharptheorem}). As in the homogeneous setting, the idea is to write Duhamel's formula for the profile $f := e^{-itH} u$ in the distorted Fourier space. This gives after a short computation
\begin{equation}
\label{pinguin}
f^\sharp(t,\xi) =  u_0^\sharp(\xi) -i \int_0^t \iint e^{-is\phi(\xi,\eta,\zeta)} \overline{f^\sharp}(s,\eta) \overline{f^\sharp}(s,\zeta) \overline{\mathcal{M}}(\xi,\eta,\zeta) \,d\eta\,d\zeta \,ds
\end{equation}
where
$$
\mathcal{M}(\xi,\eta,\zeta) = \frac{1}{(2\pi)^{9/2}} \int {e(x,\xi)} e(x,\eta) e(x,\zeta) \,dx \quad \mbox{and} \quad \phi(\xi,\eta,\zeta) = |\xi|^2 + |\eta|^2 + |\zeta|^2.
$$

Observe that the flat case, $\mathcal{M}(\xi,\eta,\zeta) = \delta_{\R^3} (\xi+\eta+\zeta)$, and thus the above integral reads
\begin{equation}\label{flatDuhamel}
\int_0^t \int e^{is\phi(\xi,\eta,-\xi-\eta)} \overline{\widehat{f}(s,\eta) \widehat{f}(s,-\xi-\eta)} \,d\eta \,ds.
\end{equation}
This simpler formula is of course a manifestation of the convolution identity $\widehat{fg} = \widehat{f} * \widehat{g}$, which does not hold anymore for the distorted Fourier transform. As we argue below, this is the source of tremendous difficulties in the analysis if $V\neq 0$ as it makes manipulations on the (distorted) Fourier side much more delicate. But it is also responsible for an interesting new phenomenon: all frequencies $\xi$, $\eta$ and $\zeta$ interact, not only the ones that add up to zero. However, we expect $\mathcal{M}$ to be nicer away from the hyperplane $\{ \xi + \eta + \zeta = 0 \}$; on this
set, a singularity occurs, which seems to be in general more complicated than the $\delta$-function we see in the flat case.

We now give a brief outline of the space-time resonance method (see~\cite{germain2011resonance} for a more complete discusssion in the case $V=0$).
The idea is to consider the integral on the right-hand side of~(\ref{pinguin}) as an oscillatory integral. 
Proving global existence and scattering essentially amounts to ensuring that this integral remains bounded 
(in a sense that we keep vague for the moment) as $t$ goes to infinity. In the flat setting,  obstructions to this behavior come only from the stationary points (in the $(s, \eta)$ integral) of the phase $\phi$ in \eqref{flatDuhamel}. This leads us to define the time and space resonant sets as the sets where $\phi$ and $\partial_\eta \phi$ vanish, respectively. In contrast, in the distorted setting \eqref{pinguin}, there are three integration variables $(s, \eta, \zeta)$ and a distribution $\MM$ that is singular on the set $\{\xi+\eta+\zeta=0\}$, a fact which poses a challenge even in defining the space resonance set in the first place. More importantly, it is crucial to point out that there are specific directions (or vector fields) along which one can differentiate $\mathcal{M}$ without increasing the severity of its singularity (as quantified by its order as a distribution, say\footnote{Recall that a distribution $\Gamma$ is said to be of order $m$ if $\langle \Gamma, \phi\rangle$ can be bounded in terms $\|\phi\|_{C^m}$ for any test function $\phi$. Distributions of order 0, like $\MM$, are measures \cite{hormander2007analysis}.}). This phenomenon is even lurking in the flat case:  applying the vector fields $\partial_{\eta^\alpha}- \partial_{\zeta^\alpha}$  to $\delta(\xi+\eta+\zeta)$ gives a distribution of order 0, whereas applying $\partial_{\eta^\alpha}$ yields a distribution of order $1$. 

Understanding and coming to terms with these ``admissible" directions plays a central role both in defining the space resonant set below, and, more generally, in bounding expressions involving derivatives of the distribution $\mathcal M$. In Section \ref{derivatives of M section}, we show that, similar to the case $V=0$, the directions given by the vector fields $\partial_{\eta^\alpha} - \partial_{\zeta^\alpha}$ are, \emph{morally speaking}, admissible. More precisely, though, the analysis suggests the need for certain non-local modifications of those vector fields (cf. Section \ref{understanding M} and identity \eqref{derivative of M}), but let us gloss over this important fact for now.

With this information in hand, one can finally mimic the flat scenario and make the following definitions.
\begin{itemize}
\item The {\em time-resonant set} $\mathcal{T}$ consists of points that are stationary in $s$ (i.e. $\partial_s (s \phi)=0$): it is thus $\mathcal{T} = \{ \phi = 0 \}$. This corresponds to resonant interactions in the ODE sense.
\item The {\em space-resonant set} $\mathcal{S}$ consists of points that are stationary in $(\eta, \zeta)$ {when} we differentiate the phase function along the admissible directions. Therefore, $\mathcal{S} = \{ (\partial_\eta - \partial_\zeta) \phi = 0 \}$. This corresponds to wave packets traveling at the same group velocity and agrees with the previous definition when $V=0$.
\item Finally, the {\em space-time resonant set} $\mathcal{R}$ is the intersection of the space and time resonant sets: $\mathcal{R} = \mathcal T \cap \mathcal S$.
\end{itemize}

This classification of stationary points gives a clear intuitive understanding of the resonant interactions at hand, and this understanding helps to answer the next question: How can we obtain the required estimates?
\begin{itemize}
\item Away from the time resonant set, it is possible to integrate by parts in $s$ using the identity $\frac{1}{i\phi} \partial_s e^{is\phi} = e^{is\phi}$. 
This corresponds to a normal form transform, which effectively increases the degree of the nonlinearity.
\item Away from the space resonant set, it is possible to integrate by parts in $(\eta, \zeta)$ using the identity $\frac{1}{is  (\partial_{\eta^\alpha} - \partial_{\zeta^\alpha}) \phi} (\partial_{\eta^\alpha} - \partial_{\zeta^\alpha}) e^{is\phi} = e^{is\phi}$. Of course, the gain here comes from the factor of $s$ in the denominator which improves the convergence chances of the integral. Note that the use of an admissible vector field is essential to control the terms involving derivatives of $\mathcal{M}$ which are inadvertently produced when we integrate by parts. For the same reason, we require bounds on derivatives of $f^\sharp$, or, equivalently, weighted-norm estimates on $f$.
\item Finally, there remains the space-time resonant set, for which none of the above strategies applies. Here the analysis becomes more problem-dependent. One possibility would be to simply split the above integral by cutting off a shrinking neighborhood of the space-time resonant set. To estimate the piece thus removed, one can use its shrinking size to gain integrability; as for the rest, the above integrations by parts apply.
\end{itemize}

For equation \eqref{QNLS}, a simple computation shows that $\mathcal{R} = \mathcal{T} = \{ (0,0, 0) \}$,
and thus space resonances can be ignored for this problem. Nonetheless, even in such a simplified context, a quantitative understanding of the ``admissible derivatives" of the distribution $\MM$ is indispensable if one needs to prove weighted estimates on the solutions. We will treat equations with $\mathcal S \not\subset \mathcal T$ in a future work using the harmonic analytic framework developed here.

\subsection{Assumptions on the potential $V$} Before stating our main results, let us now describe the assumptions we shall impose on $V$.
%
\subsubsection{A qualitatative description}
\begin{itemize}
\item {\bf H1} \emph{Existence of distorted Fourier Analysis} \cite{agmon1975spectral}
\item {\bf H2} \emph{Absence of discrete spectrum for $-\Delta + V$} 
\item {\bf H3} \emph{$L^p$ boundedness of the wave operator
$\Omega:=\lim_{t\to -\infty} e^{itH}e^{-it\Delta}$} \cite{beceanu2011waveop,yajima1993waveop}. 
\end{itemize}

In total, {\bf H1}, {\bf H2}, {\bf H3} amount to some regularity and decay requirements on $V$.  For further discussion of {\bf H1} and {\bf H2}, including  explicit sufficient conditions, 
see Section \ref{prelim:schrodingersec} (namely Theorem \ref{prelim:Fsharptheorem} and Remark \ref{prelim:H12remark}).
Assumption {\bf H3} has been proved to hold for $1\leq p \leq \infty$ in all dimensions $d\geq 2$ for so-called \emph{generic-type} potentials $V$ (i.e. with no null resonances) with  sufficient decay at infinity. 
The literature on this topic is quite lengthy, and we give a somewhat detailed account in Section \ref{prelim:waveopsec}.

At various points in our analysis, we will need to take a deeper look at the wave operator $\Omega$.  In those cases, we will assume one or both of the following:
\begin{itemize}
\item {\bf H3*\phantom{*}} \emph{$W^{1,p}$ boundedness of the wave operator $\Omega$}~\cite{beceanu2011waveop,yajima1993waveop}
\item {\bf H3**} \emph{An explicit structure theory for $\Omega$}~\cite{beceanu2011waveop,yajima1993waveop} \end{itemize}
Each of these is closely related to {\bf H3}; in particular most proofs of {\bf H3} and {\bf H3*} rely precisely on an explicit structure theory as assumed in {\bf H3**}.  
The reader is again directed to Section \ref{prelim:waveopsec} for a more detailed presentation. 
\subsubsection{Quantitative description on $\R^3$}
On $\R^3$ (where we set our PDE application), Beceanu provides in \cite{beceanu2011waveop} the most refined result guaranteeing assumptions {\bf H3, H3*, H3**} with minimal decay conditions on $V$ (e.g., $V \in \langle x \rangle^{-1/2-\epsilon}L^2$ for {\bf H3, H3**}). However, at one juncture in our analysis (Section \ref{commutator subsection}) we find it easier to use Yajima's slightly more explicit structure theory for $\Omega$ in \cite{yajima1993waveop}.  This comes at the cost of imposing more decay conditions on the potential in Theorem \ref{intro:gwptheorem} below (cf. Remark \ref{prelim:H3remark}). A very simple condition (but very far from optimal) which ensures that all these hypotheses are met on $\R^3$ is the following
\begin{itemize}
\item The operator $-\Delta + V$ does not have discrete spectrum or null-resonances.
\item There exists $\epsilon>0$ such that 
\[ \displaystyle |V(x)|  \lesssim \langle x \rangle^{-\frac{11}{2}-\epsilon}.\]
\end{itemize}

\subsection{Obtained results: Multilinear distorted Fourier analysis} 
Our first set of theorems lays the harmonic analytic groundwork for the space-time resonance method.  We hasten to point out, though, that they are of interest in their own right.  While it may seem that transitioning from the flat to the distorted regime is straightforward --- all one must do is exchange $\mathcal{F}$ for $\mathcal{F}^\sharp$ --- in practice, it requires some new ideas, and a few highly technical arguments.  In that sense, these theorems  should be thought of as one of the principal contributions of this paper.   

The oscillatory integrals encountered in Section \ref{intro:spacetimeres} lead one to study multilinear operators of the type
\begin{equation}\label{bilinear00}
T(f,g):={\mathcal{F}_\xi^\sharp}^{-1} \iint_{\R^d \times \R^d} m(\xi, \eta, \zeta) f^\sharp(\eta) g^\sharp(\zeta) \mathcal M(\xi,\eta,\zeta) d\zeta d\eta.
\end{equation}
These are the generalizations of the classical pseudo-product operators of Coifman and Meyer~\cite{coifman1978dela} to the inhomogeneous setting. Typically, one would like to prove H\"older-type estimates, i.e. that $T$ is bounded from $L^q\times L^p \to L^{r'}$ with $\frac{1}{r'}=\frac{1}{q}+\frac{1}{p}$ which hold trivially if $m=1$. A natural condition on $m$ mimics the Coifman--Meyer class of symbols by requiring
\begin{equation}\label{CMC0}
|\partial_{\xi_1}^{\alpha_1}\ldots \partial_{\xi_n}^{\alpha_n}m(\xi_1, \ldots, \xi_n)|\leq C_{\alpha_1, \ldots, \alpha_n} (|\xi_1|+\ldots+|\xi_n|)^{-(|\alpha_1|+\ldots+|\alpha_n|)}.
\end{equation}
Any $m$ satisfying the above homogeneous bounds is said to be a Coifman--Meyer symbol (in $n$ variables).
  
Our first result is the distorted analog of the celebrated Coifman--Meyer theorem~\cite{coifman1978dela}.
\begin{theorem}\label{intro:CMtheorem}
Let $V \in L^{d/2}(\R^d)$ ($d\geq 2$) be a potential satisfying $\bf{H1}$, $\bf{H2}$, and $\bf{H3}$. Suppose that $m = m (\xi, \eta, \zeta)$ is either a Coifman--Meyer symbol in three variables as in \eqref{CMC0}, or $m(\xi, \eta, \zeta)=m_0(\eta, \zeta)$ where $m_0$ is a Coifman--Meyer symbol in two variables.
\begin{enumerate}
\item[(i)] For any $p, q, r \in (1, \infty)$ satisfying $\frac{1}{p}+\frac{1}{q}+\frac{1}{r}=1$, we have
\begin{equation}\label{trilinearest00}
\|T(f,g)\|_{L^{r'}(\R^d)} \lesssim_{m, V} \|f\|_{L^q(\R^d)} \|g\|_{L^p(\R^d)}
\end{equation}
provided that the Riesz transform $\mathfrak{R}=\nabla (-\Delta+V)^{-1/2}$ is bounded on $L^p, L^q,$ and $L^r$.
\item[(ii)] Suppose instead that $V$ satisfies assumption ${\bf H3^*}$. If $q,p,r, \tilde q, \tilde p \in (1,\infty)$ satisfy $\frac{1}{r^\prime}=\frac{1}{p}+\frac{1}{q}=\frac{1}{\tilde p}+\frac{1}{\tilde q}-\epsilon$ for some $\epsilon>0$,
\begin{equation}\label{trilinearest01}
\|T(f,g)\|_{L^{r'}(\R^d)} \lesssim_{m, V, \epsilon} \|f\|_{L^q(\R^d)} \|g\|_{L^p(\R^d)}+\|f\|_{L^{\tilde q}(\R^d)} \|g\|_{L^{\tilde p}(\R^d)},
\end{equation}
whenever $f\in L^q\cap L^{\tilde q}$ and $g\in L^p\cap L^{\tilde p}$.
\end{enumerate}
\end{theorem}

The boundedness of the Riesz transform associated to a potential $V$ on $L^p(\R^d)$ is clearly equivalent to the embedding of the distorted Sobolev space 
$\dot W_{\sharp}^{1, p}:=(-\Delta+V)^{-1/2} L^p \hookrightarrow \dot W^{1,p}$, a question that has been an active subject of research for some time. 
We touch on the issue briefly in Section \ref{prelim:riesz}.  However, for our PDE applications, the imperfect estimate \eqref{trilinearest01} suffices, and so we do not need to make such strong assumptions on $V$.  
The proof of Theorem \ref{intro:CMtheorem} starts by following the same strategy as that of the classical Coifman--Meyer theorem. 
This gives the needed bounds on the  diagonal interactions (where the highest two frequencies are comparable) once one proves maximal and square function estimates for the distorted setting 
(cf. Lemma \ref{square and maximal lemma}). In contrast to the flat setting, however, non-diagonal interactions are present and very problematic when there is a potential. 
To deal with them, we use a correlation identity involving three generalized eigenfunctions (the same idea was used in \cite{hani2012global} for similar reasons). 
What results is a trilinear term involving the Riesz transform; this is precisely why we must require that $\mathfrak{R}$ is bounded on $L^p$, or else sacrifice some $\epsilon$ of integrability. 
Related work, dealing with a more restricted class of operators, namely para-products instead of pseudo-products, can be found in~\cite{bernicot20121}.

Obtaining H\"older-type estimates for $T(f,g)$ turns out to be merely part of the story. 
Indeed, to successfully prosecute the space-time resonance method, we will need weighted estimates for $T(f,g)$, or, equivalently, $L^q\times L^p \to L^{r'}$ estimates 
for \eqref{bilinear00} with a $\partial_{\xi}$ derivative falling on the integral. This leads inexorably to the issue of admissible directions discussed in Section \ref{intro:spacetimeres}.  
Ultimately, we are able to obtain the needed estimates bounding $\|\langle x \rangle T(f,g)\|_{L^{r'}}$ in terms of $\|\langle x \rangle f\|_{L^q}$ and $\|g\|_{L^p}$ 
for symbols that are slightly more regular than Coifman--Meyer ones and for potentials satisfying hypothesis {\bf H3**}. This is the content of Theorem \ref{colibri} which we do not transcribe here for brevity.

\subsection{Obtained results: Global existence for a quadratic Schr\"odinger equation}
With the above tools in hand, we are able to prove the following.  
\begin{theorem} \label{intro:gwptheorem} Let $V$ be a potential satisfying assumptions {\bf H1, H2, H3*}, and {\bf H3**} and consider the quadratic nonlinear Schr\"odinger equation in $\mathbb{R}^3$ \eqref{QNLS}. Let $X$ be the space defined by the norm
\[ \| f \|_X := \| f \|_{L_t^\infty H_x^1} + \| \partial_\xi f^\sharp \|_{L_t^\infty L_\xi^2}.\]
There exists $\epsilon_0 > 0$ such that, for any initial data $u_0$ with $\| u_0 \|_{H^1} + \|\partial_\xi u_0^\sharp \|_{L^2} < \epsilon_0$, 
there exists a solution $u \in X$ to \eqref{QNLS} defined for all time $t \geq 0$. In particular,
$$\|u(t)\|_{L^p} \lesssim \frac{1}{t^{\frac{3}{2}(1-\frac{2}{p})}}\qquad \text{for all } 2\leq p \leq 6.$$
Finally, $u$ scatters in $L^2$, namely $e^{it\Delta} u(t)$ has a limit in $L^2$ as $t \rightarrow \infty$.
\end{theorem}

The proof of this result is presented in Section \ref{ubar2} and is meant as a first application of the ``inhomogeneous" space-time resonance strategy outlined in Section \ref{intro:spacetimeres}.  
The energy space $X$ is natural in light of the dispersive estimates for the linear propagator $e^{itH}$ (cf. Proposition \ref{prelim:weightedispersiveprop}).  
In particular, we note that in the flat setting, global existence and scattering was recently proved in a similar space by Laillet \cite{laillet2011thesis}, who simplified the argument based on the space-time resonance approach originally given in \cite{germain2009global}.

Finally, we remark that our choice of NLS with a $\bar{u}^2$ nonlinearity is meant to provide a somewhat generic example of how to proceed when the unperturbed homogeneous problem can be treated using the normal forms approach. 
Many other homogeneous equations fall into this category, for example Klein--Gordon \cite{shatah1985normal,Kosecki1992}, water waves \cite{wu2011global,germain2012global}, Klein--Gordon--Zakharov \cite{Ozawa1995}, and 
Euler--Poisson for the electrons~\cite{Guo1998,germain2012non} or the ions~\cite{guo2011global}. 
To see the full strength of the space-time resonance method, though, one must consider problems where the normal forms approach on its own is not enough.  This we do in a forthcoming paper.
 
\subsection*{Plan of the article}
In {Section \ref{toolssection}}, we present several known results on the spectral theory of Schr\"odinger operators $-\Delta + V$ which are needed in the rest of the paper. In {Section 3}, we develop the theory of multilinear distorted Fourier multipliers, that is to say the analog of pseudo-products when the distorted Fourier transform is substituted to the Fourier transform. This is applied in {Section 4} to prove global existence and scattering for a quadratic Schr\"odinger equation in $\mathbb{R}^3$.

\subsection*{Notation} We write $A\lesssim B$ to signify that there is a constant $C>0$ such that $A\le C B$. We also write $A\sim B$ when $A\lesssim B\lesssim A$. If the constant $C$ involved has some explicit dependency, we emphasize it by a subscript. Thus $A\lesssim_uB$ means that $A\le C(u)B$ for some constant $C(u)$ depending on $u$. In some instances, we use the notation $A \ll B$ to signify that the implicit constant $C$ is large. Also, we use the Japanese bracket convention where $\langle x \rangle := \sqrt{1 + |x|^2}$.

\section{Tools from linear distorted Fourier Analysis} \label{toolssection}

In this section, we review a number of important topics from the theory of distorted Fourier analysis that will form the basis for our work in the remainder of the paper.   

\subsection{Schr\"odinger operators} \label{prelim:schrodingersec}  For a given potential $V : \mathbb{R}^d \to \mathbb{R}$,  consider the associated Schr\"odinger operator  $H := -\Delta + V$.  Of particular interest is the situation where $H$ is  a  perturbation of $H_0 := -\Delta$.    For instance, when  $V \in L^2(\mathbb{R}^d)$,  $H$ can be realized as a self-adjoint  operator on $L^2(\mathbb{R}^d)$ with domain $D(H) = H^2(\mathbb{R}^d)$.  We may then ask: For what $V$ do the spectral properties of $H$ resemble those of $H_0$?  A natural starting point is to impose a compactness condition on the multiplication operator associated with $V$.  With that in mind, we say that  $V$ is  \emph{short-range} (or, \emph{of class SR}) provided that 
\[ u \in H_x^2(\mathbb{R}^d) \mapsto (1+|x|)^{1+\epsilon} V u \in L_x^2(\mathbb{R}^d) \qquad \textrm{is a compact operator,}  \]
 for some $\epsilon > 0$.  It was shown by Agmon (cf. \cite{agmon1975spectral}) that, for $V$ of class SR, $\sigma(H) = \{ \lambda_j \}_{j \in J} \cup [0,\infty)$;  
the continuous spectrum being $[0,\infty)$, and the discrete spectrum consisting of a countable set of real eigenvalues $\{ \lambda_j \}$, each of finite multiplicity.
 Furthermore, we have the orthogonal decomposition \[ L^2(\mathbb{R}^d) = L^2_{\textrm{ac}}(\mathbb{R}^d) \oplus L^2_{\textrm{p}}(\mathbb{R}^d),\] where $L^2_{\textrm{p}}(\mathbb{R}^d)$ is the span of the eigenfunctions corresponding to the eigenvalues $\{ \lambda_j \}$, and $L^2_\textrm{ac}$ is the absolutely continuous subspace for $H$.  
 
 It should be noted that numerous versions of this result predate the work of Agmon, only for different classes of potential.   Indeed, \cite{agmon1975spectral} is a particularly significant milestone in a long and continuing effort to discern the optimal conditions on $V$ under which these and other spectral properties of $H$ hold (cf. the references contained in \cite{agmon1975spectral} for a summary of earlier works, and, e.g., \cite{ionescu2006agmon} for an important recent improvement).
 
 In order to identify and study the resonances in PDEs like \eqref{eqdebaseE} or \eqref{QNLS}, it is  best formulate them in frequency space. For this, we need a ``well-behaved'' eigenfunction expansion that diagonalizes $H$ (on $L_{\textrm{ac}}^2$). Of course, this is in turn predicated on the existence of a natural class of generalized eigenfunctions that serve as analogs of the plane waves $e^{ix.\xi}$ for $H_0:=-\Delta$.    For each $\xi \in \mathbb{R}^d\setminus\{0\}$, we know that $|\xi|^2$ is in the continuous spectrum of $H$;  the associated eigenfunction is the \emph{distorted plane wave} $e(\cdot; \xi)$ defined as the solution of 
\beq H e(\cdot; \xi) = |\xi|^2 e(\cdot; \xi), \label{prelim:eigenvalueprob} \eeq
with the asymptotic condition 
\[ v(x; \xi) := e(x; \xi) - e^{ix\cdot \xi} = O(|x|^{-1}) \qquad \textrm{as } |x| \to \infty,\]
 and the Sommerfeld radiation condition 
\[ \lim_{r \to \infty} r ( \partial_r -i|\xi|) v = 0.\]
This can be expressed in a more convenient way via the resolvent:  for $z \in \mathbb{C}\setminus \sigma(H)$, define $R_V(z) := (H-z)^{-1}$ and consider 
\[ R_V^{-}(z) := \lim_{\epsilon \to 0+}  R_V(z-i\epsilon).\]
The limit here is taken in the uniform operator norm topology on $\mathcal{L}(\langle x \rangle^{-s} L^2, \langle x \rangle^{s} H^{2})$ for $s>\frac{1}{2}$; its convergence is the so-called \emph{limiting absorption principle} (cf. \cite[Theorem 4.2]{agmon1975spectral}).  The eigenfunction problem \eqref{prelim:eigenvalueprob} can then be recast as the Lippman--Schwinger equation:
\beq e(\cdot; \xi) = e_\xi - R_V^-(|\xi|^2)V e_\xi, \qquad e_\xi(x) := e^{ix \cdot \xi}.\label{prelim:resolventeq} \eeq
It can be shown that there exists a unique solution to \eqref{prelim:resolventeq} for any $\xi \in \mathbb{R}^d \setminus \{0\}$  provided that $V = O(|x|^{-1-\epsilon})$ as $|x| \to \infty$, for some $\epsilon > 0$ (cf. \cite{agmon1975spectral}).  Under this assumption, the distorted plane waves are relatively smooth in $x$, but have very little regularity in $\xi$.  More precisely, for fixed $\xi \in \mathbb{R}^d \setminus \{0\}$,  
\beq e(\cdot; \xi) \in \langle x \rangle^{s} H_x^2, \textrm{ for any $s > (d+1)/2$,}\label{regeigenfunction} \eeq
however, the map $(x, \xi) \mapsto e(x; \xi)$ is merely measurable. One can improve this by requiring additional decay and regularity of $V$ (cf., e.g., \cite{ikebe1960eigen}).   
  
In view of the Fourier transform, we expect that the family $\{ e(\cdot; \xi)\}$ forms a basis for the absolutely continuous  subspace of $H$.  This is indeed true, as was first proved by Ikebe \cite{ikebe1960eigen} and later generalized by several authors.  For consistency of presentation, we give here the version due to Agmon (cf. \cite[Theorem 6.2]{agmon1975spectral}).  Before that, let us now impose assumption {\bf H2}, namely that $H$ has no discrete spectrum. However, we remark that many results in this paper (especially those in Section \ref{multestimates section}) can be directly generalized to potentials with discrete eigenvalues by simply projecting on the absolutely continuous subspace $L^2_{\textrm{ac}}$ throughout.

That said, the result is the following.
\begin{theorem}[Ikebe, Alsholm--Schmidt, Agmon] \label{prelim:Fsharptheorem} Consider the Schr\"odinger operator $H$ with potential $V$ satisfying {\bf H2} and 
\beq   ( 1+|x|)^{2(1+ \epsilon)} \int_{{B_1(x)}} |V(y)|^2 |y-x|^{-d+\theta} \, dy \in L_x^\infty(\mathbb{R}^d), ~ \textrm{for some $\epsilon > 0$, $0 < \theta< 4$.} \label{prelim:agmoncond} \eeq
Define the distorted Fourier transform $\mathcal{F}^\sharp$ by 
\beq (\mathcal{F}^\sharp f)(\xi) := f^\sharp(\xi) := \frac{1}{(2\pi)^{d/2}} \lim_{R \to \infty} \int_{B_R} \overline{e(x; \xi)} f(x) \, dx,\label{prelim:distortedFdef} \eeq
where $B_R$ is the ball or radius $R$ centered at the origin in $\mathbb{R}^d$. Then $\mathcal{F}^\sharp$ is an isometric isomorphism on $L^2(\mathbb{R}^d)$ with inverse formula
\beq f(x) = ({\mathcal{F}^\sharp}^{-1} f^\sharp)(x) := \frac{1}{(2\pi)^{d/2}} \lim_{R \to \infty} \int_{B_R} e(x; \xi) f^\sharp(\xi) \, d\xi . \label{prelim:distortedFinvdef} \eeq
Moreover, $\mathcal{F}^\sharp$ diagonalizes $H$ in the sense that, for all $f \in H^2(\mathbb{R}^d)$, 
\beq Hf = {\mathcal{F}^\sharp}^{-1}M \mathcal{F}^\sharp f,\label{prelim:Hdiag} \eeq
where $M$ is the multiplication operator $u \mapsto |x|^2 u $.
\end{theorem}

\begin{remark}[Sufficient conditions for assumptions {\bf H1} and {\bf H2}]\label{prelim:H12remark}
We are now, at last, able to give a precise meaning to assumption {\bf H1}: we say that {\bf H1} is satisfied provided that (i) the family of eigenfunctions $\{ e(\cdot, \xi) \}$ exists with the regularity stated in \eqref{regeigenfunction}, and (ii) the operator $\mathcal{F}^\sharp$ defined by \eqref{prelim:distortedFdef} exists and exhibits the properties described in Theorem \ref{prelim:Fsharptheorem}.     It follows that sufficient conditions for {\bf H1} are that $V$ satisfies {\bf H2}, \eqref{prelim:agmoncond}, and $V = O(|x|^{-1-\epsilon})$ as $|x| \to \infty$, for some $\epsilon > 0$.  If we remove assumption {\bf H2}, of course, we require only that $\mathcal{F}^\sharp$ be a unitary partial isometry with range $L^2_{\mathrm{ac}}(\mathbb{R}^d)$.  Note that, by imposing \eqref{prelim:agmoncond}, we rule out the existence of nonnegative eigenvalues.  In order for {\bf H2} to hold, we must require additionally that there are no negative eigenvalues, which is guaranteed, e.g., if the negative part of $V$ is not very large (for example if $d\geq 3$, Hardy's inequality implies that the condition $V\geq - (d-2)^2/4|x|^2$ is sufficient to rule out both non-positive eigenvalues and resonances at $0$ as defined in \eqref{prelim:norescond} below).
\end{remark}

\subsection{The wave operator $\Omega$} \label{prelim:waveopsec}  One consequence of Theorem \ref{prelim:Fsharptheorem} is that $H$ and $H_0$ are unitarily equivalent.  To see this, note that by \eqref{prelim:Hdiag}, 
\[ H = \Omega H_0 \Omega^*, \qquad\mbox{with} \quad \Omega := {\mathcal{F}^\sharp}^{-1} \mathcal{F}.\]
The operator $\Omega$ is called the {\em wave operator}. It can alternatively be defined by
\beq \Omega := \slim_{t \to -\infty} e^{itH} e^{-itH_0} 
 \label{prelim:defOmega}  \eeq
in the strong operator topology. Note that some authors denote this as $\Omega_+$, with $\Omega_-$ being the result of taking the limit $t \to +\infty$.  Under the assumption that $V$ is of class SR and {\bf H2} holds, these limits exist, and $\Omega$ is a unitary operator on $L^2$.  This fact is often referred to as \emph{asymptotic completeness}; for potentials satisfying \eqref{prelim:agmoncond}, it is originally due to Agmon (cf. \cite[Theorem 7.1]{agmon1975spectral}).

For us, the importance of $\Omega$ lies in the intertwining relations 
\beq e^{itH} = \Omega e^{itH_0} \Omega^*, \qquad \mathcal{F}^\sharp \Omega = \mathcal{F}. \label{prelim:intertwineidentities} \eeq 
In other words, $\Omega$ allows us to translate back and forth between the flat and distorted cases.  
Clearly, then, information about the structure and boundedness properties of $\Omega$ is extremely valuable.   The foundational work in this direction is due to Yajima, who first proved the $W^{k,p}$ boundedness of $\Omega$ and $\Omega^*$, under the assumption of sufficient smoothness and decay for the potential.  We paraphrase his results below.

\begin{theorem}[Yajima,  \cite{yajima1993waveop,yajima1995Wkp3, finco2006p,yajima2006Lp}] \label{prelim:yajimatheorem} Let $k \in \mathbb N$ and consider the Schr\"odinger operator $H$ with real potential $V : \mathbb{R}^d \to \mathbb{R}$ for $d \geq 3$.  Fix $p_0, k_0$ as follows:
\[ \left\{ \begin{array}{cl} p_0 = 2, ~ k_0 = 0 & \textrm{if } d = 3 \\
p_0 > d/2,~k_0 := \lfloor (d-1)/2 \rfloor &  \textrm{if } d \geq 4. \end{array} \right.\]
Assume that for some $\delta > (3d/2)+1$,
\beq   \langle x \rangle^\delta \| \partial^\alpha V \|_{L_y^{p_0}(|x-y| \leq 1)} \in L_x^\infty(\mathbb{R}^d), \qquad \textrm{for all $\alpha$ with $|\alpha| \leq k+k_0$.} \label{prelim:yajimacond} \eeq
Then $V$ is of class SR and so $\Omega$ and $\Omega^*$ are well defined as operators on $L^2(\mathbb{R}^d) \cap W^{k,p}(\mathbb{R}^d)$.  If we additionally assume that $V$ is of 
\beq \textrm{\textbf{Generic-type:} there is no $u \in \langle x \rangle^{\theta} L_x^2(\mathbb{R}^d)$ solving $H u = 0$, for any $\theta > 1/2$.} \label{prelim:norescond} \eeq
Then $\Omega$ and $\Omega^*$  may be extended to  bounded operators defined on $W^{k,p}(\mathbb{R}^d)$.
\end{theorem}
\begin{remark}  \label{prelim:yajimaremark}The assumptions here can be weakened somewhat when $V$ is small or nonnegative.  For example, in any dimension $d \geq 3$,  hypothesis \eqref{prelim:yajimacond} can be replaced by the following: for all $|\alpha| \leq k$, 
\beq \| \mathcal{F} \langle x \rangle^{\sigma} D^\alpha V\|_{L^{\frac{d-1}{d-2}}(\mathbb{R}^d)} \textrm{ is sufficiently small for some $\sigma > 2(d-2)/(d-1)$.}  \label{prelim:yajimasmallV} \eeq
See  \cite[Theorem 1.1, Theorem 1.2]{yajima1995Wkp3} and the surrounding remarks for further discussion.
\end{remark}

Condition \eqref{prelim:norescond} asks for the absence of resonance at zero. It is indeed necessary for $W^{k,p}$ boundedness, though other forms of boundedness can be salvaged without it (cf., e.g.,  \cite{yajima1995Wkp,yajima2006Lp}). We should also note that a similar result also holds on $\R^2$ \cite{jensen2002remark}.

Recently, Beceanu \cite{beceanu2011waveop} was able to extend a number of Yajima's results on $\R^3$, reaching scale invariant class of potentials. Namely, he was able to prove boundedness of the wave operator on 
$L^p(\mathbb{R}^3)$ for a class of potentials $B$ such that $\langle x \rangle^{-1/2-\epsilon} L^2 \subset B \subset L^{3/2,1}$. He also proved the $W^{1,p}(\mathbb{R}^3)$ boundedness of $\Omega$ under the assumption that:
\beq V \in \left\{ \begin{array}{ll} B & p < 3 \\
B \cap L^{3/2+\epsilon} & p = 3 \\
B \cap L^p & p > 3. \end{array}   \right. \label{prelim:beceanuW1pcond} \eeq
See \cite[Corollary 1.5]{beceanu2011waveop}.     

The above mentioned proofs are based on an asymptotic expansion of the wave operator as a Born series (obtained by repeated application of Duhamel's formula). 
The main difference between Yajima's and Beceanu's approaches is in bounding high order terms in this expansion. 
While Yajima resorts to a direct computation to estimate their contribution (see Lemma \ref{structurelemma}), Beceanu uses an abstract version of Wiener's theorem which allows him to work under much lower 
decay assumptions on the potential $V$.


Both Yajima and Beceanu prove their boundedness results by obtaining an ``explicit description" of the wave operator. 
We will need such a description in Section \ref{commutator subsection}, where we study commutators of $\Omega$ and position operators in $\R^3$. 
We elect to use Yajima's version, which is slightly more explicit than Beceanu's. 
In effect, this requires us to impose Yajima's more restrictive conditions on $V$ in Theorem \ref{commutator bound} (and hence Theorem \ref{colibri}, and Theorem \ref{intro:gwptheorem}). 
We believe that the weaker assumptions of Beceanu are sufficient, but showing this would entail considerable additional technical work.    

\begin{lemma}[Structure of the Wave operator \cite{yajima1995Wkp,beceanu2011waveop}]\label{structurelemma}
Under the assumptions of Theorem~\ref{prelim:yajimatheorem}, the adjoint $\Omega^*$ of $\Omega$ can be written as:
\begin{equation}\label{Omega expression0.1}
\Omega^* f= f -W_1 f +W_2 f - W_3 f +L f
\end{equation}
where $W_1, W_2, W_3 $ and $L$ are bounded operators on $L^p$ for all $1\leq p \leq \infty$ and have the following form:
\begin{align}
\label{W_j} W_j f=& \int_{[0,\infty)^{j-1}\times I \times \Sigma^j} \widetilde K_j(t_1, \ldots, t_{j-1}, \tau, \omega_1, \ldots, \omega_j) f(\overline x +\rho) dt_1 \ldots dt_{j-1}d\tau  d\omega_1 \ldots d\omega_j\\
\label{L} L f=& \int_{\R^3} L(x,y) f(y)\,  dy.
\end{align}
Here $\Sigma = S^2$ denotes the unit sphere in $\R^3$; $I=(-\infty,-2\omega_j(x+t_1\omega_1 +\ldots+t_{j-1}\omega_{j-1}))$ is the range of integration of the variable $\tau$; and we denote furthermore
$$\overline y:= y-2(\omega_j.y)\omega_j, \qquad \rho:= t_1 \overline \omega_1 +\ldots +t_{j-1}\overline{\omega_{j-1}}-\tau \omega_j.$$
The kernels $\widetilde{K}_1$, $\widetilde{K}_2$, $\widetilde{K}_3$,  and $L(x,y)$ are described more explicitly in Section \ref{commutator subsection} and they satisfy
\begin{equation}\label{K_j and L}
\begin{split}
\|\widetilde K_j\|_{L^1([0,\infty)^j, L^1(\Sigma^n))} \lesssim_V 1, \\
\sup_{y \in \R^3} \int_{\R^3} |L(x,y)| \, dx +\sup_{x\in \R^3} \int_{\R^3}|L(x,y)| \, dy \lesssim_V 1.
\end{split}
\end{equation}

\end{lemma}

\begin{remark}\label{boundedbypositive}
It is easy to see from the above explicit description (particularly estimate \eqref{K_j and L}) that the operator $\Omega$ is point-wise majored by a positive operator\footnote{Recall that a positive operator is one such 
that $T f \geq 0$ whenever $f\geq 0$.} that is bounded on $L^p(\R^3)$ for all $1\leq p \leq \infty$. This fact is also (even more) evident in Beceanu's explicit description of $\Omega$ in \cite{beceanu2011waveop}, 
so we shall assume it throughout our work.
\end{remark}

\begin{remark}[Sufficient condition for assumptions {\bf H3}, {\bf H3*}, and {\bf H3**}] \label{prelim:H3remark} 
We say that {\bf H3} holds provided $\Omega$ is bounded on $L^p(\mathbb{R}^d)$, for any $1\leq p \leq \infty$.  
Sufficient conditions for this are $V \in B$ (considered by Beceanu) and the generic-type assumption \eqref{prelim:norescond} when $d = 3$, or more generally \eqref{prelim:yajimacond} and \eqref{prelim:norescond} with $k = 0$, for $d \geq 4$.  
(See also Remark \ref{prelim:yajimaremark}).
{\bf H3*} holds provided $\Omega$ is bounded on $W^{1,p}(\mathbb{R}^d)$.  In $\mathbb{R}^3$, the optimal known conditions implying this are given by \eqref{prelim:beceanuW1pcond}; 
for other dimensions, the hypotheses of Theorem \ref{prelim:yajimatheorem} with $k = 1$ seem to be the weakest currently available.    
We say that {\bf H3**} is satisfied provided $\Omega$ can be written as a Born series of the form \eqref{Omega expression0.1}. 
This is satisfied whenever the hypotheses of Theorem \ref{prelim:yajimatheorem} holds (with $k=0$).  
\end{remark}

\subsection{Distorted Fourier multipliers} Suppose $V$ satisfies the hypotheses of Theorem \ref{prelim:Fsharptheorem}.  Then, for any function $m : \mathbb{R}^d \to \mathbb{C}$, we define the \emph{distorted Fourier multiplier} $m(D^\sharp)$ to be the operator
\beq m(D^\sharp) := {\mathcal{F}^\sharp}^{-1} m(\xi) \mathcal{F}^\sharp,\label{prelim:defmDsharp} \eeq
where $m(\xi)$ denotes the multiplication operator $u \mapsto m(\xi) u$.   This is an analogue  of the well-studied Fourier multipliers $m(D)$ given by
\[ m(D) := \mathcal{F}^{-1} m(\xi) \mathcal{F}. \]
The intertwining relation
\begin{equation}\label{intertwining}
m(D^\sharp) = \Omega m(\nabla) \Omega^*
\end{equation}
makes the boundedness of $m(D^\sharp)$ on a given space equivalent to that of $m(\nabla)$ as soon as the wave operator is bounded. 
This immediately gives crucial boundedness properties for some standard (distorted) Fourier multipliers.
 
For instance, we may define the (distorted) Littlewood-Paley operators as follows.  Pick a smooth function $\Psi$ supported on $B_2(0)$, and equal to $1$ on $B_1(0)$. Let $\Phi := \Psi(2 \cdot) - \Psi$. The Littlewood-Paley operators are then given by:
$$
P_N = \Phi\left( \frac{D^\sharp}{N} \right) \quad \mbox{and} \quad P_{<N} = \Psi \left( \frac{D^\sharp}{N} \right)\qquad N \in 2^\Z.
$$
Sometimes, we simply denote $f_N$ instead of $P_N f$ so that the decomposition of a function $f$ in dyadic frequency pieces reads
\begin{equation}\label{LittlewoodPaley}
f = \sum_{N \in 2^\mathbb{Z}} P_N f = \sum_{N \in 2^\mathbb{Z}} f_N.
\end{equation}
Finally, if $M$ is not a power of 2, it can be rounded to a power of 2, say $N$, and we set 
$P_M = P_N$.

We also consider the distorted Sobolev spaces $\dot W_{\sharp}^{s,p}(\R^d):=|D^\sharp|^{-s}L_x^p(\R^d)$ (for $-\infty<s<d$) and $W_\sharp^{s,p}(\R^d)=\langle D^\sharp\rangle^{-s} L_x^p(\R^d)$ (for all $s\in \R$). 
Thanks to the intertwining property and part (c) of Lemma \ref{square and maximal lemma} below, these spaces inherit the square function characterization from their Euclidean counterparts: for $1<p<\infty$,
$$
\|f\|_{\dot W_\sharp^{s, p}}\sim_d   \Vert  (\sum_{N \in 2^\Z} N^{2s}|P_N f|^2)^{1/2}\Vert_{L_x^p},\qquad \|f\|_{W_\sharp^{s, p}}\sim_d  \Vert (\sum_{N \in 2^\mathbb{N}} N^{2s}|P_N f|^2 )^{1/2}\Vert_{L_x^p}.
$$
Here we adopt the convention that when $P_1$ appears in a sum over dyadic blocks with $N\in \mathbb N$, it should be understood as $P_{\leq 1}$. The question of the equivalence (or not) of the distorted Sobolev spaces and their homogeneous counterparts is addressed in Section \ref{prelim:riesz}.

The following proposition will be useful throughout the paper and its proof is direct by using the intertwining property \eqref{intertwining}.

\begin{proposition}  \label{Bernsteinprop} \label{fractionalintprop} Assume {\bf H3} is satisfied. 
\begin{enumerate}
\item[(a)] {\em(Distorted Bernstein inequality)} Let $s > 0$ be given.  Then, for each $1 < p\leq q \leq \infty$, and $f \in L^q(\mathbb{R}^d)$ , 
\[ \|P_{1/\sqrt{s}} f \|_{L^q(\R^d)} + \| P_{\leq 1/\sqrt{s}} f \|_{L^q(\R^d)} \lesssim s^{-\frac{d}{2}( \frac{1}{p} - \frac{1}{q})} \| f \|_{L^p(\R^d)}.\]
\item[(b)] {\em(Distorted Sobolev inequality)} If $1 \leq p \leq q < \infty$ and $\frac{1}{p} - \frac{1}{q} = \frac{\alpha}{d}$, then
\[ \| |D^\sharp|^{-\alpha} f\|_{L^q} \lesssim \| f\|_{L^p}\quad \Leftrightarrow\quad \|f\|_{L^q} \lesssim \|f\|_{\dot{W}^{\alpha,p}_\sharp}.\]
\item [(c)] {\em(Distorted fractional integration)} For each $\alpha \geq 0$ and $t > 0$, define the (distorted) multipliers
\[ \Lambda_t^{-\alpha} := t^{\alpha/2} \varphi(\sqrt{t} |D^\sharp|)^{\alpha},\]
where $\varphi$ is a smooth function satisfying
\[ \varphi(k) = \left\{ \begin{array}{ll} 1/k, & k \geq 2, \\ 1 & k \leq 1.\end{array} \right. \]
If $1 \leq p \leq q \leq \infty$ and $\frac{1}{p} - \frac{1}{q} < \frac{\alpha}{d}$ or if $1 \leq p \leq q <\infty$ and $\frac{1}{p} - \frac{1}{q} = \frac{\alpha}{d}$, then
\[ \|\Lambda_t^{-\alpha} f \|_{L^q} \lesssim t^{\frac{\alpha}{2} - \frac{d}{2} ( \frac{1}{p}-\frac{1}{q})} \| f \|_{L^p}.\]
Finally, we will sometimes need to use this inequality in conjunction with the linear group: under the same conditions as above, and assuming furthermore $q \leq 2 \leq p$,
$$
\|\Lambda_t^{-\alpha} e^{-itH} f \|_{L^q} \lesssim t^{\frac{\alpha}{2} - \frac{d}{2} ( \frac{1}{p}-\frac{1}{q})} \| f \|_{L^p}.
$$
\end{enumerate}
\end{proposition}

\medskip
\subsection{Boundedness of the Riesz transforms}\label{prelim:riesz} The Riesz transform associated with the elliptic operator $-\Delta+V$ is given by
\begin{equation}\label{Riesz transform} 
\mathfrak {R}f:=\nabla (-\Delta+V)^{-1/2} f=:\frac{\nabla}{|D^\sharp|}f
\end{equation}
The question of boundedness of the operator above has been the subject of intensive study in the past twenty years \cite{Shen1995, Auscher2007maximal, Assaad2011, guillarmou2008resolvent} as it is equivalent to boundedness of the embedding $\dot {W_\sharp}^{1,p} \hookrightarrow \dot W^{1,p}$. It is the analogue of the well-known question originally raised by Strichartz concerning the $L^p(M)$ boundedness of Riesz transforms $\frac{\nabla}{\sqrt{-\Delta_g}}$ on complete Riemannian manifolds $(M,g)$ \cite{strichartz1983analysis}. The above mentioned references give certain classes of potentials for which the Riesz transform \eqref{Riesz transform} is bounded on various $L^p(\R^d)$ spaces. We are not concerned here with the state of the art in that field, though, since we choose not to impose these conditions on $V$ (especially that for a large class of potentials of interest to us, the operator $\mathfrak R$ is only bounded if $1<p<d$ \cite{guillarmou2008resolvent}). In fact, for the purposes of our PDE applications, we are content with the ``imperfect'' estimate \eqref{trilinearest01}.


Finally, note that if we regularize $\mathfrak{R}$ at low distorted frequencies, we directly get a bounded operator on $L^p(\R^d)$ for all $1<p<\infty$ under our assumptions on the potential. 
Indeed, for  $V$ satisfying {\bf H2} and ${\bf H3^*}$, the operator 
\begin{equation}\label{Bessel}
\mathfrak B : f \mapsto \nabla (I-\Delta+V)^{-1/2}\,f= \nabla \langle D^\sharp\rangle^{-1} \,f
\end{equation}
is bounded on $L^p(\mathbb{R}^d)$ for all $1<p<\infty$.  This follows directly by noticing that $\langle D^\sharp\rangle^{-1} =\Omega \langle \nabla \rangle^{-1}  \Omega^*$, using assumption ${\bf H3^*}$ and the boundedness of $\nabla \langle \nabla \rangle^{-1}$.

\subsection{Dispersive estimastes on $\R^3$}

Dispersive estimates for the semigroup of a Schr\"odinger operator have been studied extensively; let us mention in particular the work of Journ\'e, Sogge, and Soffer~\cite{journe1991decay}, and
the more recent papers of Rodnianski and Schlag~\cite{rodnianski2004decay}, and Beceanu~\cite{beceanu2011estimate}.  A recent survey is given in \cite{schlag2007dispersive}.

Observe in any case that the $L^p$ boundedness of the wave operator $\Omega$ corresponding to $H$ implies that the group $e^{itH}$ enjoys the same dispersive and Strichartz estimates as $e^{it\Delta}$.   For reference, we record here the dispersive estimates that we will require in Section \ref{ubar2}.  Note that these are specialized to the three-dimensional case.  

\begin{proposition}[Weighted dispersion estimates] \label{prelim:weightedispersiveprop} Let $V: \mathbb{R}^3 \to \mathbb{R}$ be a potential satisfying {\bf H1}, {\bf H2}, and {\bf H3}.   Then 
\beq \| e^{itH} f \|_{L^6_x(\mathbb{R}^3)} \lesssim \frac{1}{t} \| \langle x \rangle \Omega^* f \|_{L^2_x(\mathbb{R}^3)},  \label{ubar2:basicdispersionestimate} \eeq
and more generally
\beq \| e^{itH} f \|_{L_x^p(\mathbb{R}^3)} \lesssim \frac{1}{\displaystyle t^{3(\frac{1}{2} - \frac{1}{p})}} \| \langle x \rangle \Omega^* f\|_{L^2_x(\mathbb{R}^3)}, \qquad \textrm{for } p \in [2, 6]. \label{ubar2:interpdispersionestimate} \eeq
\end{proposition}
\begin{proof}  The second statement \eqref{ubar2:interpdispersionestimate} follows by interpolating \eqref{ubar2:basicdispersionestimate} with the conserved $L^2$ norm.  To prove \eqref{ubar2:basicdispersionestimate}, we note that
\begin{align*}
\| e^{itH} f \|_{L^6} &= \| \Omega e^{-it\Delta} \Omega^* f \|_{L^6} \\
& \lesssim \| e^{-it\Delta} \Omega^* f \|_{L^6}\\
& \lesssim \frac{1}{t} \| \langle x \rangle \Omega^* f \|_{L^2}. \end{align*}
The first inequality is a result of the boundedness of $\Omega$ discussed in Section \ref{prelim:waveopsec}, while the second is a standard dispersive estimate for the Schr\"odinger semigroup. Note that $\| \langle x \rangle \Omega^* f \|_{L^2}\sim \|f\|_{L^2}+\|\partial_\xi f^\sharp\|_{L^2}$.
\end{proof}

\section{Multilinear harmonic analysis with potential} \label{multestimates section}
In this section, we develop the multilinear analysis needed to build a space-time resonance theory in the inhomogeneous setting. 
We start in Section \ref{multilinearsection} with the basics, 
and proceed to prove in section~\ref{secproof} the analogues of Coifman-Meyer theorems for the distorted Fourier transform on $\R^d$ ($d\geq 2$) as stated in Theorem \ref{intro:CMtheorem}. 
Then, in Section \ref{commutator subsection}, we prove estimates for the commutators of the wave operator and position operators, such as $[|x|, \Omega]$. 
Here, we restrict ourselves to $\R^3$ as we rely on Yajima's explicit description of the wave operator $\Omega$ in \cite{yajima1995Wkp}.  
We expect similar results to hold in other dimensions using the work in \cite{yajima1995Wkp3,jensen2002remark}. 
These estimates will be instrumental in Section \ref{derivatives of M section} where we attempt to understand the behavior of derivatives of the distribution $\mathcal M$ that was
 discussed in Section \ref{intro:spacetimeres}.

\subsection{Definitions and first results} \label{multilinearsection}

We start by considering pseudo-product operators of the form:
\begin{equation}\label{CMC0}
T(f,g)(x):=\int_{\R^3}\int_{\R^3} m_0(\xi_1, \xi_2) f^\sharp(\xi_1) g^\sharp(\xi_2) e(x,\xi_1) e(x, \xi_2) \, d\xi_1 \, d\xi_2
\end{equation}
Note that the case $m=1$ corresponds (up to a constant factor) to the product of $f$ and $g$. We say that the multiplier $m_0$ satisfies Coifman-Meyer type bounds if the following homogeneous bounds hold for sufficiently many multi-indices $\alpha$ and $\beta$:

\begin{equation}\label{CMC1}
|\partial_{\xi_1}^\alpha \partial_{\xi_2}^\beta m_0(\xi_1, \xi_2)|\leq C (|\xi_1|+|\xi_2|)^{-(|\alpha|+|\beta|)}
\end{equation}

We would like to prove estimates of the form:
\begin{equation}\label{CME1}
\|T(f,g)\|_{L^{r'}(\R^d)}\lesssim_{m_0, V} \|f\|_{L^p(\R^d)}\|g\|_{L^q(\R^d)}
\end{equation}
whenever $f,g \in \mathcal S(\R^d)$ and $\frac{1}{r'}=\frac{1}{p}+\frac{1}{q}$. By duality, this  is equivalent to proving that, for all $f,g,h \in \mathcal S(\R^d)$ and $p,q,r$ satisfying $1=\frac{1}{p}+\frac{1}{q}+\frac{1}{r}$ (and $\frac{1}{r'}=1-\frac{1}{r}$),
\begin{align}\label{CME2}
\Lambda(f,g,h) := \int_{\R^d} T(f,g)(x) h(x) dx \lesssim \|f\|_{L^p(\R^d)} \|g\|_{L^q(\R^d)}\|h\|_{L^r(\R^d)}.
\end{align}
The  left-hand side above can be written
\begin{align*}
\int_{\R^d} T(f,g) h \, dx=&\int_{\R^d}\int_{\R^d} m(\xi_1, \xi_2) f^\sharp(\xi_1) g^\sharp(\xi_2) \int_{\R^3_x} h(x) e(x,\xi_1) e(x,\xi_2) \, dx \,  d\xi_1\,  d\xi_2\\
=:& \int_{\R^d}\int_{\R^d} m(\xi_1, \xi_2) f^\sharp(\xi_1) g^\sharp(\xi_2) \langle \mathcal{M}(\xi_1,\xi_2), h^\sharp\rangle\,  d\xi_1 d\xi_2
\end{align*}
where we denoted by $\MM(\xi_1, \xi_2)\in \mathcal{S}^\prime(\R^d)$ the distribution determined 
by
\begin{equation}\label{def of M}
\langle \MM(\xi_1, \xi_2), \phi\rangle =\int_{\R^d} e(x,\xi_1) e(x,\xi_2) {\FF^\sharp}^{-1}\phi(x) \, dx, \qquad \phi \in \mathcal{S}(\mathbb{R}^d).
\end{equation}
One way to see that this is well-defined is to recall from \eqref{prelim:intertwineidentities} that ${\mathcal{F}^\sharp}^{-1} = \Omega \mathcal{F}^{-1}$, and hence by the boundedness of $\Omega$ on $L^1$, ${\mathcal{F}^\sharp}^{-1} \phi \in L^1(\R^d)$ whenever $\phi \in \mathcal S(\R^d)$.  We remark also that if $V=0$, then $\MM(\xi_1, \xi_2)=\delta (\cdot +\xi_1+\xi_2)$. For convenience, we will abuse notation throughout and denote
$$
\langle \MM(\xi_1, \xi_2), \phi \rangle=\int_{\R^d} \phi(\xi_3) \MM(\xi_1, \xi_2 , \xi_3)\, d\xi_3.
$$
As a result, we will write
\begin{equation}\label{CME3}
\Lambda(f,g,h) = \int_{\R^d} T(f,g) h \, dx=\int_{\R^d}\int_{\R^d} m(\xi_1, \xi_2) f^\sharp(\xi_1) g^\sharp(\xi_2) h^\sharp(\xi_3) \mathcal{M}(\xi_1,\xi_2, \xi_3) \, d\xi_1 \,  d\xi_2 \, d\xi_3
\end{equation}
with the understanding explained above.

We will also be interested in generalizations of $\eqref{CME3}$ given by
\begin{equation}\label{CME}
\Lambda(f, g,h) :=\int_{\R^d}\int_{\R^d} m(\xi_1, \xi_2,\xi_3) f^\sharp(\xi_1) g^\sharp(\xi_2) h^\sharp(\xi_3) \mathcal{M}(\xi_1,\xi_2, \xi_3) d\xi_1 d\xi_2 d\xi_3.
\end{equation}
where $m$ is also a Coifman-Meyer multiplier in the three variables $\xi_1, \xi_2, \xi_3$ in the sense that 
\begin{equation}\label{CMC2}
|\partial_{\xi_1}^\alpha \partial_{\xi_2}^\beta \partial_{\xi_3}^\gamma m(\xi_1, \xi_2,\xi_3)|\leq C_{\alpha\beta\gamma} (|\xi_1|+|\xi_2|+|\xi_3|)^{-(|\alpha|+|\beta|+|\gamma|)}.
\end{equation}
Our aim will be to prove the following two estimates on $\Lambda(f,g,h)$: For $p, q, r, \tilde p, \tilde q, \tilde r \in (0, \infty)$ satisfying $\frac{1}{q}+\frac{1}{p}+\frac{1}{r}=1$ and $\frac{1}{\tilde q}+\frac{1}{\tilde p}+\frac{1}{\tilde r}=1+\epsilon$, we have 
\begin{equation}\label{trilinearest0}
|\Lambda(f,g,h)| \lesssim_{C_{\alpha\beta\gamma}, V, \epsilon} \|f\|_{L^q(\R^d)} \|g\|_{L^p(\R^d)}\|h\|_{L^r(\R^d)} 
\end{equation}
assuming the boundedness of $\mathfrak R$ on $L^q, L^p,$ and $L^r$; and
\begin{equation}\label{trilinearest}
|\Lambda(f,g,h)| \lesssim_{C_{\alpha\beta\gamma}, V, \epsilon} \|f\|_{L^q(\R^d)} \|g\|_{L^p(\R^d)}\|h\|_{L^r(\R^3)}+\|f\|_{L^{\tilde q}(\R^d)} \|g\|_{L^{\tilde p}(\R^d)}\|h\|_{L^{\tilde r}(\R^d)} 
\end{equation}
under no boundedness assumption on $\mathfrak R$. This directly gives Theorem \ref{intro:CMtheorem} by duality.

\begin{remark}
As will be evident from the proof, the second term on the right-hand side of \eqref{trilinearest} can be further refined by projecting at least two of the functions $f,g,h$ onto frequencies $|D^\sharp|\leq 1$.
\end{remark}

\begin{remark} \label{msCMremark}
Suppose that $m$ is a Coifman--Meyer symbol as in \eqref{CMC2}.  For each $s > 0$, we may define a symbol $m_s$ by
\[ m_s(\xi_1, \xi_2, \xi_3) := m(\sqrt{s} \xi_1, \sqrt{s} \xi_2, \sqrt{s} \xi_3).\]
Then $m_s$ is itself Coifman--Meyer with the same constants $C_{\alpha\beta\gamma}$.  We may therefore replace $m$ by $m_s$ in either of the previous two results, and the resulting estimates will hold uniformly in $s > 0$. 
\end{remark}

The remainder of this subsection is devoted to the proof of estimates \eqref{trilinearest0} and \eqref{trilinearest}.  Before we begin, though, we will need the following maximal and square function estimates.

\begin{lemma} \label{square and maximal lemma}
Suppose that $W$ is an operator that is point-wise bounded by an $L^p$-bounded positive operator, i.e. satisfying the point-wise bound 
\[ |Wf(x)|\leq C\widetilde W|f|(x), \qquad \textrm{for all } f \in L^p(\mathbb{R}^d),~x \in \mathbb{R}^d,\]  for some positive operator $\widetilde W$ that is bounded on $L^p(\R^d)$ for $1\leq p \leq \infty$.
\begin{enumerate} 
\item[(a)] Let $\psi\in C_0^\infty(\R^d)$.  For each $n \in \mathbb{R}^d$, the operators
$$
f\mapsto \sup_{N \in 2^\Z}|We^{2\pi i \frac{n \cdot \nabla}{N}} \psi(\frac{\nabla}{N}) f| \quad\hbox{ and }\quad f\mapsto \sup_{\substack{N_1, N_2 \in 2^\Z\\N_1\geq N_2}}|We^{2\pi i \frac{n \cdot \nabla}{N_1}} \psi(\frac{\nabla}{N_2}) f|
$$
are bounded on $L^p$ for all $1<p\leq \infty$ with a bound $\lesssim \langle n \rangle^d$. 

\item[(b)]  For each $n \in \mathbb{R}^d$, the operators
\begin{equation}\label{max1}
f\mapsto \sup_{N \in 2^\Z} |e^{2\pi i \frac{n \cdot D^\sharp}{N}}\psi(\frac{D^\sharp}{N})f| \quad \hbox{ and }\quad f\mapsto \sup_{\substack{N_1,N_2 \in 2^\Z\\N_1\geq N_2}} |e^{2\pi i \frac{n \cdot D^\sharp}{N_1}}\psi(\frac{D^\sharp}{N_2})f| 
\end{equation}
are bounded on $L^p(\R^d)$ for all $1<p\leq \infty$ with a bound $\lesssim \langle n\rangle^d$.

\item[(c)] Let $U$ be any bounded operator on $L^p$ for some $1\leq p < \infty$ and suppose that $\{f_n\} \subset L^p(\mathbb{R}^d)$ is a sequence of functions. Then
\begin{equation*}
\|(\sum_{n \in \Z} |U f_n|^2)^{1/2}\|_{L^p(\R^d)} \lesssim  \|(\sum_{n\in \Z} |f_n|^2)^{1/2}\|_{L^p(\R^d)},
\end{equation*}
whenever the right-hand side is finite.

\item[(d)] Moreover, if $\phi$ is smooth and supported on an annulus, the operator
\begin{equation}\label{FStype}
f \mapsto (\sum_{N_2\in 2^{\Z}} \sup_{N_1 \geq N_2} |e^{2\pi i \frac{n \cdot D^\sharp}{N_1}} \phi(\frac{D^\sharp}{N_2})f|^2 )^{1/2}
\end{equation}
is bounded on $L^p(\R^d)$ for all $1<p<\infty$ with bound $\lesssim \langle n \rangle^d$.
\end{enumerate}
\end{lemma}
\begin{proof} (a) Since $|W f|\lesssim \widetilde W |f|$, this follows directly from the fact that the operators
$$
f \mapsto \sup_{N} |e^{2\pi i \frac{n \cdot \nabla}{N}}\psi(\frac{\nabla}{N}) f|\quad \hbox{ and }\quad f \mapsto \sup_{N_1\geq N_2} |e^{2\pi i \frac{n \cdot \nabla}{N_1}}\psi(\frac{\nabla}{N_2}) f|
$$
are bounded on $L^p(\R^d)$ for $1 < p \leq \infty$. The latter follows from the point-wise inequality 
\begin{equation}\label{pm1}
|\psi(\frac{\nabla}{N})f(x-y)|\lesssim \langle N|y|\rangle^d Mf(x),
\end{equation}
where $Mf$ is the Hardy-Littlewood maximal function. The proof of this inequality is elementary, but we include it for completeness. Assuming without loss of generality that $f\geq 0$ and picking a radial majorant  for $|\mathcal F^{-1} \psi|\lesssim \sum_{j\leq 0} 2^j \chi_{B_{R_j}(0)}$ with $\sum_{j\leq 0} 2^j R_j^d <\infty$, we can bound:
\begin{align*}
\psi(\frac{\nabla}{N})f(x-y)&:=  N^d\int_{\R^d} \mathcal F^{-1} \psi(Nz) f(x-y-z)dz= N^d\int_{\R^d} \mathcal F^{-1} \psi(N(z-y)) f(x-z)dz\\
&\lesssim \sum_{j\leq 0}2^j R_j^d \fint_{B_{R_j/N(y)}} f(x-z) dz\\
&\lesssim 
\begin{cases} \displaystyle \sum_{j\leq 0}2^j R_j^d \fint_{B_{1000R_j/N}(0)} f(x-z) dz\lesssim Mf(x) & \hbox{if } |y|\leq 100R_j/N \\
\displaystyle \sum_{j\leq 0}2^j ({N|y|})^d \fint_{B_{2|y|}(0)} f(x-z) dz\lesssim (N|y|)^dMf(x) & \hbox{if } |y|\geq 99R_j/N \\
\end{cases}\\
&\lesssim \langle N|y|\rangle^d Mf(x).
\end{align*}
\\

(b) Recall that by the intertwining property \eqref{intertwining}, $m(D^\sharp)f =\Omega m(\nabla) g$ where $g=\Omega^* f$. Since $\Omega^*$ is bounded on $L^p$ for all $1\leq p \leq \infty$, it is enough to show that the operator
$$
f\mapsto \sup_{N\in 2^\Z}|\Omega e^{2\pi i \frac{n \cdot \nabla}{N}}\psi(\frac{\nabla}{N})g|
$$
is bounded on $L^p(\R^d)$ for all $1<p\leq \infty$. But this follows from part (a) and Remark \ref{boundedbypositive}. \\

(c) Let $U$ be a bounded operator on $L^p$, for $1 \leq p \leq \infty$.  By density, we may assume that  $\{f_n\}$ is a finite sequence with $n \in \{-L, \ldots, L\}$. Let $\{\epsilon_{n}\}$ be a sequence of i.i.d. random variables with $P(\epsilon_n = \pm 1) = 1/2$ for $n \in \{-L,\ldots, L\}$. By Khinchine's inequality we have the point-wise equivalence
\begin{align*}
( \sum_{n} |U f_n|^2)^{p/2}\sim \mathbb{E}[  |\sum_{n}\epsilon_n Uf_n|^p ] =\mathbb{E} [  |U\sum_{n} \epsilon_n f_n |^p ].
\end{align*}
Here $\mathbb{E}$ denotes the expected value.  Now, integrating both sides over $\R^d$, we get using Fubini's theorem that
\begin{align*}
 \| ( \sum_{n} |U f_n|^2 )^{1/2} \|_{L^p}^p & \lesssim \mathbb{E}\int_{\R^d} |U\sum_{n} \epsilon_n f_n |^p \,dx \\
& \lesssim \mathbb{E}\int_{\R^d}  |\sum_{n} \epsilon_n f_n |^p \,dx \sim \|( \sum_{n} | f_n|^2)^{1/2}\|_{L^p}^p.
\end{align*}\\

(d) We start by noticing that the operator
\begin{equation}\label{int001}
f \mapsto \sum_{N_2 \in 2^\Z} ( \sup_{N_1\geq N_2} |e^{\frac{2\pi i n \cdot \nabla}{N_1}}\phi(\frac{\nabla}{N_2}) f|^2 )^{1/2}
\end{equation}
is bounded from $L^p(\R^d)$ to $L^p(\R^d)$. This follows by writing for any $N_1 \geq N_2$:
\[|e^{\frac{2\pi i n \cdot \nabla}{N_1}}\phi(\frac{\nabla}{N_2}) f|=|e^{\frac{2\pi i n \cdot \nabla}{N_1}}\psi(\frac{\nabla}{N_1})f_{N_2}|
\lesssim \langle n\rangle^d M f_{N_2}\]
by part (a). The boundedness of $\left(\sum_{N_2}|Mf_{N_2}|^2\right)^{1/2}$ is the Fefferman-Stein maximal inequality (cf. \cite{Fefferman1971, stein1993harmonicbook}).  The statement in (d) is now a direct consequence of part (c) with $U=\widetilde{\Omega}$, where $\tilde \Omega$ is the positive operator majoring $\Omega$ as in Remark \ref{boundedbypositive}.
\end{proof}

\subsection{Proof of Theorem \ref{intro:CMtheorem}}
\label{secproof}

With the maximal and square function estimates established, we are now prepared to prove the first main theorem.

\begin{proof}[Proof of Theorem \ref{intro:CMtheorem}] We will only consider the situation when $m$ satisfies \eqref{CMC2}.  The other case is similar (in fact easier).

\subsubsection*{Step 1: decomposition of $\Lambda$}  We start by decomposing $f, g, h$ into Littlewood-Paley pieces with respect to the distorted Fourier transform as in \eqref{LittlewoodPaley},
\[ f=\sum_{N_1 \in 2^\Z} f_{N_1}, \qquad g=\sum_{N_2 \in 2^\Z} g_{N_2},\qquad h=\sum_{N_3 \in 2^\Z} h_{N_3}.\]
As a result, we get that
$$
\Lambda(f,g,h)=\sum_{N_1,N_2, N_3} \Lambda(f_{N_1}, g_{N_2}, h_{N_3}).
$$

By symmetry, we may assume without any loss of generality that $N_3 \leq N_2 \leq N_1$.  Abusing notation somewhat, we will continue to denote the corresponding sum by $\Lambda$. 
Next, we let $\tilde \phi \in C_0^\infty(\R)$ be given such that $\tilde \phi \phi=\phi$.  Define $\tilde m^{N_1}$ by the relation
$$
m(\xi_1, \xi_2, \xi_3)\tilde \phi(\frac{\xi_1}{N_1})\tilde \phi(\frac{\xi_2}{N_1})\tilde \phi(\frac{\xi_3}{N_1})=:\tilde m^{N_1}(\frac{\xi_1}{N_1}, \frac{\xi_2}{N_1}, \frac{\xi_3}{N_1}).
$$
Then $\tilde m^{N_1} \in C_c^{10}([-K/2,K/2]^{3d})$ for some constant $K$ depending on $\tilde \phi$ with a uniformly bounded $C^{10d+10}$ norm\footnote{Of course, this much regularity in $m$ is not necessary, but it is sufficient for our purposes.} (independent of $N_1, N_2, N_3$). 
Expanding $\tilde m^{N_1}$ in a Fourier series, we can write, if $(\xi_1,\xi_2,\xi_3) \in [-K/2,K/2]^{3d}$,
$$
\tilde m^{N_1}(\xi_1, \xi_2, \xi_3) =\sum_{n_1, n_2, n_3 \in \Z^d} a^{N_1}(n_1,n_2, n_3) e^{\frac{2\pi i}{K}(n_1.\xi_1+n_2.\xi_2+n_3.\xi_3)}
$$
where $a^{N_1}$ satisfies the bound:
\begin{equation}\label{bound on a}
|a(n_1, n_2,n_3)|\lesssim_{m, \phi} (1+|n_1|+|n_2|+|n_3|)^{-10d}.
\end{equation}
Consequently,
\begin{equation}\label{expansion of Lambda}
\Lambda(f, g, h)=\sum_{N_3 \leq N_2 \leq N_1} \sum_{n_1,n_2,n_3 \in \Z^d} a^{N_1}(n_1, n_2, n_3) \Xi_{N_1, N_2, N_3}^{n_1, n_2, n_3}
\end{equation}
where
\begin{align}\label{def of Xi}
\Xi_{N_1, N_2, N_3}^{n_1, n_2,n_3}:=&\int \!\!\! \int \!\!\! \int \left(e^{\frac{2\pi i}{N_1K}n_1 \cdot \xi_1}f^\sharp_{N_1}(\xi_1)\right)\left(e^{\frac{2\pi i}{N_1K}n_2 \cdot \xi_2}g^\sharp_{N_2}(\xi_2)\right)
\left(e^{\frac{2\pi i}{N_1K}n_3 \cdot \xi_3}h^\sharp_{N_3}(\xi_3)\right) \nonumber \\ 
&\qquad {} \hspace{7cm} \MM(\xi_1, \xi_2, \xi_3) \, d\xi_1 \, d\xi_2 \, d\xi_3\\
=&\int \!\!\! \int \!\!\! \int f^\sharp_{N_1, n_1}(\xi_1) g^\sharp_{N_2, n_2, N_1} (\xi_2) h^\sharp_{N_3,n_3,N_1}(\xi_3)\MM(\xi_1, \xi_2, \xi_3)\, d\xi_1\, d\xi_2\, d\xi_3\nonumber \\
=&\int \!\!\! \int  f^\sharp_{N_1, n_1}(\xi_1) g^\sharp_{N_2, n_2, N_1} (\xi_2) \Gamma(h_{N_3, n_3, N_1})(\xi_1, \xi_2)\, d\xi_1 d\xi_2.
\end{align}
Here we have denoted  \[
f_{N_1, n_1}^\sharp := e^{\frac{2\pi i}{N_1K}n_1 \cdot \xi_1}f^\sharp_{N_1}, \qquad
 g_{N_2, n_2,N_1}^\sharp := e^{\frac{2\pi i}{N_1K}n_2 \cdot \xi_2}g^\sharp_{N_2}, \qquad
 h_{N_3, n_3, N_1}^\sharp :=e^{\frac{2\pi i}{N_1K}n_3 \cdot \xi_3}h^\sharp_{N_3},\]
 and
$$
\Gamma(h_{N_3, n_3, N_1})(\xi_1, \xi_2) :=\langle \mathcal{M}(\xi_1,\xi_2),h_{N_3, n_3, N_1}^\sharp \rangle =\int_{\R^d} e(x,\xi_1) e(x,\xi_2) h_{N_3,n_3}(x) \, dx.
$$
Writing 
\[ h_{\leq N_2, n_3, N_1} =\sum_{N_3 \leq N_2} h_{N_3, n_3, N_1}, \qquad  \Xi_{N_1, N_2}^{n_1, n_2,n_3}=\sum_{N_3 \leq N_2}\Xi_{N_1, N_2, N_3}^{n_1, n_2,n_3}\]
we may further simplify \eqref{expansion of Lambda} to obtain
\begin{align}\label{Lambda estimate}
\nonumber |\Lambda (f, g, h)| &\leq \sum_{N_2 \leq N_1}\sum_{n_1, n_2, n_3 \in \Z^d}|a(n_1, n_2, n_3)| |\Xi_{N_1, N_2}^{n_1, n_2,n_3}|\\ 
&\lesssim_{m, \phi}  \sum_{n_1, n_2,n_3 \in \Z^d} (1+|n_1|+|n_2|+|n_3|)^{-10d-10} \sum_{N_2 \leq N_1} |\Xi_{N_1, N_2}^{n_1, n_2,n_3}|.
\end{align}

As a result, we are reduced to proving the following estimate:
\begin{equation}\label{poly growth}
\sum_{N_2 \leq N_1} |\Xi_{N_1, N_2}^{n_1, n_2, n_3}| \lesssim \langle n_1\rangle^{d} \langle n_2\rangle^{d} \langle n_3\rangle^{d}
\begin{cases}\|f\|_p \|g\|_q \|h\|_r \qquad &\hbox{for part (i)},\\
\|f\|_p \|g\|_q \|h\|_r + \|f\|_{\tilde{p}} \|g\|_{\tilde{q}} \|h\|_{\tilde{r}} \qquad &\hbox{for part (ii)}.\\
\end{cases}
\end{equation}

\subsubsection*{Step 2: Proof of \eqref{trilinearest0}}
We start by computing using Green's formulae:
\begin{equation}\label{Green's trick}
\begin{split}
|\xi_1|^2\Gamma(\xi_1, \xi_2)&= \int_{\R^d} He(x,\xi_1) e(x, \xi_2) h_{\leq N_2, n_3, N_1}(x) \, dx\\
&=\int_{\R^d} e(x,\xi_1)H[e(x, \xi_2) h_{\leq N_2, n_3, N_1}(x)] \, dx\\
&= |\xi_2|^2 \Gamma(\xi_1, \xi_2) +\int_{\R^d} e(x,\xi_1)e(x, \xi_2) H  [h_{\leq N_2, n_3, N_1}(x)] \, dx\\
&\qquad-\int_{\R^d} e(x,\xi_1)e(x, \xi_2) V(x)h_{\leq N_2, n_3, N_1}(x)\, dx \\ 
& \qquad +2\int_{\R^d} e(x,\xi_1)\nabla e(x, \xi_2) \nabla h_{\leq N_2, n_3, N_1}(x) \, dx\\
\end{split}
\end{equation}
Thus,
\begin{align}
\Xi_{N_1, N_2}^{n_1, n_2,n_3}=&\int \!\!\! \int \frac{1}{|\xi_1|^2}f^\sharp_{N_1, n_1}(\xi_1) |\xi_2|^2g^\sharp_{N_2, n_2, N_1} (\xi_2) \Gamma(h_{\leq N_2, n_3, N_1})(\xi_1, \xi_2) \, d\xi_1 \, d\xi_2 \label {x01}\\
&+\int \!\!\! \int \frac{1}{|\xi_1|^2}f^\sharp_{N_1, n_1}(\xi_1) g^\sharp_{N_2, n_2, N_1} (\xi_2) \Gamma(Hh_{\leq N_2, n_3, N_1})(\xi_1, \xi_2) \, d\xi_1 \, d\xi_2 \label{x02}\\
&+\int \!\!\! \int\frac{1}{|\xi_1|^2}f^\sharp_{N_1, n_1}(\xi_1) g^\sharp_{N_2, n_2, N_1} (\xi_2) \Gamma( V h_{\leq N_2, n_3, N_1})(\xi_1, \xi_2) \, d\xi_1 \, d\xi_2\label{x03}\\
&+2\int \!\!\! \int \!\!\! \int \frac{1}{|\xi_1|^2}f^\sharp_{N_1, n_1}(\xi_1) g^\sharp_{N_2, n_2, N_1} (\xi_2)   e(x,\xi_1) \nabla e(x,\xi_2) \nabla h_{\leq N_2, n_3, N_1}\, dx \, d\xi_1\, d\xi_2 \label{x04}.
\end{align}

We start with the contribution of $\eqref{x01}$: Denoting by 
\begin{equation}\label{tildeg}
\widetilde g_{N_2, n_2, N_1}:= {\FF^\sharp}^{-1}\frac{|\xi|^2}{N_2^2} \FF^\sharp g_{N_2, n_2, N_1}={\FF^\sharp}^{-1} e^{2\pi i \frac{n_2. \xi}{KN_1}}\widetilde \phi(\frac{\xi}{N_2}) \FF^\sharp g,\qquad \hbox{with } \widetilde \phi (\xi):=|\xi|^{2} \phi(\xi) 
\end{equation}
and 
\beq \label{und f}
\underline f_{N_1, n_1}:= {\FF^\sharp}^{-1}\frac{N_1^2}{|\xi|^2} \FF^\sharp f_{N_1, n_1}={\FF^\sharp}^{-1}  e^{2\pi i \frac{n_1\cdot \xi}{KN_1}} \underline \phi(\frac{\xi}{N_1}) \FF^\sharp f,\qquad \hbox{with } \underline \phi(\xi):=|\xi|^{-2} \phi(\xi).\eeq
We now estimate
\begin{align*}
|\sum_{ N_2 \leq N_1} \eqref{x01}| \leq& \sum_{N_2 \leq N_1} \left(\frac{N_2}{N_1}\right)^2 \left|\int_{\R^d}\int_{\R^d}\underline f^\sharp_{N_1, n_1}(\xi_1) \tilde g^\sharp_{N_2, n_2, N_1} (\xi_2) \Gamma(h_{\leq N_2, n_3, N_1})(\xi_1, \xi_2) \, d\xi_1 \, d\xi_2\right| \\
=&\sum_{ N_2 \leq N_1} \left(\frac{N_2}{N_1}\right)^2 \left|\int_{\R^d} \underline f_{N_1, n_1}(x) \tilde g_{N_2, n_2, N_1} (x) ( e^{2\pi i \frac{n_3 \cdot D^\sharp}{KN_1}} \psi(\frac{D^\sharp}{N_2}) h(x)) \, dx\right| \\
\lesssim& \int_{\R^d} (\sum_{N_1} |\underline f_{N_1, n_1}|^2)^{1/2}  (\sum_{N_2} \sup_{N_1 \geq N_2}|\tilde g_{N_2, n_2, N_1}|^2)^{1/2} (\sup_{N_1\geq N_2} e^{2\pi i \frac{n_3 \cdot D^\sharp}{KN_1}} \psi(\frac{D^\sharp}{N_2}) h(x)) \, dx \\
\lesssim& \langle n_1\rangle^{d} \langle n_2\rangle^{d} \langle n_3\rangle^{d} \|f\|_{L^q(\R^d)}\|g\|_{L^p(\R^d)}\|h\|_{L^r(\R^d)}.
\end{align*}
The last inequality above follows from Lemma \ref{square and maximal lemma}. The bound on \eqref{x02} is similar: Denoting by $\tilde \psi(\xi):=|\xi|^2 \psi(\xi)$, 
\begin{align*}
|\sum_{N_2 \leq N_1} \eqref{x02}| \leq& \sum_{N_2 \leq N_1} N_1^{-2} \left|\int_{\R^d}\int_{\R^d}\underline f^\sharp_{N_1, n_1}(\xi_1) g^\sharp_{N_2, n_2, N_1} (\xi_2) \Gamma((-\Delta +V)h_{\leq N_2, n_3, N_1})(\xi_1, \xi_2) d\xi_1 d\xi_2\right| \\
\leq&\sum_{ N_3 \leq N_2 \leq N_1} \left(\frac{N_3}{N_1}\right)^2 \left|\int_{\R^d} f_{N_1, n_1}(x) g_{N_2, n_2, N_1} (x) h_{N_3, n_3, N_1}(x) dx\right| \\
\lesssim& \int_{\R^d} (\sum_{N_1} |\underline f_{N_1, n_1}|^2)^{1/2}  (\sup_{N_2\leq N_1} |g_{N_2, n_2, N_1}|)\left(\sum_{N_3\leq N_1} |h_{N_3, n_3, N_1}|^2\right)^{1/2} dx \\
\lesssim& \langle n_1\rangle^{d} \langle n_2\rangle^{d} \langle n_3\rangle^{d} \|f\|_{L^q(\R^d)}\|g\|_{L^p(\R^d)}\|h\|_{L^r(\R^d)},
\end{align*}
by Lemma \ref{square and maximal lemma}.

We move on to studying the contribution of $\eqref{x03}$. Let $s_1, s_2, s_3>0$ satisfy $s_1 +s_2+s_3=2$, $s_1 \leq d/q$, $s_2\leq d/p$, and $s_3 \leq d/r$.
Then,
\begin{align*}
|\sum_{N_2 \leq N_1} \eqref{x03}| & \leq \sum_{N_2 \leq N_1} N_1^{-2} \left|\int \!\!\! \int \underline f^\sharp_{N_1, n_1}(\xi_1) g^\sharp_{N_2, n_2, N_1} (\xi_2) \Gamma(Vh_{\leq N_2, n_3, N_1})(\xi_1, \xi_2) \, d\xi_1\, d\xi_2\right| \\
&=\sum_{N_2 \leq N_1} \frac{N_2^{s_2+s_3}}{N_1^{2-s_1}} \left|\int ( N_1^{-s_1}\underline f_{N_1, n_1}) (N_2^{-s_2} g_{N_2, n_2, N_1})  (N_2^{-s_3}h_{\leq N_2, n_3, N_1}) V(x) \,  dx\right| \\
&\lesssim \int (\sum_{N_1}N_1^{-2s_1}|\underline f_{N_1, n_1}|^2 )^{1/2} (\sum_{N_2}N_2^{-2s_2}\sup_{N_1\geq N_2}|g_{N_2, n_2, N_1}|^2 )^{1/2} \\
& \qquad \cdot  (\sup_{N_2, N_1} N_2^{-s_3}e^{2\pi i \frac{n_3 \cdot D^\sharp}{KN_1}} \psi(\frac{D^\sharp}{N_2})h ) V(x) \, dx\\
&\lesssim \langle n_1 \rangle^{d}\langle n_2 \rangle^{d}\langle n_3 \rangle^{d}\|f\|_{\dot W_\sharp^{-s_1, \tilde q}(\R^3)}\|g\|_{\dot W_\sharp^{-s_2, \tilde p}(\R^d)}\|h\|_{\dot W_\sharp^{-s_3,\tilde r}(\R^d)} \|V\|_{L^{d/2}(\R^d)}\\
&\lesssim \langle n_1 \rangle^{d}\langle n_2 \rangle^{d}\langle n_3 \rangle^{d}\|f\|_{L^q(\R^d)}\|g\|_{L^p(\R^d)}\|h\|_{L^r(\R^d)}
\end{align*}
where we denoted $\frac{1}{\tilde q}=\frac{1}{q}-\frac{s_1}{d}, \frac{1}{\tilde p}=\frac{1}{p}-\frac{s_2}{d}, \frac{1}{\tilde r}=\frac{1}{r}-\frac{s_3}{d}$ and used Lemma \ref{square and maximal lemma} 
and Sobolev embedding. We also remark here that in order to control $(\sup_{N_2, N_1} N_2^{-s_3}e^{2\pi i \frac{n_3 \cdot D^\sharp}{KN_1}} \psi(\frac{D^\sharp}{N_2})h)$ in terms of $\|h\|_{\dot W_\sharp^{-s_3,r}(\R^3)}$, we repeated the proof of part (a) of Lemma \ref{square and maximal lemma} and bounded:
$$
|\frac{|\nabla|^s}{N_2^s} \psi(\frac{\nabla}{N_2})f| =|\sum_{N_3 \leq N_2} (\frac{N_3}{N_2})^{s_3} \tilde \phi^\prime(\frac{\nabla}{N_3})f| \lesssim Mf(x)
$$
where $\tilde \phi^\prime:= |\xi|^s \left(\psi (\xi)-\psi(2\xi)\right) \in C_0^\infty (\R^3)$ (recall $|\tilde \phi(\nabla/N)f|\lesssim Mf$).

Finally, we bound the contribution of \eqref{x04} under the assumption that $\mathfrak R$ in \eqref{Riesz transform} is bounded on $L^p(\R^d)$. 

\begin{align*}
| \sum_{N_2 \leq N_1} \eqref{x04}| \leq& \sum_{N_2 \leq N_1} N_1^{-2} \left|\int_{\R^d} \underline f_{N_1, n_1}(x) \nabla g_{N_2, n_2,N_1} \nabla h_{\leq N_2, n_3, N_1} dx \right| \\
\leq& \sum_{N_3\leq N_2 \leq N_1} N_1^{-2} \left|\int_{\R^d} \underline f_{N_1, n_1}(x) \nabla g_{N_2, n_2,N_1} \nabla h_{N_3, n_3, N_1} dx \right|\\
=& \sum_{N_3\leq N_2 \leq N_1} \frac{N_2 N_3}{N_1^{2}} \left|\int_{\R^d} \underline f_{N_1, n_1}(x) \mathfrak{R} g_{N_2, n_2,N_1} \mathfrak{R} h_{N_3, n_3, N_1} dx \right|\\
=& \sum_{\substack{N_3\leq N_2\\ J \in 2^\N}} \frac{N_3}{J^2 N_2} \left|\int_{\R^d} \underline f_{JN_2, n_1}(x) \mathfrak{R} g_{N_2, n_2,JN_2} \mathfrak{R} h_{N_3, n_3, JN_2} dx \right|\\
=& \sum_{\substack{N_2\in 2^\Z\\ J, L \in 2^\N}} \frac{1}{J^2 L} \left|\int_{\R^d} \underline f_{JN_2, n_1}(x) \mathfrak{R} g_{N_2, n_2,JN_2} \mathfrak{R} h_{N_2/L, n_3, JN_2} dx \right|\\
\leq& \sum_{\substack{J, L \in 2^\N}} \frac{1}{J^2 L} \int_{\R^d} \sup_{N}|\underline f_{N, n_1}| (\sum_{N_2}|\mathfrak{R} g_{N_2, n_2,JN_2}|^2 )^{1/2} ( \sum_{N_2}|\mathfrak{R} h_{N_2, n_3, JLN_2}|^2)^{1/2} dx.
\end{align*}
Using the boundedness of $\mathfrak R$ along with Lemma \ref{square and maximal lemma}, we conclude from above that
\beq | \sum_{N_2 \leq N_1} \eqref{x04}|\lesssim \sum_{\substack{J, L \in 2^\N}} \frac{1}{J^2 L} \langle n_1 \rangle^d \langle n_2 \rangle^d \langle n_3 \rangle^d \|f\|_{L^q}\|g\|_{L^p}\|h\|_{L^r}. \label{bound on 3.24R} \eeq
This finishes the proof of part (i).

\subsubsection*{Step 3: proof of part (ii)} Here we make no assumption on the boundedness  of the Riesz transform $\mathfrak R$ associated to $V$. Instead, we assume that the potential $V$ satisfies assumption ${\bf H3^*}$, which implies the $L^p(\R^d)$ boundedness of the operator $\mathfrak B$ defined in \eqref{Bessel}.

We split the analysis into three cases, depending on the size of $N_1$ and $N_2$:

\textbf{Case 1: $N_1 \geq 1$ and $N_2 \geq  1$.}  In this case, we compute using Green's formula as in \eqref{Green's trick}. The only departure from the proof of part (i) is in bounding \eqref{x04}, for which we write:

\begin{align*}
| \sum_{ \substack{N_2 \leq N_1\\ N_2 \geq 1}} \eqref{x04}| \leq& \sum_{1\leq N_3\leq N_2 \leq N_1} N_1^{-2} \left|\int_{\R^3} \underline f_{N_1, n_1}(x) \nabla g_{N_2, n_2,N_1} \nabla( h^*_{N_3, n_3, N_1}) dx \right| \\
\leq& \sum_{1\leq N_3\leq N_2 \leq N_1} \frac{N_2 N_3}{N_1^{2}} \left|\int_{\R^3} \underline f_{N_1, n_1}(x) \mathfrak B g_{N_2, n_2,N_1} \mathfrak B h^*_{N_3, n_3, N_1} dx \right| \\
\end{align*}
where we denoted $h^*_{N_3, n_3, N_1}=h_{N_3, n_3, N_1}$ when $N_3 >1$ and $h^*_{1, n_3, N_1}=h_{\leq 1, n_3, N_1}$. The result now follows in this case exactly as in \eqref{bound on 3.24R}.
\textbf{Case 2: $N_1 \geq 1$ and $N_2 < 1$.}  In this case, we bound (with $\tilde g=\langle D^\sharp\rangle P_{\leq1} g, \tilde h=\langle D^\sharp\rangle P_{\leq1} h$)

\begin{align*}
| \sum_{ \substack{N_2 \leq N_1}} \eqref{x04}| \leq& \sum_{N_3\leq N_2 \leq N_1} N_1^{-2} \left|\int_{\R^3} \underline f_{N_1, n_1}(x) \nabla g_{N_2, n_2,N_1} \nabla( h_{N_3, n_3, N_1}) dx \right| \\
=& \sum_{N_3\leq N_2<1 \leq N_1} N_1^{-2} \left|\int_{\R^3} \underline f_{N_1, n_1}(x) \mathfrak B \tilde g_{N_2, n_2,N_1} \mathfrak B \tilde h_{N_3, n_3, N_1} dx \right| \\
=& \sum_{N_3\leq N_2<1 \leq N_1} N_1^{-2}\left(\frac{N_3}{N_2}\right)^\epsilon \left|\int_{\R^3} \underline f_{N_1, n_1}(x) N_2^\epsilon \mathfrak B \tilde g_{N_2, n_2,N_1} N_3^{-\epsilon}\mathfrak B \tilde h_{N_3, n_3, N_1} dx \right| \\
\leq& \sum_{N_1\geq 1} N_1^{-2} \int_{\R^3} |\underline f_{N_1, n_1}|  ( \sum_{N_2}N_2^{2\epsilon} |\mathfrak B \tilde g_{N_2, n_2,N_1}|^2 )^{1/2}  ( \sum_{N_3} N_3^{-2\epsilon}|\mathfrak B \tilde h_{N_3, n_3, N_1}|^2 )^{1/2} dx \\
\lesssim& \|\sup_{N}|\underline f_{N, n_1}|\|_{L^{\tilde q}} \sup_{K \geq 1} \|  ( \sum_{N_2\leq 1}N_2^{2\epsilon} |\tilde g_{N_2, n_2,KN_2}|^2 )^{1/2}  \|_{L^{\tilde p}}  \| ( \sum_{N_3} N_3^{-2\epsilon}|\tilde h_{N_3, n_3, KN_3}|^2 )^{1/2} \|_{L^{\tilde r^*}} 
\end{align*}
where we have used Lemma \ref{square and maximal lemma} (part c)) and chose $\tilde r^*$ such that:
$$
\frac{1}{\tilde r^*}=\frac{1}{\tilde r}- \frac{\epsilon}{d}\quad\hbox{ and } \frac{1}{\tilde q}+\frac{1}{\tilde p}+\frac{1}{\tilde r^*}=1.
$$
Using Lemma \ref{square and maximal lemma} and Sobolev embedding we arrive at the needed estimate:
$$
|\sum_{ \substack{N_2 \leq N_1}} \eqref{x04}| \lesssim  \langle n_1 \rangle^{d}\langle n_2 \rangle^{d}\langle n_3 \rangle^{d} \|f\|_{L^{\tilde q}} \|g\|_{L^{\tilde p}} \| h\|_{L^{\tilde r}}.
$$

\textbf{Case 3: $N_1 \leq 1$.} Here we refrain from doing the integration by parts after \eqref{Lambda estimate} and estimate:
\begin{align*}
|\Xi_{N_1, N_2}^{n_1,n_2, n_3}| \lesssim& \sum_{N_2 \leq N_1\leq 1} (\frac{N_2}{N_1})^\epsilon\int_{\R^3} (N_1^\epsilon |f_{N_1,n_1}|)(N_2^{-\epsilon} \sup_{N_1 \geq N_2} |g_{N_2, n_2, N_1}|) (\sup_{N_1 \geq N_2} |h_{\leq N_2, n_3, N_1}|) \, dx \\
\lesssim& \int_{\R^3} (\sum_{N_1\leq 1}N_1^{2\epsilon}|f_{N_1, n_1}|^2)^{1/2}(\sum_{N_2} N_2^{-2\epsilon}\sup_{N_1 \geq N_2} |g_{N_2, n_2, N_1}|^2)^{1/2} (\sup_{N_1 \geq N_2} |h_{\leq N_2, n_3, N_1}|)\, dx\\
\lesssim&  \langle n_1 \rangle^{d}\langle n_2 \rangle^{d}\langle n_3 \rangle^{d} \|{(D^\sharp)}^\epsilon f_{\leq 64}\|_{L^{\tilde q}}\| {(D^\sharp)}^{-\epsilon} g\|_{L^{(\frac{1}{\tilde p} -\frac{\epsilon}{d})^{-1}}}\|h\|_{L^{\tilde r}}\\
\lesssim& \langle n_1 \rangle^{d}\langle n_2 \rangle^{d}\langle n_3 \rangle^{d} \|f_{\leq 64}\|_{L^{\tilde q}}\| g\|_{L^{\tilde p}}\|h\|_{L^{\tilde r}}. \qedhere
\end{align*}
\end{proof}\subsection{Commutators of position and wave operators}\label{commutator subsection}

The purpose of this section is to prove the following theorem.

\begin{theorem}\label{commutator bound}
Assume that $V$ satisfies assumption {\bf H3**} as in Yajima \cite{yajima1995Wkp} (cf.  Remark \ref{prelim:H3remark}.)
Then for any radial function $a(x)$ satisfying $|\nabla a|\lesssim 1$, the wave operator $\Omega$ satisfies: 
\begin{equation}\label{commboundeq}
[a(x), \Omega]: L^p(\R^3) \to L^q(\R^3)\quad \hbox{ for any $1\leq p < q\leq \infty$ with } \frac{1}{p}-\frac{1}{q}< \min(\frac{1}{6}, \frac{1}{p});
\end{equation}
In particular, $[a(x), \Omega]$ is bounded from $L^p$ to $L^{p+\epsilon}$ whenever $1\leq p <\infty$ and $\epsilon$ is sufficiently small. 
\end{theorem}
\begin{remark}
Typical examples for $a(x)$ are $|x|$ and $\langle x \rangle$. In the latter case, the above theorem gives the endpoint case of Corollary 1.6 of \cite{beceanu2011waveop}, where $\Omega$ is proved to be bounded on $\langle x\rangle^{-\beta} L^p(\R^3)$ for $0\leq \beta <1$.\end{remark}
\begin{corollary}
Assume that $V$ satisfies assumption {\bf H3**}. The wave operator $\Omega$ is bounded on weighted spaces $\langle x\rangle^{-\beta} L^p(\R^3)$ for $0\leq \beta \leq 1$.
\end{corollary}
\begin{remark}
Since $[a(x), \Omega]^*=[\Omega^*, \overline{a(x)}]$, the above theorem is equivalent to the boundedness of $[a(x), \Omega^*]$ with the same range of exponents and same conditions on $a$. In fact, it is the latter statement that we prove below.
\end{remark}
 \begin{proof}
By the remark above, the boundedness of $[a(x), \Omega]$ in the $(p,q)$ range specified in \eqref{commboundeq} would follow by interpolation once we show that
\begin{align}
\label{comaim1}[a(x),\Omega^*]: L^1(\R^3) &\to L^p(\R^3)\\
\label{comaim2}[a(x),\Omega^*]: L^q(\R^3) &\to L^\infty(\R^3)
\end{align}
are bounded for any $p\in (1,6/5)$ and any $q>6$. We start with the first estimate \eqref{comaim1}.   For this we need the expansion of the wave operator as
\begin{equation}\label{Omega expression}
\Omega^* f= f -W_1 f +W_2 f - W_3 f +L f
\end{equation}
where $W_1, W_2, W_3 $ and $L$ are bounded operators on $L^p$ for all $1\leq p \leq \infty$ and will be described in more detail below. 

\subsubsection*{Step 1: Bounding $[a(x), W_j]$} Recall from \cite[Equation (1.28)]{yajima1995Wkp} that $W_j$ takes the following form:
\begin{align}
\label{cormoran}
W_j f= \int_{[0,\infty)^{j-1}\times I \times \Sigma^j} \widetilde K_j(t_1, \ldots, t_{j-1}, \tau, \boldsymbol{\omega}) f(\overline x +\rho) dt_1 \ldots dt_{j-1}d\tau  \, d\boldsymbol{\omega} 
\end{align}
where $\Sigma= S^2$ (unit sphere in $\R^3$); $\boldsymbol{\omega} = (\omega_1, \ldots, \omega_j) \in \Sigma^j$;   $I=(-\infty, -2\omega_j(x+t_1\omega_1 +\ldots+t_{j-1}\omega_{j-1})$ is the range of integration of the variable $\tau$; $\overline y= y-2(\omega_j \cdot y)\omega_j$ is an isometry; 
and $\rho= t_1 \overline \omega_1 +\ldots +t_{j-1}\overline{\omega_{j-1}}-\tau \omega_j$. Furthermore, the kernels $\widetilde K_j$ ($j=1,2,3$) are given by:
\begin{equation}\label{explicitK_j}
\widetilde K_j(t_1, \ldots, t_j, \boldsymbol{\omega})=\int_{[0, \infty)^n}e^{-i\sum_{l=1}^j t_ls_l/2}(s_1 \ldots s_j)K_j(s_1\omega_1, \ldots, s_j \omega_j) ds_1 \ldots ds_j
\end{equation}
where 
$$
K_j(k_1, \ldots, k_j)=i^j (2\pi)^{-3j/2}2^{-j}\prod_{l=1}^j\widehat V(k_l-k_{l-1}),\quad \textrm{and} \qquad k_0:=0. 
$$

For instance, 
$$
\widetilde K_1(t_1, \omega_1)=c\int_{[0, \infty)}e^{-it_1s_1/2}s_1\widehat V(s_1\omega_1)ds_1
$$
which can be written by integrating by parts (using the identity $e^{its/2}=\frac{2}{it}\partial_s e^{its/2}$) as:
\begin{align*}
\widetilde K_1 (t_1, \omega_1)=\frac{c'\widehat V(0)}{t^2} -\frac{c'}{t^2}\int_0^\infty e^{-it_1s/2}\partial_s^2(s \widehat V(s\omega_1))ds.
\end{align*}
Since, by Plancherel's theorem, $\widetilde K_1$ is also obviously bounded in $L^2_{t_1, \omega_1}$, we can write that:
\begin{align}\label{shape of K_1}
\langle t_1\rangle^2|\widehat K_1 (t_1, \omega_1)|\lesssim_V  p(t_1, \omega_1)+p'(t_1, \omega_1).
\end{align}
where $p \in L^\infty_{t_1, \omega_1}$ and $p' \in L^2_{t_1, \omega_1}$. A similar calculation for $\widetilde K_j$ shows that if we define
\begin{equation}\label{def of k_j}
k_j(t_1, \ldots, t_j, \omega_1, \ldots, \omega_j):=\left(\prod_{l=1}^j\langle t_l\rangle^2\right)|\widehat K_j (t_1,\ldots, t_j, \omega_1, \ldots, \omega_j)|,
\end{equation}
then $k_j(t_1, \ldots, t_j, \boldsymbol{\omega})$ can be bounded by a finite linear combination of products of the form:
\begin{align}\label{shape of k_j}
\prod_{l=1}^j p_l(t_l, \omega_l)
\end{align}
where each $p_l(t_l, \omega_l)$ is either in $L^\infty_{t_l, \omega_l}$ or in $L^2_{t_l, \omega_l}$. Consequently, from now on we will assume without any loss of generality that:
\begin{equation}\label{estimate on tildeK_j}
\begin{cases}
|\widetilde K_j (t_1, \ldots, t_j, \boldsymbol{\omega})|=\prod_{l=1}^j\langle t_l\rangle^{-2} k_j (t_1, \ldots, t_j,\boldsymbol{\omega})\\
k_j(t_1, \ldots, t_j, \boldsymbol{\omega}) \lesssim p_1(t_1, \omega_1)\ldots p_j(t_j, \omega_j); \qquad \quad p_i(t_i, \omega_i)\in L^\infty_{t_i, \omega_i}+L^2_{t_i, \omega_i}.
\end{cases}
\end{equation}
\begin{remark}\label{estimateremark}
It follows from the integration by parts argument leading to \eqref{estimate on tildeK_j}, that the highest order term in $K_1$ is exactly $C\widehat V(0)/t^2$. This explains the failure of the boundedness of \eqref{commboundeq} when $p=q$.
\end{remark}
Changing variables in \eqref{cormoran} by $\tau \mapsto t_j=-\tau -2\omega_j(x+t_1\omega_1+\ldots +t_{j-1}\omega_{j-1})$, we get that:
$$
W_j f(x)=\int_{[0,\infty)^j \times \Sigma^j} \widetilde K_j(t_1, \ldots, t_{j-1}, -t_j-\sigma, \boldsymbol{\omega}) f(x+\tilde \rho) \,d\mathbf{t} \, d \boldsymbol{\omega} 
$$
where $\sigma=2\omega_j(x+t_1\omega_1 +\ldots+t_{j-1}\omega_{j-1})$, $\mathbf{t} = (t_1, \ldots, t_j)$, and $\tilde \rho=t_1\omega_1 +\ldots+t_j\omega_j$. Consequently, we have the following expression for $[a(x), W_j]$:
\begin{equation}\label{commeq1}
[a(x), W_j]=\int_{[0,\infty)^j \times \Sigma^j} [a(x)-a(x+\tilde \rho)]\widetilde K_j(t_1, \ldots, t_{j-1}, -t_j-\sigma, \boldsymbol{\omega}) f(x+\tilde \rho) \,d\mathbf{t} \, d \boldsymbol{\omega} 
\end{equation}

Now, notice that since $a$ is radial and $z\mapsto z-2(z \cdot \omega_j)\omega_j$ is an isometry (reflection with respect to a plane), we have that $a(x)=a\left(x-2(x \cdot \omega_j)\omega_j\right)$ which gives along with the 
bound on $\nabla a$ that:
\begin{align*}
|a(x)-a(x+\tilde \rho)|&\lesssim |\tilde \rho+2(x \cdot \omega_j)\omega_j|=|\tilde \rho +2(x \cdot \omega_j) \omega_j -2[\left(\tilde \rho +2(x \cdot \omega_j) \omega_j\right) \cdot \omega_j]\omega_j|\\
&=|t_1 \omega_1 +\ldots +t_{j-1}\omega_{j-1} -(t_j+\sigma) \omega_j|\leq |t_1|+\ldots+|t_{j-1}|+|t_j+\sigma|.
\end{align*}

We remark that it is important here to get a bound in terms of $|t_j+\sigma|$ rather than $|t_j|$, since the kernel $\widetilde K_j(t_1, \ldots, t_{j-1}, -t_j-\sigma, \boldsymbol{\omega})$ in \eqref{commeq1} does 
not decay in $t_j$, but in $t_j+\sigma$. This leads to the following remark: 

\begin{remark} \label{commutatorunboundedremark} 
From a technical point of view, if $a$ was not assumed to be radial (e.g. $a(x)=x$), then $|a(x)-a(x+\tilde \rho)|$ can only be bounded by $|\tilde \rho|$ which would not give a bounded contribution from $L^p$ to $L^q$ even for the $[a(x),W_1]$. In fact, an asymptotic expansion of $\widetilde K_1$ in \eqref{explicitK_j} shows that its first order term is $t^{-2}$, for which the commutator $[x, \Omega]$ does not seem to be bounded from $L^p \to L^q$ whenever $q\geq p$, unless $f$ is assumed to lie in some weighted space.

This failure makes sense physically: $\Omega$ maps the initial data of the Schr\"odinger equation with operator $H$ to its (free) scattering data at time $+\infty$. From the quantum-classical correspondence, one can think of the solution at time zero as a particle coming from far away towards the potential; it can be deflected for positive time and emerge from the interaction  behaving like a free particle that started from a different direction. Thus $\Omega$ has no obligation to preserve directions, and hence can only commute well with radial localization operators.
\end{remark}

As a result, we have:
\begin{align}\label{numfull}
\left|[a(x), W_j] f\right|\leq&\int_{[0,\infty)^j \times \Sigma^j} \left(|t_1|+\ldots+|t_{j-1}|+|t_j +\sigma|\right)|\widetilde K_j(t_1, 
\ldots, t_j+\sigma)||f(x+\tilde \rho)| \,d\mathbf{t} \, d \boldsymbol{\omega} \nonumber \\
\lesssim&\int_{[0,\infty)^j \times \Sigma^j} \frac{|t_1|+\ldots+|t_{j-1}|+|t_j +\sigma|}{\langle t_1\rangle^2\ldots \langle t_{j-1}\rangle^2 \langle t_j+\sigma \rangle^2}k_j(t_1, \ldots, t_j+\sigma, \boldsymbol{\omega})|f(x+\tilde \rho)| \,d\mathbf{t} \, d \boldsymbol{\omega} .
\end{align}
We will only prove the bounds \eqref{comaim1} and $\eqref{comaim2}$ for the operator
\begin{align}
\label{numt_j}\mathfrak G_1 f:=\int_{[0,\infty)^j \times \Sigma^j} \frac{|t_j +\sigma|}{\langle t_1\rangle^2\ldots \langle t_{j-1}\rangle^2 \langle t_j+\sigma \rangle^2}k_j(t_1, \ldots, t_j+\sigma, \boldsymbol{\omega}) \,|f(x+\tilde \rho)| \, d\mathbf{t}\, d\boldsymbol{\omega}
\end{align}
as the others are similar. We assume without loss of generality that $f\geq 0$, and argue by duality. For any non-negative test function $h$,
\begin{align}
\nonumber \langle \mathfrak G_1f, h\rangle =&\int_{\R_x^3}\int_{[0,\infty)^j \times \Sigma^j} \frac{1}{\langle t_1\rangle^2\ldots \langle t_{j-1}\rangle^2 \langle t_j+\sigma \rangle}k_j(t_1, \ldots, t_j+\sigma, \boldsymbol{\omega}) f(x+\tilde \rho) h(x) \,d\mathbf{t} \, d \boldsymbol{\omega} \,dx\\
&\lesssim [\int_{\R_x^3}\int_{[0,\infty)^{j-1}\times \R \times\Sigma^j} f( x +\tilde \rho) \left(\prod_{l=1}^{j-1}p_l (t_l, \omega_l)^p\langle t_l\rangle^{-p}\right) \langle t_j+\sigma \rangle^{-p} |t_j|^{2-2p}p_j(t_j+\sigma, \omega_j)^p \,d\mathbf{t} \, d \boldsymbol{\omega}  \, dx]^{1/p}\label{holder0.1} \\
&\times [\int_{\R_x^3}h( x)^{p'} \int_{[0,\infty)^{j-1}\times \R\times \Sigma^j} f(x +\tilde \rho) (\prod_{l=1}^{j-1}\langle t_l\rangle^{-p'})  |t_j|^2\,d\mathbf{t} \, d \boldsymbol{\omega}  \, dx]^{1/p'}. \label{holder0.2}
\end{align}
We start by estimating \eqref{holder0.1}. Recall that $\sigma=2\omega_j(x+t_1\omega_1 +\ldots+t_{j-1}\omega_{j-1})$ and $\tilde \rho=t_1\omega_1 +\ldots+t_j\omega_j$, so
\begin{align*}
\eqref{holder0.1}^p &=  \int_{[0,\infty)^{j-1}\times\Sigma^j}\bigg(\prod_{l=1}^{j-1}p_l(t_l, \omega_l)^p\langle t_l\rangle^{-p}\bigg)\int_{\R_{t_j}} \int_{\R_x^3}f(x+\tilde \rho)  \frac{|t_j|^{2-2p}}{\langle t_j +\sigma\rangle^{p}} p_j(t_j+\sigma, \omega_j)^p \, dx\,d\mathbf{t} \, d \boldsymbol{\omega} \\
&=\int_{[0,\infty)^{j-1}\times\Sigma^j}\bigg(\prod_{l=1}^{j-1}p_l(t_l, \omega_l)^p\langle t_l\rangle^{-p}\bigg)\int_{\R_{t_j}} \int_{\R_y^3}f(y)  \frac{|t_j|^{2-2p}}{\langle -t_j +2\omega_j .y\rangle^{p}} p_j(-t_j+2\omega_j. y, \omega_j)^p \, dy\,d\mathbf{t} \, d \boldsymbol{\omega}\\
&\lesssim \int_{\mathbb{S}^2_{\omega_j}}\int_{\R_{t_j}} \int_{\R_y^3}f(y)  \frac{|t_j|^{2-2p}}{\langle -t_j +2\omega_j .y\rangle^{p}} p_j(-t_j+2\omega_j. y, \omega_j)^p \, dy\,dt_j \, d \omega_j;\\
\end{align*}
since $p<2$ and $p_l \in L^\infty_{t_l, \omega_l}+L^2_{t_l, \omega_l}$. As a result, we get that:
\begin{align*}
\eqref{holder0.1}^p &\lesssim \|f\|_{L^1}\sup_{y \in \R^3}\int_{\mathbb{S}^2_{\omega_j}}\int_{\R_{t_j}} \frac{|t_j|^{2-2p}}{\langle -t_j +2\omega_j .y\rangle^{p}} p_j(-t_j+2\omega_j. y, \omega_j)^p \,dt_j \, d \omega_j\\
&\lesssim \|f\|_{L^1}\sup_{y \in \R^3}\int_{\mathbb{S}^2_{\omega_j}}\left[\int_{-1}^1 |t_j|^{2-2p}p_j(-t_j+2\omega_j. y, \omega_j)^p \, dy\,dt_j+\int_{|t_j|\geq 1}\frac{p_j(-t_j+2\omega_j. y, \omega_j)^p}{\langle -t_j +2\omega_j .y\rangle^{p}} \,dt_j\right]d\omega_j\\
&\lesssim \|f\|_{L^1}
\end{align*}
if $p_j\in L^\infty_{t_j, \omega_j}+L^2_{t_j, \omega_j}$ and $1<p<\frac{6}{5}$ thanks to H\"older's inequality. Consequently, we have that
$$
\eqref{holder0.1}\lesssim \|f\|_{L^1}^{1/p}.
$$

To estimate $\eqref{holder0.2}$, we observe that
\begin{align*}
&\sup_{x\in \R^3}\int_{[0,\infty)^{j-1}\times \R\times \Sigma^j} f(x +\tilde \rho) \left(\prod_{l=1}^{j-1}\langle t_l\rangle^{-p'}\right)  |t_j|^2 \,d\mathbf{t} \, d \boldsymbol{\omega} \\
&\sup_{x\in \R^3}\int_{[0,\infty)^{j-1}\times \Sigma^{j-1}}\int_{\R_y^3} f(x +t_1 \omega_1 +\ldots +t_{j-1}\omega_{j-1} +y) \left(\prod_{l=1}^{j-1}\langle t_l\rangle^{-p'}\right)  dt_1 \ldots dt_{j-1}\, d\omega_1 \ldots d\omega_{j-1} \, dy\, dx \\
&\lesssim \|f\|_{L^1},
\end{align*}
from which we can conclude  
$$
\eqref{holder0.2}\lesssim \|h\|_{L^{p'}}\|f\|_{L^1}^{1/p'}.
$$
This shows that $[a(x), W_j]:L^1 \to L^p$ for all $p\in (1,6/5)$ and hence gives a satisfactory contribution to \eqref{comaim1}. 

Next we show that this commutator is also bounded from $L^q\to L^\infty$ which verifies \eqref{comaim2} and concludes Step 1. As before, we only write the estimates for the contribution of \eqref{numt_j} to \eqref{numfull} (the other terms are treated similarly).
For this, we notice that for every $x\in \R^3$, if $1<q'<6/5$ (or equivalently $q>6$), then
\begin{align*} 
\mathfrak G_1 f(x) & \leq [\int_{[0,\infty)^{j}\times \Sigma^j}  f^q(x +t_1 \omega_1 +\ldots +\overbrace{t_{j}\omega_{j}}^{y}) \,\overbrace{t_j^2 dt_j\,d\omega_j}^{dy} \,\frac{dt_1}{\langle t_1\rangle^{q}} \ldots \frac{dt_{j-1}}{\langle t_{j-1}\rangle^{q}} \, d\omega_1 \ldots d\omega_{j-1}]^{1/q}\\
&\qquad \times [\int_{[0,\infty)^{j} \times \Sigma^j}  \frac{p_1(t_1, \omega_1)^{q'}}{\langle t_1\rangle^{q'}} \ldots \frac{p_{j-1}(t_{j-1})^{q'}}{\langle t_{j-1}\rangle^{q'}}\frac{|t_j|^{2-2q'}p_j(t_j+\sigma, \omega_j)^{q'}}{\langle t_j+\sigma \rangle^{q'}} d\boldsymbol{t} \, d\boldsymbol{\omega}]^{1/q'} \\
& \lesssim \|f\|_{L^q},
\end{align*}
by arguing exactly as when we estimated the contribution of \eqref{holder0.1} above.
\subsubsection*{Step 2: Bounding $[a(x), L]$} In order to prove the desired estimates for the integral operator $L$, we will need more information about its kernel $L(x,y)$. These are included in the following lemma:
\begin{lemma}[\cite{yajima1995Wkp}]\label{Lkernel}
The kernel $L(x,y)$ of the integral operator $L$ can be written as:
$$
L(x,y)=\frac{L_+(x,y)-L_-(x,y)}{i\pi}
$$
where $L_\pm(x,y)=\sum_{l=0}^3 Z_{\pm, l}(x,y)$ and $Z_{\pm}$ satisfy the following bounds:
\begin{align}
\label{Zest1}|Z_{\pm, l}(x,y)|\lesssim \langle |x|-|y|\rangle^{-3}\langle x\rangle^{-1}\langle y \rangle^{-1}\quad \hbox{for }l=0,2,3.\\
\label{Zest2}|Z_{+, 1}(x,y)-Z_{-,1}(x,y)|\lesssim \langle |x|-|y|\rangle^{-2}\min\left(\langle x\rangle^{-2},\langle y \rangle^{-2}\right)
\end{align}
\end{lemma}
\begin{proof}
All these estimate can be found in \cite{yajima1995Wkp} (particularly in displays (3.10), (3.27), (3.28) and Propositions 3.2, 3.3, and 3.4).
\end{proof}

The operator $[a(x), L]$ is an integral operator with kernel
$$
K(x,y):=[a(x)-a(y)]L(x,y).
$$
Since $a$ is radial with bounded gradient $|a(x)-a(y)|=|a(|x|\frac{y}{|y|})-a(y)|\lesssim ||x|-|y||$ and hence by Lemma \ref{Lkernel} we can bound:
$$
|K(x,y)|\lesssim K_1(x,y)+K_2(x,y)
$$
where $K_1$ and $K_2$ are the non-negative functions
\begin{equation}\label{estimates on K_j}
\begin{split}
K_1(x,y)=&\langle |x|-|y|\rangle^{-2}\langle x\rangle^{-1}\langle y \rangle^{-1}\\
K_2(x,y)=&\langle |x|-|y|\rangle^{-1}\min\left(\langle x\rangle^{-2},\,\langle y \rangle^{-2}\right).
\end{split}
\end{equation}

Our aim is to show that the integral operators $T_j: f(x)\mapsto \int_{\R^3}K_j(x,y) f(y) dy$ are bounded from $L^1(\R^3) \to L^p(\R^3)$ and from $L^q(\R^3)$ to $L^\infty(\R^3)$ for $1<p<6/5$ and $6<q <\infty$. Since $K(x,y)$ is symmetric, $T_j$ is formally self-adjoint and it is enough to show that it is bounded from $L^1\to L^p$ for all $1<p<6/5$. 

By Minkowski's inequality, this follows once we show that for $j=1,2$:
\begin{equation}
\sup_{y} \|K_j(\cdot,y)\|_{L^p(\R^3)} \lesssim 1\quad \hbox{whenever }1<p<\frac{6}{5},
\end{equation}
which is easily verified thanks to \eqref{estimates on K_j}.
\end{proof}

\subsection{Derivatives of multilinear estimates}\label{derivatives of M section}
In this section, we develop an understanding of how to manage estimates involving the derivatives of $\mathcal M$. As we discussed in Section \ref{intro:spacetimeres}, this understanding informs at a  fundamental level the definition of space resonance in the distorted regime.  Moreover, derivatives of $\mathcal{M}$ arise when one attempts to prove boundedness of the multilinear operators $T$ and $\Lambda$ (see \eqref{CMC0} and \eqref{CME2}) in weighted Sobolev spaces.  In particular, we will see that in PDE applications, it is often important to be able to bound quantitates like $x T(f,g)$ in $L^p$, say.  Since 
\[ x T(f,g) =\Omega^* ({\FF^{\sharp}}^{-1} (-i\nabla_\xi) \FF^{\sharp} ) \Omega T(f,g),\]
 and $\Omega$ is bounded on $L^p$, it suffices to study 
\beq {\mathcal F^{\sharp}}^{-1}\partial_{\xi_3} \int \!\!\! \int m(\xi_1, \xi_2, \xi_3) f^\sharp(\xi_1) g^{\sharp}(\xi_2) \mathcal M(\xi_1, \xi_2, \xi_3) \, d\xi_1 \,d\xi_2. \label{wnorm:dxiLambda} \eeq
Naturally, this necessitates having a distributional interpretation for $\partial_{\xi_3} \mathcal{M}$.

To develop an intuition for what's going on, let us look at the case $V=0$. 
There $\MM=\delta(\xi_1+\xi_2+\xi_3)$ and the above expression is equal to 
\begin{align*}
{\mathcal F}^{-1}\left\{- \int m(-\xi_2-\xi_3, \xi_2, \xi_3) \partial_{\xi}f^\sharp(-\xi_3-\xi_2) g^{\sharp}(\xi_2)\,d\xi_2\right.\\
\left.\quad \quad+\int \partial_{\xi_3}m(-\xi_2-\xi_3, \xi_2, \xi_3) f^\sharp(-\xi_3-\xi_2) g^{\sharp}(\xi_2)\,d\xi_2\right\}.
\end{align*}
These are standard convolution integrals that one can estimate in $L^p$ in terms of Lebesgue  norms of $f, g$ and $xf$, for appropriate multipliers $m$. A moment's thought shows that this is possible due to the identity
$$
\partial_{\xi_3}\delta(\xi_1+\xi_2+\xi_3)=\partial_{\xi_1}\delta(\xi_1+\xi_2+\xi_3)
$$
which enables us to integrate by parts in $\xi_1$.  Consequently, we might hope that for $V\not \equiv 0$, 
\begin{equation}\label{hoperelation}
\partial_{\xi_3}\MM(\xi_1, \xi_2, \xi_3)=\partial_{\xi_1}\MM(\xi_1,\xi_2, \xi_3)+\textrm{remainder},
\end{equation}
where this remainder is a distribution of order $0$.
This is equivalent to saying that the vector field $\partial_{\xi_3}-\partial_{\xi_1}$ preserves the order of the distribution $\MM$.
To test the validity chances of \eqref{hoperelation}, we inspect the case $m=1$, for which  \eqref{hoperelation} would predict that
$$
\Omega x \Omega^* (fg)=(\Omega x \Omega^*f) g+\textrm{remainder}.
$$
The relation 
$$\Omega x \Omega^* (fg)=xfg+[\Omega, x]\Omega^*(fg)=(\Omega x \Omega^*f)g+\left([\Omega, x] \Omega^*f\right) g+[\Omega, x]\Omega^*(fg),$$
shows that the validity of the previous line is tied to the boundedness of the commutator $[x, \Omega]$ on $L^p$. Remark \ref{commutatorunboundedremark} tells us, however, that $[x,\Omega]$ \emph{cannot} be bounded since $x$ is not radial.   It is then more reasonable to expect that 
$$
\Omega |x| \Omega^* (fg)=(\Omega |x| \Omega^*f) g+\textrm{remainder},
$$
with the remainder bounded in $L^p$ in terms of Lebesgue norms of $f$ and $g$\footnote{It seems important for our future work  (and much cleaner for our current application) that estimates on this remainder do not involve weights falling on $f$ or $g$.}. In other words, we seek to relate $|\partial_{\xi_3}|\MM:=|\nabla_{\xi_3}| \MM$ and $|\partial_{\xi_1}|\MM=|\nabla_{\xi_1}|\MM$, rather than something of the form of \eqref{hoperelation}. This strategy turns out to be the right one, and we are able to exhibit precisely such a relation, both qualitatively in \eqref{derivative of M2}, and quantitatively in Theorem \ref{colibri}.

\begin{remark}
In what follows, we will use the summation convention that repeated indices are summed over from 1 to 3. For instance, $|\nabla_x|=R^\alpha \partial_{x^\alpha}$ where $R^\alpha=\frac{\partial_{x^\alpha}}{|\nabla|}$ is the Euclidean Riesz transform.
\end{remark}
\subsubsection{Understanding $\partial_\xi \mathcal{M}$}\label{understanding M} Let us start by looking at \eqref{wnorm:dxiLambda} with $m=1$. By duality, integrating \eqref{wnorm:dxiLambda} against a test function $h$, we get 
\beq -\int_{\R^3\times \R^3\times \R^3} f^\sharp(\xi_1) g^{\sharp}(\xi_2) \partial_{\xi_3^\alpha}h^\sharp (\xi_3) \mathcal M(\xi_1, \xi_2, \xi_3) d\xi_1 \,d\xi_2 \,d\xi_3.  \label{wnorms:prelimexpression}\eeq
From the identity $x^\alpha f=\Omega^* ({\FF^{\sharp}}^{-1} (-i\partial_{\xi^\alpha}) \FF^{\sharp} ) \Omega f$ we have (dropping for convenience the constant $-i$)
\begin{align*}
-\eqref{wnorms:prelimexpression} &= \int_{\R^3} f(x) g(x) \Omega|x|\Omega^*(\Omega\frac{x^\alpha}{|x|}\Omega^*) h (x) \, dx\\
&=\int_{\R^3} |x| f(x) g(x)(\Omega\frac{x^\alpha}{|x|}\Omega^*) h (x) \, dx +\int_{\R^3} f(x) g(x) [\Omega, |x|]\Omega^* (\Omega\frac{x^\alpha}{|x|}\Omega^*) h (x) \, dx\\
&=\int_{\R^3} (\Omega|x|\Omega^*) f(x) g(x)(\Omega\frac{x^\alpha}{|x|}\Omega^*) h (x) \, dx +\int_{\R^3} \left([|x|, \Omega]\Omega^* f(x)\right) g(x)(\Omega\frac{x^\alpha}{|x|}\Omega^*) h (x)\, dx\\
&\qquad -\int_{\R^3} f(x) g(x) [|x|,\Omega]\Omega^* (\Omega\frac{x^\alpha}{|x|}\Omega^*) h (x) \, dx.
\end{align*}

Define the operators  
\[ \RR^{\alpha}:=\Omega\frac{x^\alpha}{|x|}\Omega^*, \qquad \mathcal E=[|x|,\Omega]\Omega^*,\]
and their Fourier space manifestations
\begin{equation}\label{mcandnomc}
R^\alpha:=\FF^\sharp \RR^\alpha {\FF^\sharp}^{-1}, \qquad E=\FF^\sharp \mathcal E {\FF^\sharp}^{-1}.
\end{equation}
so that $\FF^\sharp \RR^\alpha f= R^\alpha \FF^\sharp f$ and similarly for $\mathcal E$ and $E$. Notice that since $\FF^\sharp \Omega=\mathcal F$, $R^\alpha$ is nothing but the Euclidean Riesz transform $\frac{\partial_{\xi^\alpha}}{|\nabla|}$. As a result, we may write
\begin{align}
-\eqref{wnorms:prelimexpression} & =\int_{\R^3} \RR^\beta(\Omega x_\beta \Omega^*)f(x) g(x)\RR^\alpha h (x) dx \nonumber \\
&\qquad+\int_{\R^3} \mathcal E f(x) g(x)\RR^\alpha h (x)\, dx
-\int_{\R^3} f(x) g(x) \mathcal E \RR^\alpha h (x)\, dx\label{distsensephys}\\
&=\int_{\R^3} R_{\xi_1}^\beta\partial_{\xi_1^\beta}f^\sharp(\xi_1) g^\sharp(\xi_2) R^\alpha_{\xi_3} h^\sharp (\xi_3) \MM(\xi_1, \xi_2, \xi_3) d\xi_1 \, d\xi_2\, d\xi_3 \label{distsense1}\\
&\qquad +\int_{\R^3}  Ef^\sharp(\xi_1) g^\sharp(\xi_2)R^\alpha_{\xi_3} h^\sharp (\xi_3) \MM(\xi_1, \xi_2, \xi_3)d\xi_1 \, d\xi_2 d\xi_3 \label{distsense2}\\
&\qquad-\int_{\R^3} f^\sharp(\xi_1) g^\sharp(\xi_2)E R^\alpha_{\xi_3} h^\sharp (\xi_3) \MM(\xi_1, \xi_2, \xi_3) \, d\xi_1 \, d\xi_2 \, d\xi_3 \label{distsense3},
\end{align}
In total, this gives the following formal expression for the distribution $\partial_{\xi_3^\alpha}\MM$:
\begin{align}\label{derivative of M}
\partial_{\xi_3^\alpha}\MM=&R_{\xi_3}^\alpha R_{\xi_1}^\beta \partial_{\xi_1^\beta} \MM- R^\alpha_{\xi_3}E_{\xi_1}^*\MM+R^\alpha_{\xi_3}E^*_{\xi_3}\MM\qquad \textrm{or,}\\
|\nabla_{\xi_3}|\MM=& |\nabla_{\xi_1}| \MM- E_{\xi_1}^*\MM+E^*_{\xi_3}\MM\label{derivative of M2}
\end{align}
where we use the notation $E_{\xi_1}$ to indicate that $E$ is being applied in the $\xi_1$ variable fixing $\xi_2, \xi_3$.  Recall that $R^\alpha$ is self-adjoint and commutes with differentiation. 
\begin{remark}
All the entities in  \eqref{derivative of M} and \eqref{derivative of M2} can be interpreted as distributions on $\R^9$ as follows: Clearly, one can make sense of them when applied to test functions in $\mathcal S(\R^9)$ that are linear combinations of tensored functions (i.e. of the form $\phi_1(\xi_1) \phi_2(\xi_2) \phi_3(\xi_3)$). In fact this is explicitly done in \eqref{distsense1}, \eqref{distsense2}, and \eqref{distsense3} and the physical space manifestation in \eqref{distsensephys}. Since this gives a continuous linear functional on a dense subset of $\mathcal S(\R^9)$, the relations \eqref{derivative of M} and \eqref{derivative of M2} are well-defined in $\mathcal S^\prime (\R^9)$ by density.
\end{remark}

\subsubsection{Multilinear estimates} 
After this qualitative formulation, we are now ready to obtain quantitative estimates on \eqref{wnorm:dxiLambda} in $L^{r^\prime}$.  By duality, integrating \eqref{wnorm:dxiLambda} against a test function $h$, we get:
\begin{equation}\label{dual0.1}
-\int_{\R^9} m(\xi_1, \xi_2, \xi_3) f^\sharp(\xi_1) g^{\sharp}(\xi_2) \partial_{\xi_3^\alpha}h^\sharp (\xi_3) \mathcal M(\xi_1, \xi_2, \xi_3) d\xi_1 \,d\xi_2 \,d\xi_3,
\end{equation}
which is equal to
\begin{align*}
&\int_{\R^9} m(\xi_1, \xi_2, \xi_3) f^\sharp(\xi_1) g^{\sharp}(\xi_2) h^\sharp (\xi_3) \partial_{\xi_3^\alpha}\mathcal M(\xi_1, \xi_2, \xi_3) d\xi_1 \,d\xi_2 \,d\xi_3\\
&\qquad+\int_{\R^9} \partial_{\xi_3^\alpha}m(\xi_1, \xi_2, \xi_3) f^\sharp(\xi_1) g^{\sharp}(\xi_2) h^\sharp (\xi_3) \mathcal M(\xi_1, \xi_2, \xi_3) d\xi_1 \,d\xi_2 \,d\xi_3\\
&=-\int_{\R^9} R_{\xi_3}^\alpha R_{\xi_1}^\beta\partial_{\xi_1^\beta}\left[m(\xi_1, \xi_2, \xi_3) f^\sharp(\xi_1) g^{\sharp}(\xi_2) h^\sharp (\xi_3)\right] \mathcal M(\xi_1, \xi_2, \xi_3) d\xi_1 \,d\xi_2 \,d\xi_3\\
&\qquad-\int_{\R^9} R^\alpha_{\xi_3}E_{\xi_1}\left[m(\xi_1, \xi_2, \xi_3) f^\sharp(\xi_1) g^{\sharp}(\xi_2) h^\sharp (\xi_3)\right] \mathcal M(\xi_1, \xi_2, \xi_3)\, d\xi_1 \,d\xi_2 \,d\xi_3\\
&\qquad+\int_{\R^9} E_{\xi_3}R^\alpha_{\xi_3}\left[m(\xi_1, \xi_2, \xi_3) f^\sharp(\xi_1) g^{\sharp}(\xi_2) h^\sharp (\xi_3)\right] \mathcal M(\xi_1, \xi_2, \xi_3)\, d\xi_1 \,d\xi_2 \,d\xi_3\\
&\qquad+\int_{\R^9} \partial_{\xi_3^\alpha}m(\xi_1, \xi_2, \xi_3) f^\sharp(\xi_1) g^{\sharp}(\xi_2) h^\sharp (\xi_3) \mathcal M(\xi_1, \xi_2, \xi_3)\, d\xi_1 \,d\xi_2 \,d\xi_3.
\end{align*}
With some rearranging, this becomes
\begin{align}
\eqref{dual0.1}&=\int_{\R^9} \left[\partial_{\xi_3^\alpha}m\,f^\sharp(\xi_1) h^\sharp (\xi_3)-R_{\xi_3}^\alpha R_{\xi_1}^\beta \left(\partial_{\xi_1^\beta}m\, f^\sharp(\xi_1) h^\sharp (\xi_3)\right)\right] g^{\sharp}(\xi_2) \mathcal M(\xi_1, \xi_2, \xi_3) d\xi_1 \,d\xi_2 \,d\xi_3\label{der on m}\\
&\qquad-\int_{\R^9} R_{\xi_3}^\alpha R_{\xi_1}^\beta \left[ m(\xi_1, \xi_2, \xi_3) \partial_{\xi_1^\beta} f^\sharp(\xi_1) g^{\sharp}(\xi_2) h^\sharp (\xi_3) \right] \mathcal M(\xi_1, \xi_2, \xi_3) d\xi_1 \,d\xi_2 \,d\xi_3\label{der on f}\\
&\qquad-\int_{\R^9} R^\alpha_{\xi_3}E_{\xi_1}\left[m(\xi_1, \xi_2, \xi_3) f^\sharp(\xi_1) g^{\sharp}(\xi_2) h^\sharp (\xi_3)\right] \mathcal M(\xi_1, \xi_2, \xi_3) d\xi_1 \,d\xi_2 \,d\xi_3\label{rem1}\\
&\qquad+\int_{\R^9} E_{\xi_3}R^\alpha_{\xi_3}\left[m(\xi_1, \xi_2, \xi_3) f^\sharp(\xi_1) g^{\sharp}(\xi_2) h^\sharp (\xi_3)\right] \mathcal M(\xi_1, \xi_2, \xi_3) d\xi_1 \,d\xi_2 \,d\xi_3.\label{rem2}
\end{align}
This suggests looking at the following operators
\begin{subequations}
\begin{align}
\label{platypus1}
& \Lambda_1: (f,g,h) \mapsto \int_{\mathbb{R}^9} m(\xi_1, \xi_2, \xi_3) f^\sharp(\xi_1) g^\sharp(\xi_2) h^\sharp(\xi_3) R_{\xi^3}^\alpha R_{\xi_1}^\beta \mathcal{M}(\xi_1,\xi_2,\xi_3) \\
\label{platypus2}
& \Lambda_2: (f,g,h) \mapsto \int_{\mathbb{R}^9} m(\xi_1, \xi_2, \xi_3) f^\sharp(\xi_1) g^\sharp(\xi_2) h^\sharp(\xi_3) E_{\xi_1}^* R_{\xi^3}^\alpha \mathcal{M}(\xi_1,\xi_2,\xi_3) \\
\label{platypus3}
& \Lambda_3: (f,g,h) \mapsto \int_{\mathbb{R}^9} m(\xi_1, \xi_2, \xi_3) f^\sharp(\xi_1) g^\sharp(\xi_2) h^\sharp(\xi_3) R_{\xi^3}^\alpha E_{\xi_3}^* \mathcal{M}(\xi_1,\xi_2,\xi_3) 
\end{align}
\end{subequations}

\begin{theorem} \label{colibri}
Suppose that $V$ is a real potential satisfying {\bf H1}, {\bf H2}, {\bf H3*}, and {\bf H3**}. Let 
$m$ be a symbol of the form 
\begin{equation}\label{bettermult}
m(\xi_1, \xi_2, \xi_3)=\langle |\xi_1|+|\xi_2|+|\xi_3|\rangle^{-\delta}m_0(\xi_1, \xi_2, \xi_3),
\end{equation} 
for a Coifman--Meyer symbol $m_0$ and some $\delta> 0$. For any $q,p,r,\tilde{q},\tilde{p},\tilde{r}$ in $(1,\infty)$ such that
$$
\frac{1}{p} + \frac{1}{q} + \frac{1}{r} = 1 \quad \mbox{and} \quad \frac{1}{\tilde p} + \frac{1}{\tilde q} + \frac{1}{\tilde r} > 1,
$$
we have for all $i = 1,2,3$
$$
|\Lambda_i(f,g,h)| \lesssim \|f\|_p \|g\|_q \|h\|_r + \|f\|_{\widetilde{p}} \|g\|_{\widetilde{q}} \|h\|_{\widetilde{r}}.
$$
\end{theorem}

\begin{proof}
We only bound $\Lambda_1$, the argument for $\Lambda_2$ and $\Lambda_3$ being similar. We first split $f$, $g$ and $h$ into their distorted Littlewood-Paley pieces and write:
\begin{equation*}
\eqref{platypus1}=\sum_{N_1, N_2, N_3 \in 2^\Z}\int_{\R^9} R_{\xi_3}^\alpha R_{\xi_1}^\beta  \left[ m(\xi_1, \xi_2, \xi_3) f_{N_1}^\sharp(\xi_1) g_{N_2}^{\sharp}(\xi_2) h_{N_3}^\sharp (\xi_3) \right] 
\mathcal M(\xi_1, \xi_2, \xi_3) d\xi_1 \,d\xi_2 \,d\xi_3
\end{equation*}

We only present the calculation for the case when $N_1\geq N_2 \geq N_3$. The other cases can be treated analogously.  Expanding $m_0$ as before, we get
\begin{equation*}
m_0(\xi_1, \xi_2, \xi_3) \psi(\frac{\xi_1}{2N_1})\psi(\frac{\xi_2}{2N_1})\psi(\frac{\xi_3}{2N_1})=\sum_{n_1, n_2, n_3 \in \Z^3} a^{N_1}(n_1, n_2, n_3) e^{\frac{2\pi i}{K N_1}(n_1\cdot \xi_1+n_2 \cdot \xi_2 +n_3\cdot \xi_3)},
\end{equation*}
where $|a^{N_1}(n_1, n_2, n_3)|\lesssim (1+|n_1|+|n_2|+|n_3|)^{-100}$ uniformly in $N_1$. Consequently, the contribution of the case $N_1 \geq N_2 \geq N_3$ is bounded by
\begin{multline*}
\sum_{n_1,n_2, n_3} (1+|n_1|+|n_2|+|n_3|)^{-100}\\
\times \sum_{N_2 \leq N_1} \bigg|\int_{\R^9}  R_{\xi_1}^\beta [e^{\frac{2\pi i\, n_1\cdot \xi_1}{K N_1}} \langle \xi_1 \rangle^{-\delta}f_{ N_1}^\sharp(\xi_1)] (e^{\frac{2\pi i\, n_2\cdot \xi_2}{K N_1}}g_{N_2}^{\sharp}(\xi_2)) R_{\xi_3}^\alpha (e^{\frac{2\pi i\, n_3\cdot \xi_3}{K N_1}}h_{\leq N_2}^\sharp (\xi_3))\\
\mathcal M(\xi_1, \xi_2, \xi_3) \,d\xi_1 \,d\xi_2 \,d\xi_3\bigg|.
\end{multline*}
It therefore suffices to prove that
\begin{align}\label{Aneededbd1}
\sum_{N_2 \leq N_1} \left|\int_{\R^9}  \RR^\beta {\langle D^\sharp \rangle}^{-\delta}f_{N_1,n_1}(x) \,g_{N_2,n_2, N_1}(x) \,\RR^\alpha h_{\leq N_2, n_3, N_1}(x)\, dx\right|\lesssim (\langle n_1\rangle\langle n_2\rangle \langle n_3 \rangle)^3\|f\|_{L^p}\|g\|_{L^q}\|h\|_{L^r}.
\end{align}
Here we have denoted, following the previous convention, 
\begin{gather*} f_{N_1,n_1} := {\FF^{\sharp}}^{-1}e^{\frac{2\pi i\, n_1\cdot \xi_1}{K N_1}} f_{N_1}^\sharp(\xi_1) \\
g_{N_2,n_2, N_1}(x):={\FF^\sharp}^{-1}e^{\frac{2\pi i\, n_2\cdot \xi_2}{K N_1}}g_{N_2}^{\sharp}(\xi_2) \\ h_{\leq N_2, n_3, N_1}(x):={\FF^\sharp}^{-1}e^{\frac{2\pi i\, n_3\cdot \xi_3}{K N_1}}h_{\leq N_2}^\sharp (\xi_3).\end{gather*}
We divide the above sum into three pieces: first $N_1,N_2 \geq 1$, second $N_1,N_2 \leq 1$ and finally $N_2 \leq 1 \leq N_1$.  In fact, the last piece can be treated by combining in an obvious fashion the arguments for the first two, so we shall skip it here.  For the first piece, we use (an elementary consequence of) Lemma \ref{square and maximal lemma} to obtain
\begin{equation*}
\begin{split}
&\sum_{1\leq N_2 \leq N_1} \left|\int_{\R^9}  \RR^\beta {\langle D^\sharp \rangle}^{-\delta}f_{N_1,n_1}(x) \,g_{N_2,n_2, N_1}(x) \,\RR^\alpha h_{\leq N_2, n_3, N_1}(x)\, dx\right| \\
& \qquad \qquad \lesssim \sum_{N_2\leq N_1} N_1^{-\delta} \| f_{N_1,n_1} \|_p \| g_{N_2,n_2, N_1} \|_q \| h_{\leq N_2, n_3, N_1} \|_r \lesssim (\langle n_1\rangle\langle n_2\rangle \langle n_3 \rangle)^3\|f\|_{L^p}\|g\|_{L^q}\|h\|_{L^r}.
\end{split}
\end{equation*}
Turning to the second piece, let $Q$ be defined by $\frac{1}{Q} = 1 - \frac{1}{\tilde p} - \frac{1}{\tilde r} < \frac{1}{\tilde q}$.  Then from Bernstein's inequality we have that
\begin{equation*}
\begin{split}
&\sum_{N_2 \leq N_1 \leq 1} \left|\int_{\R^9}  \RR^\beta {\langle D^\sharp \rangle}^{-\delta}f_{N_1,n_1}(x) \,g_{N_2,n_2, N_1}(x) \,\RR^\alpha h_{\leq N_2, n_3, N_1}(x)\, dx\right| \\
& \qquad \qquad \lesssim \sum_{N_2 \leq N_1 \leq 1} \| f_{N_1,n_1} \|_{\tilde p} \| g_{N_2,n_2, N_1} \|_Q \| h_{\leq N_2, n_3, N_1} \|_{\tilde r} \\
& \qquad \qquad \lesssim \sum_{N_2 \leq N_1 \leq 1} \| f_{N_1,n_1} \|_{\tilde p} N_2^{3(\frac{1}{\tilde q} - \frac{1}{Q})} \| g_{N_2,n_2, N_1} \|_{\tilde q} \| h_{\leq N_2, n_3, N_1} \|_{\tilde r} \\
& \qquad\qquad \lesssim (\langle n_1\rangle\langle n_2\rangle \langle n_3 \rangle)^3\|f\|_{L^p}\|g\|_{L^q}\|h\|_{L^r}.
\end{split}
\end{equation*}
This completes the proof.
\end{proof}

\section{Small-data scattering for a quadratic NLS} \label{ubar2} 

\subsection{The a priori estimate}

The proof of Theorem~\ref{intro:gwptheorem} essentially consists of an a priori estimate in the space $X$ given by the norm 
\beq \| f \|_X :=  \| x \Omega^* f(t,x) \|_{L^\infty_t L^2_x} + \|f \|_{L^\infty_t H^1_x}\sim \| |x| \Omega^* f(t,x) \|_{L^\infty_t L^2_x} + \|\Omega^* f \|_{L^\infty_t H^1_x},
\eeq
thanks to assumption {\bf H3*}. Due to the intertwining property of $\Omega$ and (distorted) Plancherel identity, the $X$ norm is equivalent to the following:
$$
\| f \|_X \sim \|f^\sharp\|_{L^\infty_t L^2_\xi} + \|\xi f^\sharp\|_{L^\infty_t L^2_\xi} + \|\partial_\xi f^\sharp\|_{L^\infty_t L^2_\xi}.
$$
Observe that Duhamel's formula for equation \eqref{QNLS}  can be written
\[ f = u_0 - i B(f, f),\]
where
\begin{equation}\label{Bff0} [B(f, f)]^\sharp(t,\xi) := \int_0^t \!\! \int \!\!\! \int e^{-is \phi(\xi,\eta,\zeta)} \overline{f^\sharp}(s,\eta) \overline{f^\sharp}(s,\zeta) \overline{\mathcal{M}}(\xi,\eta,\zeta)  \, d\eta \, d\zeta \, ds,\end{equation}
with $\phi(\xi,\eta,\zeta) = |\xi|^2 + |\eta|^2 + |\zeta|^2$. The key point is then to prove the following result.

\begin{proposition}[A priori bound]
\label{hummingbird}
If $f$ solves $f = u_0 - iB(f,f)$ and belongs to $X$ with a norm less than a fixed constant $\epsilon$, then
\[ \left\| B(f,f) \right\|_X \lesssim \|f\|_X^2 \]
\end{proposition}

The above proposition entails the inequality $\|f\|_X \lesssim \|u_0\|_X + \|f\|_X^2$, which gives an a priori control on $f$ in $X$. This leads to the main theorem by a continuous induction argument.  
In particular, thanks to the $L^2$-convergence of the integral in \eqref{Bff0}, which follows from \eqref{ubar2:basicdispersionestimate},  we also get that $f$ converges in $L^2$ as $t\to \infty$ to some $\phi_\infty$.  This is equivalent to the scattering of $u(t)$ to a linear solution $e^{itH}\phi_{\infty}$. By the boundedness of $\Omega$ on $L^2$, one also obtains scattering to a free solution $e^{it\Delta} \phi'_{\infty}$ for some $\phi'_\infty \in L^2$.

\subsection{Decomposition of $B(f,f)$}
We will first decompose $B$ into high and low frequencies:  call $\mathbf{I}$ the piece corresponding to $|(\xi,\eta,\zeta)|\gtrsim 1/\sqrt{s}$, and $\mathbf{II}$
the piece corresponding to $|(\xi,\eta,\zeta)|\lesssim 1/\sqrt{s}$. We then perform a normal form transformation on the high frequency piece (which is localized away from the time resonant set).
Finally, we estimate separately the resulting terms $\mathbf{I}_1$, $\mathbf{I}_2$, $\mathbf{I}_3$, $\mathbf{I}_4$.

To be more specific, let $\chi$ be a smooth cutoff function satisfying
 \[ \chi(\xi,\eta,\zeta) = \left\{ \begin{array}{ll} 1 &  \textrm{if } |\xi|^2 + |\eta|^2 + |\zeta|^2 < 1, \\
 0 &  \textrm{if } |\xi|^2 + |\eta|^2 + |\zeta|^2 \geq 4. \end{array} \right. \]
For any $s \geq 0$, we further define $\chi_s = \chi(\sqrt{s} \cdot)$; observe that the support of $\chi_s$ will be a subset of $B_{2/\sqrt{s}}(0)$.  
Using $\chi$ to partition the frequency domain, $B(f,f)$ becomes
$$
B(f,f)^\sharp =  \mathbf{I} + \mathbf{II} 
$$
with
\begin{align}
\displaystyle \mathbf{I} &:= \int_0^t  \!\!\int  \!\!\! \int (1-\chi_s)  e^{-is\phi } \overline{f^\sharp}(\eta) \overline{f^\sharp}(\zeta) \overline{\mathcal{M}}(\xi, \eta, \zeta) \, d\eta \, d\zeta \, ds \label{apriori:defI} \\
\displaystyle \mathbf{II} &:= \int_0^t \!\! \int  \!\!\! \int \chi_s  e^{-is\phi } \overline{f^\sharp}(\eta) \overline{f^\sharp}(\zeta) \overline{\mathcal{M}}(\xi,\eta,\zeta)  \, d\eta \, d\zeta \, ds. \label{apriori:defII}
\end{align}
We now perform a normal form transformation on $\mathbf{I}$ (integration by parts using the identity $-\frac{1}{i\phi} \partial_s e^{-is\phi} = e^{-is\phi}$). This gives
$$
\mathbf{I} = \mathbf{I}_1 +\mathbf{I}_2 + \mathbf{I}_3 + \mathbf{I}_4
$$
where
\begin{align}
\displaystyle \mathbf{I}_1 &:=  \int \!\!\! \int \frac{1-\chi_t}{i\phi} e^{-it\phi} \overline{f^\sharp}(\eta) \overline{f^\sharp}(\zeta) \overline{\mathcal{M}}  \, d\eta \, d\zeta \label{apriori:defI1} \\
\displaystyle \mathbf{I}_2 &:=  - \int_0^t \!\! \int  \!\!\! \int \frac{\partial_s \chi_s}{i\phi} e^{-is\phi} \overline{f^\sharp}(\eta) \overline{f^\sharp}(\zeta) \overline{\mathcal{M}}  \, d\eta \, d\zeta \, ds  \label{apriori:defI2} \\
\displaystyle \mathbf{I}_3 &:=  \int_0^t \!\! \int  \!\!\! \int \frac{1- \chi_s}{i\phi} e^{-is\phi} \partial_s \overline{f^\sharp}(\eta) \overline{f^\sharp}(\zeta) \overline{\mathcal{M}}  \, d\eta \, d\zeta \, ds \label{apriori:defI3}\\
\displaystyle \mathbf{I}_4 &:= \int_0^t \!\! \int  \!\!\! \int \frac{1- \chi_s}{i\phi} e^{-is\phi} \overline{f^\sharp}(\eta) \partial_s\overline{f^\sharp}(\zeta) \overline{\mathcal{M}}  \, d\eta \, d\zeta \, ds  \label{apriori:defI4}
\end{align}
Here we have suppressed some dependencies in the interest of readability; we will continue using this convention throughout this section.  

Finally, it will often be necessary to distinguish which of the frequencies is largest.  This will be the role of the cutoff functions $\Phi^{(j)} = \Phi^{(j)}(\xi,\eta,\zeta)$, $j=1,2,3$ that we now define.  These functions are chosen to be smooth, and homogeneous of degree $0$ outside of the ball of radius $1/10$. Furthermore, they add up to $1$ away from $(0,0,0)$:
\begin{subequations} \label{apriori:defPhi123}
\beq
\Phi^{(1)} + \Phi^{(2)} + \Phi^{(3)} = 1, \textrm{ for all } (\xi,\eta,\zeta) \textrm{ with } |\xi|^2+ |\eta|^2 +|\zeta|^2 \geq 1
\eeq
and are such that
\beq \begin{split}
\mbox{On } B_1(0)^c \cap \operatorname{Supp} \Phi^{(1)}, & \qquad |\xi| \gtrsim |\eta| + |\zeta| \\
\mbox{On } B_1(0)^c \cap \operatorname{Supp} \Phi^{(2)}, &\qquad |\eta| \gtrsim |\xi| + |\zeta| \\
\mbox{On } B_1(0)^c \cap \operatorname{Supp} \Phi^{(3)}, & \qquad|\zeta| \gtrsim |\xi|+  |\eta|.
\end{split} \eeq
\end{subequations}
In line with our convention, denote $\Phi_s^{(i)} = \Phi^{(i)}(\sqrt{s} \cdot)$, for each $s \geq 0$.

\subsection{A priori estimates: proof of Proposition~\ref{hummingbird}} \label{ubar2apriorisubsec}

In this section, we prove Proposition~\ref{hummingbird}. We begin with two simple observations. First, for any $g \in X$, 
\begin{align*}
\| g \|_{L^2(\mathbb{R}^3)} & \lesssim \| \Omega^* g \|_{L^2(\mathbb{R}^3)}  \lesssim \| x \Omega^* g \|_{L^2(B_1^c)} + \| \Omega^* g \|_{L^2(B_1)} \\
& \lesssim \|  x  \Omega^* g \|_{L^2(\mathbb{R}^3)} + \| \Omega^* g \|_{L^6(\mathbb{R}^3)} \lesssim \|  x  \Omega^* g \|_{L^2(\mathbb{R}^3)} + \| \nabla \Omega^* g \|_{L^2(\mathbb{R}^3)} \\
& \lesssim \|  x  \Omega^* g \|_{L^2(\mathbb{R}^3)} + \| D^\sharp g \|_{L^2(\mathbb{R}^3)}.
\end{align*}
Thus, to control the full $X$ norm, we need only prove $\| x \Omega^* B(f,f) \|_{L^2}, \| D^\sharp B(f,f) \|_{L^2} \lesssim \|f\|_X^2$.   Second, for $6/5 < p \leq 6$, we can bound the $L_t^\infty L_x^p$ norm of a function in terms of the $X$ norm by noticing that
\beq \| f \|_{L_t^\infty L_x^p} \lesssim \| \Omega^* f \|_{L_t^\infty L_x^p} \lesssim \left\{ \begin{array}{ll} \| \langle x \rangle \Omega^*  f \|_{L_t^\infty L_x^2} & \textrm{if } 6/5 < p \leq 2 \\
\| D^\sharp f \|_{L_t^\infty L_x^2}^\theta \| f \|_{L_t^\infty L_x^2}^{1-\theta} &  \textrm{if } 2 \leq p \leq 6 \end{array} \right\} \leq \| f \|_X.  \label{LptoweightedL2estimate} \eeq
Here, the estimate for $6/5 < p \leq 2$ is simply from H\"older's inequality, and for $2 \leq p \leq 6$ we interpolate between $L_x^2$ and $L_x^6$, controlling the latter with the (Euclidean) Sobolev inequality.

Even at this stage, it is clear that our proof will require considering an inordinate number of terms:   each of \eqref{apriori:defI1}--\eqref{apriori:defI4} may be carved into multiple pieces using the cutoff functions.  Indeed, we will compound the problem by, at certain points,  differentiating various terms as well.   Thankfully, the vast majority of these quantities can be estimated using slight variations on a few basic strategies.  With that in mind, we give a full treatment of each class of representative term, as well as the most problematic ones.  In hopes of keeping our exposition to a readable length, however, when a quantity can be controlled using a straightforward modification of an already stated argument, we will only provide a sketch.  

\begin{proof}[Proof of Proposition \ref{hummingbird}] By the preceding remarks, we know that we must control $B(f,f)$ in the weighted $L^2$ and $\dot{H}_\sharp^1$ norms in terms of $\|f \|_X^2$.  Let us estimate each of the pieces $\mathbf{I}_1, \ldots, \mathbf{I}_4, \mathbf{II}$ in these spaces one at a time. 

\subsubsection{Weighted $L^2$ estimate of $\mathbf{I}_1$ for $t > 1$.} 
We begin by  controlling the long time weighted norm $\|  x \Omega^* {\mathcal{F}^\sharp}^{-1} \eqref{apriori:defI1} \|_{L^2}\sim \|\partial_\xi \eqref{apriori:defI1}\|_{L^2}$.  Differentiating $\mathbf{I}_1$ in $\xi$, and using~\eqref{derivative of M}, yields (sticking to the convention that repeated indices are summed)
\begin{subequations}
\begin{align}
\label{apple1}
\partial_\xi^\alpha \mathbf{I}_1 = & \int  \!\!\! \int \frac{t}{\sqrt{\phi}} m_t^1 e^{-it\phi} \overline{f^\sharp}(\eta) \overline{f^\sharp}(\zeta) \overline{\mathcal{M}}  \, d\eta \, d\zeta \\
\label{apple2}
& - \int  \!\!\! \int \frac{t}{\sqrt{\phi}} m_t^2 e^{-it\phi} \overline{f^\sharp}(\eta) \overline{f^\sharp}(\zeta) R_\eta^\beta R_\xi^\alpha \overline{\mathcal{M}}  \, d\eta \, d\zeta \\
\label{apple3}
& - \int  \!\!\! \int \frac{1-\chi_t}{i\phi} e^{-it\phi} \partial_{\eta^\beta} \overline{f^\sharp}(\eta) \overline{f^\sharp}(\zeta) R_\eta^\beta R_\xi^\alpha \overline{\mathcal{M}}   \, d\eta \, d\zeta \\
\label{apple4}
& + \int  \!\!\! \int \frac{1-\chi_t}{i\phi} e^{-it\phi} \overline{f^\sharp}(\eta) \overline{f^\sharp}(\zeta) \left[ R_{\xi}^\alpha E_\eta + B_\xi R^\alpha_\eta \right] \overline{\mathcal{M}} \, d\eta \, d\zeta.
\end{align}
\end{subequations}
where
\begin{equation*}
\begin{split}
& m_t^1 = \frac{(\partial_{\xi^\alpha} \chi)_t}{\sqrt{t \phi}} + (1-\chi_t)\frac{\partial_{\xi^\alpha} \phi}{t\phi \sqrt{\phi}} - (1-\chi_t)\frac{\partial_{\xi^\alpha} \phi}{ \sqrt{\phi}} \\
& m_t^2 = \frac{(\partial_{\eta^\beta} \chi)_t}{\sqrt{t \phi}} + (1-\chi_t)\frac{\partial_{\eta^\beta} \phi}{t\phi \sqrt{\phi}} - (1-\chi_t)\frac{\partial_{\eta^\beta} \phi}{\sqrt{\phi}}
\end{split}
\end{equation*}
satisfy uniform Coifman--Meyer bounds and are supported outside of a ball of  approximate radius $1/\sqrt{t}$.  Notice that \eqref{apple2}--\eqref{apple4} arise when we use \eqref{derivative of M} to evaluate $\partial_\xi \overline{\mathcal{M}}$, and then perform an integration by parts in $\eta$. To simplify the notation, we omit from now on the superscripts $\alpha$ and $\beta$ when there is no chance for ambiguity.

We start our analysis with \eqref{apple1}. For any $\delta > 0$ sufficiently small, it can be bounded by
\begin{equation*}
\begin{split}
\|\eqref{apple1}\|_{L^2_\xi} & = \| t \int\!\!\! \int \frac{|\eta|}{\sqrt{\phi}} m^1_t e^{-it\phi} \frac{1}{|\eta|} \overline{f^\sharp}(\eta) \overline{f^\sharp}(\zeta) \overline{\mathcal{M}} \,d\eta\,d\zeta \|_{L^2_\xi} \\
& \lesssim t \left[ \|e^{itH} \Lambda^{-1} f \|_{L^3} \| e^{itH} f \|_{L^6} + \|e^{itH} \Lambda^{-1} f \|_{L^{3-\delta}} \| e^{itH} f \|_{L^6} \right] \\
& \lesssim t \| f \|_{L^{5/4}} \| e^{itH} f \|_{L^6} \\
& \lesssim \|f\|_X^2.
\end{split}
\end{equation*}
Here we have used Proposition~\ref{Bernsteinprop},~\eqref{LptoweightedL2estimate}, and expressed 
\[ \frac{|\eta|}{\sqrt{\phi}}=\sum_{\gamma=1}^3 \frac{\eta^\gamma}{\sqrt \phi} \frac{\eta^\gamma}{|\eta|}\]
 in order to apply Theorem \ref{intro:CMtheorem} with the Coifman-Meyer symbols ${\eta^\gamma}/{\sqrt \phi}$.  This last step relies on the $L^p$ boundedness of the operator ${D^\sharp}^\gamma/|D^\sharp|$, which follows from \eqref{intertwining} and the boundedness of the Euclidean Riesz transform.


The term~(\ref{apple2}) can be bounded similarly, so we turn directly to~(\ref{apple3}). We need to decompose this term further, with the help of the cutoff functions $\chi$, $1-\chi$, $\Phi_t^{(1)}$, $\Phi_t^{(2)}$, $\Phi_t^{(3)}$. This gives six elementary pieces, of which we will only present the details for the ones involving $\Phi^{(2)}_t$; those with $\Phi_t^{(1)}$ and $\Phi_t^{(3)}$ can be managed similarly, with the small modification of a duality argument for $\Phi_t^{(1)}$. Thus we consider 
\begin{subequations}
\begin{align}
\label{robin1}
& \int  \!\!\! \int \frac{1-\chi_t}{i\phi} \chi \Phi^{(2)}_t e^{-it\phi} \partial_\eta \overline{f^\sharp}(\eta) \overline{f^\sharp}(\zeta) R_\eta^\beta R_\xi^\alpha \overline{\mathcal{M}}   \, d\eta \, d\zeta \\
\label{robin2}
& \qquad \qquad +\int \!\!\! \int \frac{1-\chi_t}{i\phi} (1-\chi) \Phi^{(2)}_t e^{-it\phi} \partial_\eta \overline{f^\sharp}(\eta) \overline{f^\sharp}(\zeta) R_\eta^\beta R_\xi^\alpha \overline{\mathcal{M}}   \, d\eta \, d\zeta
\end{align}
\end{subequations}
To estimate \eqref{robin1}, we use Theorem \ref{colibri} and Proposition~\ref{Bernsteinprop}: for any $\delta$ sufficiently small, we have
\begin{equation*}
\begin{split}
\| \eqref{robin1} \|_{L^2} & = \|  \int  \!\!\! \int \frac{(1-\chi_t)|\eta|^2}{i\phi} \chi \Phi^{(2)}_t e^{-it\phi} \frac{1}{|\eta|^2} \partial_\eta \overline{f^\sharp}(\eta) \overline{f^\sharp}(\zeta) R_\eta^\beta R_\xi^\alpha \overline{\mathcal{M}}   \, d\eta \, d\zeta \|_{L^2} \\
& \lesssim \|e^{-itH} \Lambda_t^{-2} \Omega x \Omega^* f \|_{L^3} \|e^{-itH} f \|_{L^6} + \|e^{-itH} \Lambda_t^{-2} \Omega x \Omega^* f \|_{L^3} \|e^{-itH} f \|_{L^{6-\delta}} \\
& \lesssim t^{3/4} \|\Omega x \Omega^* f\|_{L^2} \frac{1}{t^{7/8}} \|f\|_X \lesssim \frac{1}{t^{1/8}} \|f\|_X^2.
\end{split}
\end{equation*}
To estimate~(\ref{robin2}), also rely on Theorem~\ref{colibri} to get
\begin{equation*}
\begin{split}
\| \eqref{robin2} \|_{L^2} & = \|  \int  \!\!\! \int \frac{(1-\chi_t)|\eta|}{i\phi} (1-\chi) \Phi^{(2)}_t e^{-it\phi} \frac{1}{|\eta|} \partial_\eta \overline{f^\sharp}(\eta) \overline{f^\sharp}(\zeta) R_\eta^\beta R_\xi^\alpha \overline{\mathcal{M}}   \, d\eta \, d\zeta \|_{L^2} \\
& \lesssim \|e^{-itH} \Lambda_1^{-1} \Omega x \Omega^* f \|_{L^3} \|e^{-itH} f \|_{L^6} + \|e^{-itH} \Lambda_1^{-1} \Omega x \Omega^* f \|_{L^3} \|e^{-itH} f \|_{L^{6-\delta}} \\
& \lesssim \|\Omega x \Omega^* f\|_{L^2} \frac{1}{t^{7/8}} \|f\|_X \lesssim \frac{1}{t^{7/8}} \|f\|_X^2.
\end{split}
\end{equation*}
The remaining  term~\eqref{apple4} is treated the same way.

\subsubsection{Weighted $L^2$ estimate of $\mathbf{I}_1$ for $0 \leq t \leq  1$.}  Here the general idea is to repeat the arguments of the previous subsection, but systematically replace the dispersion estimates with the (distorted) Sobolev inequality of Proposition \ref{fractionalintprop} (b). For instance,
\begin{equation*}
\begin{split}
\|\eqref{apple1}\|_{L^2} & =  \| t \int \!\!\! \int \frac{|\eta|}{\sqrt{\phi}} m_t^{1} e^{-it\phi} \frac{1}{|\eta|}\overline{f^\sharp}(\eta) \overline{f^\sharp}(\zeta) \overline{\mathcal{M}} \,d\eta\,d\zeta \|_{L^2} \\
& \lesssim t \left[ \|e^{itH} \Lambda^{-1} f \|_{L^3} \| e^{itH} f \|_{L^6} + \|e^{itH} \Lambda^{-1} f \|_{L^{3-\delta}} \| e^{itH} f \|_{L^6} \right] \\
& \lesssim t \| f \|_{L^{5/4}} \| f \|_{H^1_\sharp} \lesssim t \|f\|_X^2.
\end{split}
\end{equation*}
Similarly, to estimate \eqref{apple2}, we employ the cutoff functions.  To bound a representative term of the form \eqref{robin1}, we note that
\begin{equation*}
\begin{split}
\| \eqref{robin1} \|_{L^2} & = \|  \int  \!\!\! \int \frac{(1-\chi_t)|\eta|^2}{i\phi} \chi \Phi^{(2)}_t e^{-it\phi} \frac{1}{|\eta|^2} \partial_\eta \overline{f^\sharp}(\eta) \overline{f^\sharp}(\zeta) R_\eta^\beta R_\xi^\alpha \overline{\mathcal{M}}   \, d\eta \, d\zeta \| \\
& \lesssim \|e^{-itH} \Lambda_t^{-2} \Omega x \Omega^* f \|_{L^3} \|e^{-itH} f \|_{L^6} + \|e^{-itH} \Lambda_t^{-2} \Omega x \Omega^* f \|_{L^3} \|e^{-itH} f \|_{L^{6-\delta}} \\
& \lesssim t^{3/4} \|\Omega x \Omega^* f\|_{L^2} \|f\|_{H^1_\sharp} \lesssim t^{3/4} \|f\|_X^2.
\end{split}
\end{equation*}
All in all, this gives the desired bound $\|\partial_\xi \mathbf{I}_1\|_2 \lesssim \| f \|_X^2$ for $t<1$.

\subsubsection{Weighted $L^2$ estimate of $\mathbf{I}_2$ for $t > 1$.} \label{apriori:I2wL2bigt} 
First, we need to split the time integral giving $\mathbf{I}_2$ into two pieces:
$$
\mathbf{I}_2 = \int_0^1 \dots + \int_1^t \dots =: \mathbf{I}_2^1 + \mathbf{I}_2^2.
$$
The term $\mathbf{I}_2^1$ will be dealt with when we consider $\mathbf{I}_2$ for $t < 1$; thus for now, we only 
consider the term $\mathbf{I}_2^1$. Differentiating it in $\xi$ gives
\begin{subequations}
\begin{align}
\label{tangerine1}
\partial_{\xi^\alpha} \mathbf{I}_2^2 & = \int_1^t \int  \!\!\! \int \sqrt{s} m_s^1 e^{-is\phi} \overline{f^\sharp}(\eta) \overline{f^\sharp}(\zeta) \overline{\mathcal{M}}  \, d\eta \, d\zeta \\
\label{tangerine2}
& \qquad+\int_1^t \int  \!\!\! \int \sqrt{s} m_s^2 e^{-is\phi} \overline{f^\sharp}(\eta) \overline{f^\sharp}(\zeta) R_\eta^\beta R_\xi^\alpha \overline{\mathcal{M}}  \, d\eta \, d\zeta \\
\label{tangerine3}
& \qquad+\int_1^t  \int  \!\!\! \int m_s^3 e^{-is\phi} \partial_{\eta^\beta} \overline{f^\sharp}(\eta) \overline{f^\sharp}(\zeta) R_\eta^\beta R_\xi^\alpha \overline{\mathcal{M}}   \, d\eta \, d\zeta \\
\label{tangerine4}
& \qquad+\int_1^t \int  \!\!\! \int m_s^3 e^{-is\phi} \overline{f^\sharp}(\eta) \overline{f^\sharp}(\zeta) \left[ R_{\xi}^\alpha E_\eta + B_\xi R^\alpha_\eta \right] \overline{\mathcal{M}} \, d\eta \, d\zeta,
\end{align}
\end{subequations}
where $\nabla$ refers to the gradient with respect to $(\xi,\eta,\zeta)$.  Here
\begin{equation*}
\begin{split}
& m_s^1 := -\partial_{\xi^\alpha} \left( \frac{(\xi,\eta,\zeta) \cdot (\nabla \chi)_s}{s \phi} \right) +i \frac{(\xi,\eta,\zeta) \cdot (\nabla \chi)_s}{\sqrt{s} \phi} \sqrt{s} \partial_{\xi^\alpha} \phi, \\
& m_s^2 := -\partial_{\eta^\beta} \left( \frac{(\xi,\eta,\zeta) \cdot (\nabla \chi)_s}{s \phi} \right) +i \frac{(\xi,\eta,\zeta) \cdot (\nabla \chi)_s}{\sqrt{s} \phi} \sqrt{s} \partial_{\eta^\beta} \phi, \\
& m_s^3 := -\frac{(\xi,\eta,\zeta) \cdot (\nabla \chi)_s}{2\sqrt{s} i\phi}
\end{split}
\end{equation*}
are supported on annuli with inner and outer radii approximately $1/\sqrt{s}$, and satisfy uniform Coifman--Meyer bounds. 
We begin with~(\ref{tangerine1}) which can be controlled with the help of Theorem~\ref{intro:CMtheorem} and Bernstein inequalities: for $6/5<p<4/3$ and $\delta > 0$ sufficiently small
\beq \begin{split}
 \| \eqref{tangerine1}\|_{L^2} & \lesssim_\delta \int_1^t \sqrt{s} \|e^{isH} P_{<1/\sqrt{s}} f \|_{L^3} \| e^{isH} f\|_{L^6}  \, ds \\
 & \qquad + \int_1^t \sqrt{s} \|e^{isH} P_{<1/\sqrt{s}} f \|_{L^3} \| e^{isH} f\|_{L^{6-\delta}}  \, ds \\
 & \lesssim \int_1^t \sqrt{s} s^{\frac{1}{2} - \frac{3}{2p}}  \| f \|_{L^{p}} \frac{\|f\|_X}{s^{7/8}}  \, ds\\
 & \lesssim \|f\|_X^2 \int_1^t \frac{1}{s^{3/(2p) - 1/8}} \, ds \lesssim \|f\|_X^2.\end{split} \eeq
The term~(\ref{tangerine2}) admits a similar treatment, thus we skip it and focus on~(\ref{tangerine3}). Due to the presence of the (Euclidean) Riesz transforms, we will need to use Theorem~\ref{colibri}.  This is permissible because, for $t > 1$, $m_s^3$ enjoys the improved   Coifman--Meyer structure of \eqref{bettermult}.  Thus, for $t > 1$, we may select $\delta > 0$ small enough so that
\begin{equation*}
\begin{split}
\| \eqref{tangerine3} \|_{L^2} &\lesssim \int_1^t \| \int  \!\!\! \int  m_s^3 e^{-is\phi} \partial_\eta \overline{f^\sharp}(\eta) \overline{f^\sharp}(\zeta) R_\eta^\beta R_\xi^\alpha \overline{\mathcal{M}}   \, d\eta \, d\zeta  \|_{L^2} \,ds \\
&\lesssim_\delta \int_1^t \left[ \| e^{isH}P_{< 1/\sqrt{s}} \Omega x \Omega^* f  \|_{L^3} \| e^{isH} f \|_{L^6} +  \| e^{isH} P_{< 1/\sqrt{s}}\Omega x \Omega^* f  \|_{L^3} \| e^{isH} f \|_{L^{6-\delta}} \right]\,ds \\
&\lesssim \int_1^t \frac{1}{s^{1/4}} \|\Omega x \Omega^* f \|_{L^2} \frac{\|f\|_X}{s^{7/8}} \,ds \\
& \lesssim \|f\|_X^2.
\end{split}
\end{equation*}
Finally, \eqref{tangerine4} can be estimated similarly, thus we skip it.

\subsubsection{Weighted $L^2$ estimate of $\mathbf{I}_2$ for $0 \leq t \leq 1$.} \label{apriori:I2wL2smallt} Once again, it suffices to follow the argument of the previous subsection merely replacing the dispersive estimates with appeals to the (distorted) Sobolev embedding inequalities.

\subsubsection{Weighted $L^2$ estimate of $\mathbf{I}_3$ for $t > 1$.}  First let us point out that, due to symmetry, the argument we give for $\mathbf{I}_3$ applies to $\mathbf{I}_4$ as well.  
The main idea in both cases is that we can use the equation satisfied by $u$ to replace the quadratic nonlinearity with a cubic one.  
Indeed, this is immediately apparent since, by the definition of the profile $f$ and the equation satisfied by $u$,
\beq \partial_s \bar{f} = \partial_s (e^{-isH} \bar{u}) = e^{-isH} \bar{u}^2.\label{ubar2apriori:Icidentity} \eeq
Using \eqref{ubar2apriori:Icidentity} will give us additional decay, but it presents a technical nuisance in the form of an asymmetry in the integrand. 
We start by splitting the integral defining $\mathbf{I}_3$ into
$$
\mathbf{I}_3 = \int_0^1 \dots + \int_1^t \dots =: \mathbf{I}_3^1 + \mathbf{I}_3^2.
$$
The estimate of $\mathbf{I}_3^1$ follows from the argument we give for $\mathbf{I}_3$ with $0 < t \leq 1$ in the following subsection; so let us focus on $\mathbf{I}_3^2$.  Using \eqref{derivative of M}, we see that differentiating it in $\xi$ gives
\begin{subequations}
\begin{align}
\partial_\xi \mathbf{I}_3^2 & =  -\int_1^t \!\! \int  \!\!\! \int \frac{s}{\sqrt{\phi}} m_s^1 e^{-i s \phi} \partial_s \overline{f^\sharp}(\eta) \overline{f^\sharp}(\zeta) \overline{\mathcal{M}}  \, d\eta \, d\zeta \, ds \label{pear1} \\
& \qquad  + \int_1^t \!\! \int  \!\!\! \int \frac{s}{\sqrt{\phi}} m_s^2 e^{-i s \phi} \partial_s \overline{f^\sharp}(\eta) \overline{f^\sharp}(\zeta) R_\zeta^\beta R_\xi^\alpha \overline{\mathcal{M}}  \, d\eta \, d\zeta \, ds \label{pear2} \\
& \qquad  + \int_1^t \!\! \int  \!\!\! \int \frac{1- \chi_s}{i\phi} e^{-is \phi} \partial_s \overline{f^\sharp}(\eta) \partial_\zeta \overline{f^\sharp}(\zeta) R_\zeta^\beta R_\xi^\alpha \overline{\mathcal{M}} \, d\eta \, d\zeta \, ds \label{pear3} \\
& \qquad + \int_1^t  \!\! \int  \!\!\! \int \frac{1- \chi_s}{i\phi} e^{-is \phi} \partial_s \overline{f^\sharp}(\eta) \overline{f^\sharp}(\zeta) [R^\alpha_\xi E_\zeta + B_\xi R^\alpha_\zeta ] \overline{\mathcal{M}} \, d\eta \, d\zeta \, ds \label{pear4} 
\end{align}
\end{subequations}
where in this case
\begin{align*}
& m_s^1 := -\frac{(\partial_\xi \chi)_s}{\sqrt{s} \sqrt{\phi}} + (1-\chi_s) \frac{\partial_\xi \phi}{s \phi \sqrt{\phi}} + (1 - \chi_s) \frac{\partial_\xi \phi}{\sqrt{\phi}} \\
& m_s^2 := -\frac{(\partial_\zeta \chi)_s}{\sqrt{s} \sqrt{\phi}} + (1-\chi_s) \frac{\partial_\zeta \phi}{s \phi \sqrt{\phi}} + (1 - \chi_s) \frac{\partial_\zeta \phi}{\sqrt{\phi}} 
\end{align*}
satisfy uniform Coifman-Meyer bounds and are supported outside of a ball of radius $\sim 1/\sqrt{s}$.

First consider \eqref{pear1}.  Introducing the functions $\Phi^{(j)}$, and arguing as before, we wish to establish the $L_\xi^2$ boundedness of terms of the form
\beq \int_1^t \!\! \int  \!\!\! \int \frac{s m_s^1}{\sqrt{\phi}} \Phi_s^{(j)} e^{-is\phi} \partial_s\overline{f^\sharp}(\eta) \overline{f^\sharp}(\zeta) \overline{\mathcal{M}}  \, d\eta \, d\zeta \, ds, \qquad j=1,2,3.  \label{ubar2apriori:Icwproto1} \eeq
We simply deal with the term involving $\Phi_s^{(2)}$, since it is possible to treat the other ones in a similar way.
\begin{align*}
\| \eqref{ubar2apriori:Icwproto1}\|_{L_\xi^2} & \lesssim \int_1^t s \| \int  \!\!\! \int \frac{ m_s^1}{\sqrt{\phi}} \Phi_s^{(2)} e^{-is\phi} \partial_s\overline{f^\sharp}(\eta) \overline{f^\sharp}(\zeta) \overline{\mathcal{M}}  \, d\eta \, d\zeta \|_{L^2_\xi} \, ds \\
& \lesssim_\delta \int_1^t s  \| \Lambda_s^{-1} e^{-isH} \partial_s \bar{f}\|_{L^3} (\| e^{isH} f\|_{L^6} + \| e^{isH} f\|_{L^{6-\delta}} ) \, ds  \\
& \lesssim \int_1^t s ( s^{1/4} \| e^{-isH} \partial_s \bar{f}\|_{L^2})  (\frac{1}{s^{7/8}} \|f\|_X) \, ds \\
& \lesssim \int_1^t s ( s^{1/4} \frac{1}{s^{3/2}} \|f\|_X^2)  (\frac{1}{s^{7/8}} \|f\|_X) \, ds \lesssim \|f\|_X^3.\end{align*}
Here we have used the fact that 
\[ \| e^{-isH} \partial_s \bar{f}\|_{L^2} = \| \overline{u}^2 \|_{L^2} = \| e^{isH} f \|_{L^4}^2  \lesssim  s^{-3/2} \|f\|_X^2.\]

The term~(\ref{pear2}) is almost identical, thus we immediately consider \eqref{pear3}. It can be split as
\begin{subequations}
\begin{align}
\eqref{pear3} = & \sum_{j=1}^3  \int_1^t  \| \int  \!\!\! \int \frac{1- \chi_s}{i\phi} (1-\chi) \Phi^{(j)}_s  e^{-is\phi}  \partial_s \overline{f^\sharp}(\eta) \partial_\zeta\overline{f^\sharp}(\zeta) R_\zeta^\beta R_\xi^\alpha \overline{\mathcal{M}}   \, d\eta \, d\zeta \|_{L^2_\xi} \, ds \label{grapefruit2a} \\
& \qquad +  \sum_{j=1}^3 \int_1^t  \| \int  \!\!\! \int \frac{1- \chi_s}{i\phi} \chi  \Phi^{(j)}_s e^{-is\phi}  \partial_s \overline{f^\sharp}(\eta) \partial_\zeta\overline{f^\sharp}(\zeta) R_\zeta^\beta R_\xi^\alpha \overline{\mathcal{M}}   \, d\eta \, d\zeta \|_{L^2_\xi} \, ds,  \label{grapefruit2b} 
\end{align} 
\end{subequations}
We single out one of these terms,~(\ref{grapefruit2b}) for $j=2$, the other being similar. It can be estimated as follows:
\begin{equation*}
\begin{split}
\left\| (\ref{grapefruit2b})_{j=2} \right\|_2 & \leq \int_1^t 
\| \int  \!\!\! \int \frac{(1- \chi_s)|\eta||\zeta|}{i\phi} \chi \Phi^{(2)}_s  e^{-is\phi} \frac{1}{|\eta|} \partial_s \overline{f^\sharp}(\eta)\frac{1}{|\zeta|} \partial_\zeta\overline{f^\sharp}(\zeta) R_\zeta^\beta R_\xi^\alpha \overline{\mathcal{M}}   \, d\eta \, d\zeta \|_{L^2_\xi} \, ds \\
& \lesssim \int_1^t \left[ \left\| \Lambda_s^{-1} e^{-isH} \partial_s f \right\|_3 + \left\| \Lambda_s^{-1} e^{-itH} \partial_s f \right\|_{3-\delta} \right] \left\| \Lambda^{-1} e^{-isH} \Omega x \Omega^* f \right\|_6 \,ds\\
& \lesssim \int_1^t \sqrt{s} \left[ \left\| u^2 \right\|_3 + \left\| u^2 \right\|_{3-\delta} \right] \|\Omega x \Omega^* f \|_2\,ds \\
& \lesssim \int_1^t \sqrt{s} \frac{1}{s^{7/4}} \|u\|_X^2 \|u\|_X \,ds \lesssim \|u\|_X^3.
\end{split}
\end{equation*}

\subsubsection{Weighted $L^2$ estimate of $\mathbf{I}_3$ for $0 \leq t \leq 1$.}  Again, in light of the distorted Sobolev embedding theorem, the small time bounds are essentially a corollary of the long time bounds.

\subsubsection{$\dot{H}_\sharp^1$ estimates of $\mathbf{I}$.} 
 Let us consider one at a time the terms $\mathbf{I}_1$ to $\mathbf{I}_4$.
Multiplying $\mathbf{I}_1$ by $\xi$, partitioning in frequency, and then taking the $L^2_\xi$ norms gives
\begin{align}
\| {\mathcal{F}^\sharp}^{-1} \mathbf{I}_1\|_{\dot{H}_\sharp^1} & \leq \| \int  \!\!\! \int \xi \frac{1-\chi_t}{i\phi} e^{-it\phi} \overline{f^\sharp}(\eta) \overline{f^\sharp}(\zeta) \Phi^{(1)}_t \overline{\mathcal{M}}   \, d\eta \, d\zeta \|_{L_\xi^2} \label{ubar2apriori:IaH1i}\\ 
& \qquad + \| \int  \!\!\! \int \xi \frac{1-\chi_t}{i\phi} e^{-it\phi} \overline{f^\sharp}(\eta) \overline{f^\sharp}(\zeta) \Phi^{(2)}_t \overline{\mathcal{M}}  \, d\eta \, d\zeta \|_{L_\xi^2} \label{ubar2apriori:IaH1ii} \\
& \qquad  + \| \int  \!\!\! \int \xi \frac{1-\chi_t}{i\phi} e^{-it\phi} \overline{f^\sharp}(\eta) \overline{f^\sharp}(\zeta) \Phi^{(3)}_t \overline{\mathcal{M}}  \, d\eta \, d\zeta \|_{L_\xi^2}. \label{ubar2apriori:IaH1iii} \end{align}
The procedure for each of these is similar, so let us consider only the second one.  For small time we may estimate by Theorem~\ref{intro:CMtheorem} and the Sobolev embedding theorem: if $t<1$,
\begin{align*}
\eqref{ubar2apriori:IaH1ii} & \lesssim \|\int  \!\!\! \int \xi |\eta| \frac{1-\chi_t}{i\phi} e^{-it\phi} \frac{1}{|\eta|} \overline{f^\sharp}(\eta) \overline{f^\sharp}(\zeta) \Phi^{(2)}_t \overline{\mathcal{M}}  \, d\eta \, d\zeta \|_{L_\xi^2} \\
& \lesssim \| \Lambda_t^{-1} e^{itH} f \|_{L^6} (\| e^{itH} f \|_{L^3}+ \| e^{itH} f \|_{L^{3-\delta}} )\\
& \lesssim \sqrt{t} \|f\|_{H^1_\sharp}^2 \lesssim \|f\|_X^2 .\end{align*}
The long time control is found in the same way only using dispersion: if $t>1$,
\begin{align*}
\eqref{ubar2apriori:IaH1ii} & \lesssim \| \Lambda_t^{-1} e^{itH} f \|_{L^6} (\| e^{itH} f \|_{L^3}+ \| e^{itH} f \|_{L^{3-\delta}} ) \\
& \lesssim \sqrt{t} \| e^{itH} f \|_{L^6} (\| e^{itH} f \|_{L^3}+ \| e^{itH} f \|_{L^{3-\delta}} ) \\
& \lesssim \sqrt{t} \frac{1}{t} \|f\|_X^2.
\end{align*}

Next, multiplying $\mathbf{I}_2$ by $\xi$, and then taking the $L^2_\xi$ norm, we arrive at   
\begin{align*}
 \| {\mathcal{F}^\sharp}^{-1} \mathbf{I}_2 \|_{\dot{H}^1_{\sharp}} & \lesssim\int_0^t  \frac{1}{\sqrt{s}} \| \!\! \int  \!\!\! \int \xi \frac{(\xi, \eta, \zeta) \cdot (\nabla \chi)_s}{i\phi} e^{-is\phi} \overline{f^\sharp}(\eta) \overline{f^\sharp}(\zeta) \overline{\mathcal{M}}  \, d\eta \, d\zeta    \|_{L^2_\xi} \, ds \\
 & =: \int_0^t  \frac{1}{\sqrt{s}} \| \!\! \int  \!\!\! \int m_s e^{-is\phi} \overline{f^\sharp}(\eta) \overline{f^\sharp}(\zeta) \overline{\mathcal{M}}  \, d\eta \, d\zeta    \|_{L^2_\xi} \, ds,  \end{align*}
 where $m_s$ is a Coifman--Meyer symbol supported on an annulus of size $\sim 1/\sqrt{s}$. For the small time control we use the (distorted) Sobolev embedding: if $t<1$
\begin{align*}
 \| {\mathcal{F}^\sharp}^{-1} \mathbf{I}_2 \|_{\dot{H}^1_{\sharp}} & \lesssim \int_0^t \frac{1}{\sqrt{s}} \|  e^{isH} f \|_{L^3} \| e^{isH} f\|_{L^6} \, ds  + \int_0^t \frac{1}{\sqrt{s}} \| e^{isH}  f \|_{L^{3-\delta}} \|  e^{isH} f\|_{L^{6}} \, ds \\
 & \lesssim \int_0^t \frac{1}{\sqrt{s}}  \|f\|_{H^1_\sharp}^2 \, ds \lesssim \|f\|_X^2. \end{align*}
The long time bounds are the proved in the same way with the dispersive estimates in place of the Sobolev inequality.  

Finally, terms $\mathbf{I}_3$ and $\mathbf{I}_4$ are equivalent by symmetry, so let us only discuss the former.   We being by splitting the time integral as 
$$
\mathbf{I}_3 = \int_0^1 \dots + \int_1^t \dots =: \mathbf{I}_{3}^1 + \mathbf{I}_{3}^2.
$$
For $\mathbf{I}_3^2$, we partition the frequency space using the cutoff functions $\Phi_s^{(j)}$, $j = 1,2,3$.  The resulting pieces are treated similarly, so we consider only the representative term with  $j=2$.  In that case, we argue as follows:  if $t > 1$, 
\begin{equation*}
\begin{split}
\| \mathcal{F}^{\sharp -1} \mathbf{I}_{3}^2 \|_{\dot{H}^1_\sharp} 
& \lesssim \int_1^t \| \int\!\!\! \int \frac{\xi}{\phi} (1 - \chi_s) \Phi_s^{(2)} e^{-is\phi} \partial_s f^\sharp (\eta) f^\sharp(\zeta) \overline{\mathcal{M}} \,d\eta\,d\zeta \|_{L^2_\xi} \,ds \\
& \lesssim \int_1^t \sqrt{s} \left[ \|e^{isH} \partial_s f\|_{L^3} \|e^{isH} f \|_{L^6} + \|e^{isH} \partial_s f\|_{L^3} \|e^{isH} f \|_{L^{6-\delta}} \right] \,ds \\
& \lesssim \int_1^t \left[  \sqrt{s} \|u^2\|_{L^3} \|u\|_{L^6} + \sqrt{s} \|u^2\|_{L^3} \|u\|_{L^{6-\delta}} \right]\,ds \\
& \lesssim \int_1^t \sqrt{s} \|f\|_X^3 \frac{1}{s} \frac{1}{s} \frac{1}{\sqrt{s}} \,ds \lesssim \|f\|_X^3.
\end{split}
\end{equation*}
For, $\mathbf{I}_3^1$, we repeat the argument above substituting as usual the Sobolev embedding theorem for the dispersive estimates.  Thus 
\begin{equation*}
\begin{split}
\|\mathcal{F}^{\sharp -1} \mathbf{I}_3^1 \|_{\dot{H}^1_\sharp} & \lesssim \int_0^1 \left[  \sqrt{s} \|u^2\|_{L^3} \|u\|_{L^6} + \sqrt{s} \|u^2\|_{L^3} \|u\|_{L^{6-\delta}} \right]\,ds \\
& \lesssim \int_0^1 \sqrt{s} \|f\|_{H^1_\sharp}^3\,ds \lesssim \|f\|_X^3.
\end{split}
\end{equation*}
Notice that the same argument works to bound $\| {\mathcal{F}^\sharp}^{-1} \mathbf{I}_3\|_{\dot{H}_\sharp^1}$ for $t < 1$.

\subsubsection{Weighted $L^2$ estimate of $\mathbf{II}$ for $t>1$.}  In contrast to the normal forms method used to control $\mathbf{I}$, when estimating $\mathbf{II}$ 
we only want to use the fact that the integrand is supported compactly due to the presence of the $\chi_s$ factor.  Recall that
\[ \mathbf{II} = \int_0^t \!\! \int  \!\!\! \int \chi_s  e^{-is\phi } \overline{f^\sharp}(\eta) \overline{f^\sharp}(\zeta) \overline{\mathcal{M}}  \, d\eta \, d\zeta \, ds,\]
We split it as
$$
\mathbf{II} = \int_0^1 \dots + \int_1^t \dots =: \mathbf{II}^1 + \mathbf{II}^2.
$$
We start with $\mathbf{II}^2$; the term $\mathbf{II}^1$ will be estimated along with the case $t<1$ in the following subsection.  Applying $\partial_\xi$ and using \eqref{derivative of M} to evaluate $\partial_\xi \mathcal{M}$ gives
\begin{subequations} 
\begin{align} \partial_\xi \mathbf{II}^2 & = \int_1^t \!\! \int  \!\!\!\int \sqrt s  m_s^1 e^{-is\phi } \overline{f^\sharp}(\eta) \overline{f^\sharp}(\zeta) \overline{\mathcal{M}} \, d\eta \, d\zeta \, ds \label{cassis1} \\
& \qquad - \int_1^t \!\! \int  \!\!\!\int \sqrt s m_s^2 e^{-is\phi } \overline{f^\sharp}(\eta) \overline{f^\sharp}(\zeta) R_\eta^\beta R_\xi^\alpha \overline{\mathcal{M}} \, d\eta \, d\zeta \, ds \label{cassis2}  \\
& \qquad - \int_1^t \!\! \int  \!\!\!\int  \chi_s   e^{-is\phi } \partial_\eta \overline{f^\sharp}(\eta) \overline{f^\sharp}(\zeta)  R_\xi R_\eta \overline{\mathcal{M}}  \, d\eta \, d\zeta \, ds \label{cassis3}\\
& \qquad + \int_0^t \!\! \int  \!\!\!\int   \chi_s   e^{-is\phi } \overline{f^\sharp}(\eta) \overline{f^\sharp}(\zeta)  (R_\xi E_\eta^* + B_\xi^* R_\xi) \overline{\mathcal{M}}  \, d\eta \, d\zeta \, ds. \label{cassis4}
\end{align}
\end{subequations}
where
\begin{align}
& m_s^1 = (\partial_\xi \chi)_s + \chi_s \sqrt s\partial_\xi \phi \\
& m_s^2 = (\partial_\eta \chi)_s + \chi_s \sqrt s \partial_\eta \phi.
\end{align}
are symbols with a support of size $\sim 1/\sqrt{s}$ and enjoy uniform Coifman-Meyer bounds. 
Let us start with~\eqref{cassis1}:
\begin{align*} \| \eqref{cassis1}\|_{L^2_\xi} & \lesssim \int_1^t \sqrt{s} \|  \int\!\!\! \int \chi_s \overline{[ e^{isH} f ]}^\sharp(\eta) \overline{[e^{isH} f]}^\sharp(\zeta) \overline{\mathcal{M}} \, d\eta \, d\zeta \|_{L^2_\xi} \, ds \\
 & \lesssim \int_1^t \sqrt{s} \left (\| e^{isH} f \|_{L^6} + \| e^{isH} f \|_{L^{6-\delta}}\right) \| e^{isH}P_{< 1/\sqrt{s}}  f \|_{L^{3}} \,ds \nonumber \\
 & \lesssim \int_1^t \sqrt{s} (\frac{1}{s^{7/8}} \|f\|_X) ( \frac{1}{\sqrt{s} } \| P_{< 1/\sqrt{s}}  f \|_{L^{3/2}}) \, ds  \lesssim  \|f\|_X^2
 \end{align*}
where we chose $p> {6}/{5}$ sufficiently small; here we used Bernstein's inequality, as well as appealed to \eqref{LptoweightedL2estimate}. The term~\eqref{cassis2} can be treated similarly to~\eqref{cassis1}, thus we skip it.  For \eqref{cassis3}, we estimate using Theorem \ref{colibri} that
  \begin{align*}
  \| \eqref{cassis3} \|_{L^2} & \lesssim \int_1^t \| \int\!\!\! \int \chi_s   e^{-is\phi } \partial_\eta \overline{f^\sharp}(\eta) \overline{f^\sharp}(\zeta)  R_\xi R_\eta \overline{\mathcal{M}}  \, d\eta \, d\zeta \|_{L^2_\xi}  \, ds  \\
  & \lesssim \int_1^t  \| P_{<1/\sqrt{s}} e^{isH} \Omega x \Omega^* f \|_{L^3} \left( \| e^{isH} f \|_{L^6} + \ \| e^{isH} f \|_{L^{24/5}} \right) \, ds  \\
  & \lesssim \int_1^t   \frac{1}{s^{1/4}} \| e^{isH} \Omega x \Omega^* f \|_{L^2} \frac{1}{s^{7/8}} \| f \|_X \, ds \lesssim \| f \|_X^2.
  \end{align*}
Here we have exploited the fact that $\chi_s$ exhibits the improved Coifman--Meyer structure (in fact, on the frequencies $\leq 1$ where $\chi_s$ is supported, the two notions of Coifman-Meyer symbols and improved Coifman-Meyer symbols coincide).  We omit the details for the final term, \eqref{cassis4}, because it is in fact a little simpler (no derivatives hit the profile).  

\subsubsection{Weighted $L^2$ estimate for $\mathbf{II}$ if $t<1$.} Once again, it suffices to recapitulate the argument from the previous subsection, using the Sobolev embedding theorem where previously we had the dispersive estimates.

\subsubsection{$\dot{H}_\sharp^1$ estimates for ${\mathcal{F}^\sharp}^{-1} \mathbf{II}$.}  Multiplying by $\xi$ and evaluating the $L^2_\xi$ norm, we see that
\begin{align*}
\| {\mathcal{F}^\sharp}^{-1} \mathbf{II} \|_{\dot{H}^1_\sharp}  & \lesssim \int_0^t  \| \int \!\! \int \xi \chi_s e^{-is\phi} \overline{f^\sharp}(\eta) \overline{f^\sharp}(\zeta) \overline{\mathcal{M}} \, d\eta \, d\zeta \|_{L^2_\xi} \, ds \\
& \lesssim \int_0^t \frac{1}{\sqrt{s}} \| \int \!\! \int  m_s e^{-is\phi} \overline{f^\sharp}(\eta) \overline{f^\sharp}(\zeta) \overline{\mathcal{M}} \, d\eta \, d\zeta \|_{L^2_\xi} \, ds.  \end{align*}
where $m_s$ is a Coifman-Meyer symbol supported on a ball of size $1/\sqrt{s}$.

For large time, the desired bound can be obtained with the help of the dispersive estimates, whereas we resort to Sobolev embedding for small time. We only illustrate the latter case: if $t<1$,
\begin{align*}
\| {\mathcal{F}^\sharp}^{-1} \mathbf{II} \|_{\dot{H}^1_\sharp}  & \lesssim \int_0^t \frac{1}{\sqrt{s}} \|  e^{isH} f \|_{L^3} \| e^{isH} f\|_{L^6} \, ds + \int_0^t \frac{1}{\sqrt{s}} \| e^{isH}  f \|_{L^{3-\delta}} \|  e^{isH} f\|_{L^{6}} \, ds \\
 & \lesssim \int_0^t \frac{1}{\sqrt{s}} \|f\|_{H^1_\sharp}^2  \, ds \lesssim \|f\|_X^3.  \end{align*}
This completes the $H_\sharp^1$ estimates for $\mathbf{II}$, and the proof of Proposition~\ref{hummingbird}. \end{proof}

\bigskip

\noindent {\bf Acknowledgement:} The authors are grateful to Fr\'ed\'eric Bernicot for very useful discussions and references on the boundedness of Riesz transforms. They also thank Jonathan Luk, Jalal Shatah, Christopher Sogge, and Daniel Tataru for stimulating discussions.

\bibliographystyle{siam}
\bibliography{vschro}

\begin{thebibliography}{10}

\bibitem{agmon1975spectral}
{\sc S.~Agmon}, {\em Spectral properties of {S}chr{\"o}dinger operators and
  scattering theory}, Ann. Scuola Norm. Sup. Pisa Cl. Sci. (4), 2 (1975),
  pp.~151--218.

\bibitem{Assaad2011}
{\sc J.~Assaad}, {\em Riesz transforms associated to {S}chr{\"o}dinger
  operators with negative potentials}, Publ. Mat., 55 (2011), pp.~123--150.

\bibitem{Auscher2007maximal}
{\sc P.~Auscher and B.~Ben~Ali}, {\em Maximal inequalities and {R}iesz
  transform estimates on {$L^p$} spaces for {S}chr{\"o}dinger operators with
  nonnegative potentials}, Ann. Inst. Fourier (Grenoble), 57 (2007),
  pp.~1975--2013.

\bibitem{BC}
{\sc D.~Bambusi and S.~Cuccagna}, {\em On dispersion of small energy solutions
  of the nonlinear klein gordon equation with a potential}, American journal of
  mathematics, 133 (2011), pp.~1421--1468.

\bibitem{beceanu2011estimate}
{\sc M.~Beceanu}, {\em New estimates for a time-dependent {S}chr{\"o}dinger
  equation}, Duke Math. J., 159 (2011), pp.~417--477.

\bibitem{beceanu2011waveop}
\leavevmode\vrule height 2pt depth -1.6pt width 23pt, {\em Structure of wave
  operators in $\mathbb{R}^3$}, Preprint, http://arxiv.org/abs/1101.0502,
  (2011).

\bibitem{bernicot20121}
{\sc F.~Bernicot}, {\em A t(1)-theorem in relation to a semigroup of operators
  and applications to new paraproducts}, Trans. of Amer. Math. Soc,  (2012).

\bibitem{blue2008decay}
{\sc P.~Blue}, {\em Decay of the maxwell field on the schwarzschild manifold},
  Journal of Hyperbolic Differential Equations, 5 (2008), pp.~807--856.

\bibitem{coifman1978dela}
{\sc R.~Coifman and Y.~Meyer}, {\em Au dela des op{\'e}rateurs
  pseudo-diff{\'e}rentiels}, Soci{\'e}t{\'e} math{\'e}matique de France, 1978.

\bibitem{cuccagna2012decay}
{\sc S.~Cuccagna, V.~Georgiev, and N.~Visciglia}, {\em Decay and scattering of
  small solutions of pure power {NLS} in $\mathbb{R}$ with $ p> 3$ and with a
  potential}, arXiv preprint arXiv:1209.5863,  (2012).

\bibitem{dafermos2010new}
{\sc M.~Dafermos and I.~Rodnianski}, {\em A new physical-space approach to
  decay for the wave equation with applications to black hole spacetimes}, in
  XVIth International Congress on Mathematical Physics, World Sci. Publ.,
  Hackensack, NJ, 2010, pp.~421--432.

\bibitem{Fefferman1971}
{\sc C.~Fefferman and E.~M. Stein}, {\em Some maximal inequalities}, Amer. J.
  Math., 93 (1971), pp.~107--115.

\bibitem{finco2006p}
{\sc D.~Finco and K.~Yajima}, {\em The {$L^p$} boundedness of wave operators
  for {S}chr{\"o}dinger operators with threshold singularities. {II.} even
  dimensional case}, Journal of Mathematical Sciences, The University of Tokyo,
  13 (2006), pp.~277--346.

\bibitem{germain2011resonance}
{\sc P.~Germain}, {\em Space-time resonances}, Preprint,
  http://arxiv.org/abs/1102.1695,  (2011).

\bibitem{germain2012non}
{\sc P.~Germain, N.~Masmoudi, and B.~Pausader}, {\em Non-neutral global
  solutions for the electron {E}uler-{P}oisson system in 3d}, arXiv preprint
  arXiv:1204.1536,  (2012).

\bibitem{germain2009global}
{\sc P.~Germain, N.~Masmoudi, and J.~Shatah}, {\em Global solutions for 3d
  quadratic {S}chr{\"o}dinger equations}, International Mathematics Research
  Notices, 2009 (2009), pp.~414--432.

\bibitem{germain2012global}
\leavevmode\vrule height 2pt depth -1.6pt width 23pt, {\em Global solutions for
  the gravity water waves equation in dimension 3}, Ann. Math, 175 (2012),
  pp.~691--754.

\bibitem{guillarmou2008resolvent}
{\sc C.~Guillarmou and A.~Hassell}, {\em Resolvent at low energy and {R}iesz
  transform for {S}chr{\"o}dinger operators on asymptotically conic manifolds.
  {I}.}, Mathematische Annalen, 341 (2008), pp.~859--896.

\bibitem{Guo1998}
{\sc Y.~Guo}, {\em Smooth irrotational flows in the large to the
  {E}uler-{P}oisson system in {$\bold R^{3+1}$}}, Comm. Math. Phys., 195
  (1998), pp.~249--265.

\bibitem{guo2011global}
{\sc Y.~Guo and B.~Pausader}, {\em Global smooth ion dynamics in the
  {E}uler--{P}oisson system}, Communications in Mathematical Physics, 303
  (2011), pp.~89--125.

\bibitem{gustafson2004asymptotic}
{\sc S.~Gustafson, K.~Nakanishi, and T.-P. Tsai}, {\em Asymptotic stability and
  completeness in the energy space for nonlinear schr{\"o}dinger equations with
  small solitary waves}, International Mathematics Research Notices, 2004
  (2004), pp.~3559--3584.

\bibitem{gustafson2009scattering}
{\sc S.~Gustafson, K.~Nakanishi, and T.-P. Tsai}, {\em Scattering theory for
  the gross--pitaevskii equation in three dimensions}, Communications in
  Contemporary Mathematics, 11 (2009), pp.~657--707.

\bibitem{hani2012global}
{\sc Z.~Hani}, {\em Global well-posedness of the cubic nonlinear
  {S}chr{\"o}dinger equation on closed manifolds}, Communications in Partial
  Differential Equations, 37 (2012), pp.~1186--1236.

\bibitem{hormander2007analysis}
{\sc L.~H{\"o}rmander}, {\em The analysis of linear partial differential
  operators I: Distribution theory and {F}ourier Analysis}, vol.~1, Springer,
  2007.

\bibitem{ikebe1960eigen}
{\sc T.~Ikebe}, {\em Eigenfunction expansions associated with the
  {S}chroedinger operators and their applications to scattering theory}, Arch.
  Rational Mech. Anal., 5 (1960), pp.~1--34 (1960).

\bibitem{ionescu2006agmon}
{\sc A.~D. Ionescu and W.~Schlag}, {\em Agmon-{K}ato-{K}uroda theorems for a
  large class of perturbations}, Duke Math. J., 131 (2006), pp.~397--440.

\bibitem{jensen2002remark}
{\sc A.~Jensen and K.~Yajima}, {\em A remark on {$L^p$}-boundedness of wave
  operators for two dimensional {S}chr{\"o}dinger operators}, Communications in
  mathematical physics, 225 (2002), pp.~633--637.

\bibitem{john1979blow}
{\sc F.~John}, {\em Blow-up of solutions of nonlinear wave equations in three
  space dimensions}, Manuscripta Mathematica, 28 (1979), pp.~235--268.

\bibitem{journe1991decay}
{\sc J.-L. Journ{{\'e}}, A.~Soffer, and C.~D. Sogge}, {\em Decay estimates for
  {S}chr{\"o}dinger operators}, Comm. Pure Appl. Math., 44 (1991),
  pp.~573--604.

\bibitem{keel2002global}
{\sc M.~Keel, H.~Smith, and C.~Sogge}, {\em Global existence for a quasilinear
  wave equation outside of star-shaped domains}, Journal of Functional
  Analysis, 189 (2002), pp.~155--226.

\bibitem{Klainerman}
{\sc S.~Klainerman}, {\em The null condition and global existence to nonlinear
  wave equations}, Nonlinear systems of partial differential equations in
  applied mathematics, Part 1 (Santa Fe, N.M., 1984), Lectures in Appl. Math.,
  23 (1986), pp.~293--326.

\bibitem{Kosecki1992}
{\sc R.~Kosecki}, {\em The unit condition and global existence for a class of
  nonlinear {K}lein-{G}ordon equations}, J. Differential Equations, 100 (1992),
  pp.~257--268.

\bibitem{laillet2011thesis}
{\sc N.~Laillet}, {\em R{\'e}sonances en espace et en temps pour l'{\'e}quation
  de {S}chr{\"o}dinger dans $\mathbb{R}^3$ avec une nonlin{\'e}arit{\'e}},
  Master's thesis, 2011.

\bibitem{mckean1991nonlinear}
{\sc H.~McKean and J.~Shatah}, {\em The nonlinear {S}chr{\"o}dinger equation
  and the nonlinear heat equation reduction to linear form}, Communications on
  Pure and Applied Mathematics, 44 (1991), pp.~1067--1080.

\bibitem{metcalfe2005global}
{\sc J.~Metcalfe, M.~Nakamura, and C.~Sogge}, {\em Global existence of
  solutions to multiple speed systems of quasilinear wave equations in exterior
  domains}, 17 (2005), p.~133.

\bibitem{metcalfe2007global}
{\sc J.~Metcalfe and C.~Sogge}, {\em Global existence of null-form wave
  equations in exterior domains}, Mathematische Zeitschrift, 256 (2007),
  pp.~521--549.

\bibitem{Ozawa1995}
{\sc T.~Ozawa, K.~Tsutaya, and Y.~Tsutsumi}, {\em Normal form and global
  solutions for the {K}lein-{G}ordon-{Z}akharov equations}, Ann. Inst. H.
  Poincar{\'e} Anal. Non Lin{\'e}aire, 12 (1995), pp.~459--503.

\bibitem{rodnianski2004decay}
{\sc I.~Rodnianski and W.~Schlag}, {\em Time decay for solutions of
  {S}chr{\"o}dinger equations with rough and time-dependent potentials},
  Invent. Math., 155 (2004), pp.~451--513.

\bibitem{schlag2007dispersive}
{\sc W.~Schlag}, {\em Dispersive estimates for Schr{\"o}dinger operators: A
  survey}, vol.~163, Mathematical Aspects of Nonlinear Dispersive Equations,
  Ann. of Math. Stud, 2007.

\bibitem{shatah1985normal}
{\sc J.~Shatah}, {\em Normal forms and quadratic nonlinear klein-gordon
  equations}, Communications on Pure and Applied Mathematics, 38 (1985),
  pp.~685--696.

\bibitem{Shen1995}
{\sc Z.~W. Shen}, {\em {$L^p$} estimates for {S}chr{\"o}dinger operators with
  certain potentials}, Ann. Inst. Fourier (Grenoble), 45 (1995), pp.~513--546.

\bibitem{SW}
{\sc A.~Soffer and M.~I. Weinstein}, {\em Resonances, radiation damping and
  instabilitym in hamiltonian nonlinear wave equations}, Inventiones
  mathematicae, 136 (1999), pp.~9--74.

\bibitem{stein1993harmonicbook}
{\sc E.~M. Stein}, {\em Harmonic analysis: real-variable methods,
  orthogonality, and oscillatory integrals}, vol.~43 of Princeton Mathematical
  Series, Princeton University Press, Princeton, NJ, 1993.
\newblock With the assistance of Timothy S. Murphy, Monographs in Harmonic
  Analysis, III.

\bibitem{strichartz1983analysis}
{\sc R.~Strichartz}, {\em Analysis of the laplacian on the complete riemannian
  manifold}, J. Funct. Anal, 52 (1983), pp.~48--79.

\bibitem{tsai2001asymptotic}
{\sc T.-P. Tsai and H.-T. Yau}, {\em Asymptotic dynamics of nonlinear
  schr{\"o}dinger equations: Resonance-dominated and dispersion-dominated
  solutions}, Communications on pure and applied mathematics, 55 (2001),
  pp.~153--216.

\bibitem{wu2011global}
{\sc S.~Wu}, {\em Global wellposedness of the 3-d full water wave problem},
  Inventiones mathematicae, 184 (2011), pp.~125--220.

\bibitem{yajima1993waveop}
{\sc K.~Yajima}, {\em The {$W^{k,p}$}-continuity of wave operators for
  {S}chr{\"o}dinger operators}, Proc. Japan Acad. Ser. A Math. Sci., 69 (1993),
  pp.~94--98.

\bibitem{yajima1995Wkp}
\leavevmode\vrule height 2pt depth -1.6pt width 23pt, {\em The
  {$W^{k,p}$}-continuity of wave operators for {S}chr{\"o}dinger operators}, J.
  Math. Soc. Japan, 47 (1995), pp.~551--581.

\bibitem{yajima1995Wkp3}
\leavevmode\vrule height 2pt depth -1.6pt width 23pt, {\em The
  {$W^{k,p}$}-continuity of wave operators for {S}chr{\"o}dinger operators.
  {III}. {E}ven-dimensional cases {$m\geq 4$}}, J. Math. Sci. Univ. Tokyo, 2
  (1995), pp.~311--346.

\bibitem{yajima2006Lp}
\leavevmode\vrule height 2pt depth -1.6pt width 23pt, {\em The {$L^p$}
  boundedness of wave operators for {S}chr{\"o}dinger operators with threshold
  singularities. {I}. {T}he odd dimensional case}, J. Math. Sci. Univ. Tokyo,
  13 (2006), pp.~43--93.

\bibitem{yang2012global}
{\sc S.~Yang}, {\em Global stability of solutions to nonlinear wave equations},
  arXiv preprint arXiv:1205.4216,  (2012).

\end{thebibliography}

\end{document}